\documentclass[11pt]{article}
\usepackage[left=3cm,right=3cm,top=2.2cm,bottom=2.2cm]{geometry}
\usepackage{amsmath,amssymb}
\usepackage{amsthm}
\usepackage{mathtools,bm}
\usepackage[shortlabels]{enumitem}
\usepackage{authblk}
\usepackage{titlefoot}
\usepackage[title]{appendix}

\usepackage[style = alphabetic, backend = biber]{biblatex}
\theoremstyle{plain}
\newtheorem{theorem}{Theorem}[section]

\newtheorem{corollary}[theorem]{Corollary}

\newtheorem{lemma}[theorem]{Lemma}
\newtheorem{proposition}[theorem]{Proposition}

\theoremstyle{definition}
\newtheorem{definition}[theorem]{Definition}

\newtheorem{remark}[theorem]{Remark}

\newenvironment{proofThm}[1]
{
	\par\vspace{\baselineskip}\noindent
	\textit{Proof of Theorem} #1.
	\noindent
}
{
	\qed\ignorespacesafterend
}

\usepackage{graphicx}
\makeatletter
\newcommand*\bigcdot{{\mathpalette\bigcdot@{.5}}}
\newcommand*\bigcdot@[2]{\mathbin{\vcenter{\hbox{\scalebox{#2}{$\m@th#1\bullet$}}}}}
\makeatother

\newcommand{\vertiii}[1]{{\left\vert\kern-0.25ex\left\vert\kern-0.25ex\left\vert #1 
    \right\vert\kern-0.25ex\right\vert\kern-0.25ex\right\vert}}

\begin{document}
	
\title{The Lions Derivative in Infinite Dimensions - Application to Higher Order Expansion of Mean-Field SPDEs}
\date{\today}
	
\author[1]{Wilhelm Stannat}
\author[2]{Alexander Vogler}
	
\affil[1]{Technische Universität Berlin, Berlin, Germany}
\affil[2]{Berlin, Germany}
	
	
\maketitle
	
\unmarkedfntext{\textit{Mathematics Subject Classification (2020) --- 46G05, 46G10, 60H15, 49N80} }
	
\unmarkedfntext{\textit{Keywords and phrases --- Lions derivative, Stochastic partial differential  
equations, Mean-Field Stochastic partial differential equations, SPDE, Vector measures, Taylor expansions, 
Functional Itô calculus}}
	
\unmarkedfntext{\textit{Mail}: \textbullet\, stannat@math.tu-berlin.de\, \textbullet\,
alexander\_vogler@hotmail.de}


\begin{abstract}
In this paper we present a new interpretation of the Lions derivative as the Radon-Nikodym derivative of a 
vector measure, which provides a canonical extension of the Lions derivative for functions taking values in 
infinite dimensional Banach spaces. This is of particular relevance for the analysis of Hilbert space valued 
Mean-Field equations. As an illustration we derive a mild Ito-formula for Mean-Field stochastic partial 
differential equations (SPDEs), which provides the basis for a higher order Taylor expansion and higher order 
numerical schemes.
\end{abstract}
	
\tableofcontents
	
\section{Introduction and Main Result}
	
The notion of $L$-differentiability for functions on spaces of probability measures is of fundamental 
importance in McKean-Vlasov stochastic differential equations with applications to mean-field games and 
optimal control theory. In the finite dimensional setting the notion of L-differentiability is based on the 
lifting of a function $f:\mathcal{P}_2(\mathbb{R}^d)\rightarrow \mathbb{R}$, defined on the space of square 
integrable probability measures $\mathcal{P}_2(\mathbb{R}^d)$, to the space of random variables
\begin{align*}
\hat{f}:L^2(\Omega, \mathcal{F}, \mathbb{P};\mathbb{R}^d)\rightarrow \mathbb{R},X\mapsto f(\mathcal{L}(X)),
\end{align*}
where $\Omega$ is a Polish space, $\mathcal{F}$ its Borel-sigma-algebra and $\mathbb{P}$ is atomless. 
In particular, a 
function $f$ is called L-differentiable in $\mu_0\in \mathcal{P}_2(\mathbb{R}^d)$, if there exists a random 
variable $X_0\in L^2(\Omega, \mathcal{F}, \mathbb{P};\mathbb{R}^d)$ with $\mathcal{L}(X_0)=\mu_0$, such that the lift $\hat{f}$ is 
Fr\'echet differentiable in $X_0$. Due to the nature of the L-derivative, it is well suited to control 
infinitesimal perturbations induced by variations on the space of random variables. These type of variations 
naturally appear in the proof of the maximum principle and the dynamic programming principle for mean-field 
control problems and mean-field game problems. In order to show that this notion of differentiability is 
intrinsic and to provide some structure of the derivative, the authors in \cite{CD18} derive a factorized 
gradient representation for the Fr\'echet derivative of the lifted function of the type
\begin{align*}
D\hat{f}(X_0)Y=\mathbb{E}\left[\langle g(X_0),Y\rangle\right],
\end{align*}
where $g\in L^2(\mathbb{R}^d,\mathcal{B} (\mathbb{R}^d), \mu_0;\mathbb{R}^d)$ depends only on $\mu_0=\mathcal{L}(X_0)$. By setting 
$\partial_{\mu}f(\mu_0)=g$ this gives a meaning to the $L$-derivative which is completely independent of the 
lifting $\hat{f}$ (i.e. independent of the probability space $(\Omega, \mathcal{F},\mathbb{P})$ and the random 
variable $X_0$). In particular, if one considers another probability space $(\tilde{\Omega},
\tilde{\mathcal{F}},\tilde{\mathbb{P}})$ such that there are random variables $\tilde{X}_0,\tilde{Y}$ with 
$\mathcal{L}(\tilde{X}_0,\tilde{Y})=\mathcal{L}(X_0,Y)$, then 
\begin{align*}
D\hat{f}(X_0)Y&=\mathbb{E}\left[\langle \partial_{\mu}f(\mu_0)(X_0),Y\rangle\right] \\
& = \int_{\mathbb{R}^d\times \mathbb{R}^d}\langle \partial_{\mu}f(\mu_0)(x),y\rangle 
   \mathcal{L}(X_0,Y)(d(x,y)) \\
& = \int_{\mathbb{R}^d\times \mathbb{R}^d}\langle \partial_{\mu}f(\mu_0)(x),y\rangle  
   \mathcal{L}(\tilde{X}_0,\tilde{Y})(d(x,y)) \\
& = \tilde{\mathbb{E}}\left[\langle \partial_{\mu}f(\mu_0)(\tilde{X}_0),\tilde{Y}\rangle\right] \\
& = D\tilde{f}(\tilde{X}_0)\tilde{Y}.
\end{align*}

\medskip 
\noindent
The factorized gradient representation crucially relies on the identification 
$$ 
L(L^2(\Omega,\mathcal{F}, \mathbb{P}; \mathbb{R}^d), \mathbb{R})\cong L^2(\Omega, \mathcal{F},  
\mathbb{P};\mathbb{R}^d), 
$$ 
by the Riesz representation theorem. Identifying $D\hat{f}(X_0)$ as an element of $L^2(\Omega, \mathbb{R}^d)$ 
it only remains to show that $D\hat{f}(X_0)$ is $\sigma(X_0)$-measurable in order to obtain the desired 
factorization. 

\medskip
\noindent
It is natural to ask if the notion of L-differentiability can be extended to the infinite dimensional setting, 
where $f:\mathcal{P}_2(H)\rightarrow U$, for a separable Hilbert space $H$ and some Banach space $U$. This 
question is of particular importance in the context of McKean-Vlasov SPDE optimal control problems, where 
infinitesimal perturbations of Hilbert space valued random variables naturally appear in the maximum principle 
and dynamic programming principle. As in the finite dimensional setting one can define $L$-differentiability of 
a function $f:\mathcal{P}_2(H)\rightarrow U$ in terms of the Fr\'echet differentiability of its lift. 
To do so we fix an atomless complete probability space $(\Omega, \mathcal{F},\mathbb{P})$, where $\Omega$ is 
Polish and $\mathcal{F}$ is its Borel $\sigma$-field. On such a probability space we can construct for any 
probability distribution $\mu$ on some metric space $E$ a random variable on $(\Omega,\mathcal{F},\mathbb{P})$ 
with distribution $\mu$. This allows us to define the lift of a map on the space of probability measures 
on $E$. 

\begin{definition}
\label{Ldiff1}
Given some Hilbert space $H$ and a Banach space $U$, we call a map $f:\mathcal{P}_2(H)\rightarrow U$ 
L-differentiable in $\mu_0\in \mathcal{P}_2(H)$, if there exists a random variable $X_0\in L^2(\Omega, 
\mathcal{F}, \mathbb{P};$ $H)$ with $\mathcal{L}(X_0)=\mu_0$, such that the lift
\begin{align*}
\hat{f}:L^2(\Omega, \mathcal{F}, \mathbb{P};H)\rightarrow U,X\mapsto f(\mathcal{L}(X))
\end{align*}
is Fr\'echet differentiable in $X_0$.
\end{definition}

\noindent 
Although we can define $L$-differentiability this way, it is not clear whether the $L$-derivative depends 
on the choice of the lifting $\hat{f}$. In order to give a meaning to the $L$-derivative independent of the 
choice for the lifting, we need to find a factorized representation of $D\hat{f}(X_0)$ similar to the finite 
dimensional setting. In \cite{CGKPR22} the authors construct a factorization of $D\hat{f}(X_0) 
\in L(L^2(\Omega, \mathcal{F}, \mathbb{P};H),U)$ in the case where $H=D([0,T], \tilde{H})$ is the space of 
$\tilde{H}$-valued càdlàg paths and $U=\mathbb{R}$, in order to derive the dynamic programming principle 
for path dependent McKean-Vlasov SDEs in the Hilbert space $\tilde{H}$. They explicitly construct the 
factorization of $D\hat{f}(X_0)\in L(L^2(\Omega, \mathcal{F}, \mathbb{P};H), \mathbb{R}) \cong 
L^2(\Omega, \mathcal{F}, \mathbb{P};H)$ in the case where $\mu_0$ is of the 
simple type $\mu_0=\sum_{i=1}^N p_i\delta_{x_i}$ and then proceed by a limit argument using the continuity of 
$D\hat{f}$ for the general case. This explicit construction is rather technical and lacks of interpretation. 
Furthermore it still crucially relies on the identification $D\hat{f}(X_0)\in L(L^2(\Omega,\mathcal{F}, 
\mathbb{P};H),\mathbb{R})\cong L^2(\Omega, \mathcal{F}, \mathbb{P};H)$.

\medskip 
\noindent 	
When $U$ is not finite dimensional anymore the situation becomes much more difficult, since the Fr\'echet 
derivative of the lifting is now an element of $L(L^2(\Omega, \mathcal{F}, \mathbb{P};H),U)$, which makes it 
difficult to identify $D\hat{f}(X_0)$ with a random variable. If $U$ is a separable Hilbert space, given an 
orthonormal basis (ONB) $(\tilde{\phi}_i)_{i\in \mathbb{N}}$ of $U$, the first naive approach would be to 
directly derive the factorized gradient representation by considering 
\begin{align*}
\hat{f}_i:L^2(\Omega, \mathcal{F}, \mathbb{P};H) \rightarrow \mathbb{R},  
\hat{f}_i(X)=\langle \hat{f}(X),\tilde{\phi}_i\rangle_U
\end{align*}
and define 
\begin{equation}\label{sumDiv}
g(u)h := \sum_{i=1}^{\infty} \langle g_i(u),h\rangle_H\tilde{\phi}_i,
\end{equation}
where $g_i\in L^2(H,\mathcal{B} (H), \mu_0;H)$ is given by the factorized gradient representation of $D\hat{f}_i$, which 
satisfies 
\begin{align*}
D\hat{f}_i(X_0)Y=\mathbb{E}\left[\langle g_i(X_0),Y\rangle_H\right],
\end{align*}
for any $Y\in L^2(\Omega, \mathcal{F}, \mathbb{P};H)$. 

\medskip 
\noindent 	
Unfortunately it is not clear if $g$ is well-defined, since we do not have any suitable a-priori estimates 
on $g_i$. In fact, in Section \ref{Examples1}, we will 
present an explicit example where the series in \eqref{sumDiv} will not converge in any suitable manner, 
although the map $f$ is $L$-differentiable in the sense of Definition \ref{Ldiff1}. In order to circumvent 
these difficulties, we will consider a completely different approach which will also give us a new 
interpretation of the $L$-derivative as the Radon-Nikodym derivative of a vector measure, which seems to be new 
even in the finite dimensional setting. The following simple result based on \cite[§13.3 Theorem~1]{Din66} 
allows us to identify $D\hat{f}(X_0)$ with a vector measure.

\begin{lemma}
Let $H$ be a Hilbert space, $U$ be a Banach space and $\hat{f} : L^2 (\Omega , \mathcal{F}, \mathbb{P}; H) 
\rightarrow U$ be Fr\'echet differentiable in $X_0\in L^2(\Omega, \mathcal{F}, \mathbb{P};H)$, then 
$$
m_{D\hat{f} (X_0)} (A) u:= D\hat{f} (X_0) (u1_A) , A\in\mathcal F , u\in H \, , 
$$
defines a vector measure, see Definition \ref{DefVecMeas},
$$ 
m_{D\hat{f} (X_0)} : \mathcal F \rightarrow L(L^2 (\Omega , \mathcal F, \mathbb{P} ; H), U). 
$$
\end{lemma}

\noindent 
In fact there is a one-to-one correspondence between vector measures and linear operators, see Theorem 
\ref{meas1}. The existing results on vector measures then provide us with the right tools to construct 
a factorized gradient of $D \hat{f}(X_0)$ as the Radon-Nikodym derivative of the image measure  
$m_{D\hat{f} (X_0)}\circ X_0^{-1}$ with respect to $\mathcal{L}(X_0)=\mu_0$. One of the key observations 
will be that $m_{\mu_0}=m_{D\hat{f} (X_0)}\circ X_0^{-1}$ is independent of the choice of the Lions lift, 
i.e. independent of $(\Omega,\mathcal{F},\mathbb{P})$ and $X_0\in L^2(\Omega, \mathcal{F}, \mathbb{P};H)$, 
as long as $\mathcal{L}(X_0)=\mu_0$. In particular, we can give meaning to the L-derivative for $U$-valued 
functions $f$, completely independent of the lifting $\hat{f}$. Our main additional assumption which is needed in order 
to provide the structure of the L-derivative will be the continuity of $D\hat{f}$ in the 2-variation norm, 
$\vertiii{\cdot}_{2,\mathbb{P}}$, which is defined below in Definition \ref{2var}. In the situation $U=\mathbb{R}$ 
the 2-variation norm coincides with the usual operator norm, see Lemma \ref{variation1d}, and our main 
assumption reduces to the usual continuity assumption on $D\hat{f}$ also imposed in \cite{CD18} or 
\cite{CGKPR22}.
	
\begin{definition}
\label{2var}
For some Hilbert space $H$ and a Banach space $U$ we define the Banach space $(\Lambda_2^{\mathbb{P}}(H,U), 
\vertiii{\cdot}_{2,\mathbb{P}})$ of continuous linear operators with bounded 2-variation by
\begin{align*}
\Lambda_2^{\mathbb{P}}(H,U):=\{L\in L(L^2(\Omega, \mathcal{F}, \mathbb{P}; H),U):\vertiii{L}_{2,\mathbb{P}}
 <\infty\},
\end{align*}
where $\Lambda_2^{\mathbb{P}}(H,U)$ is equipped with the 2-variation norm 
\begin{align*}
\vertiii{L}_{2,\mathbb{P}} 
& :=\sup\{\sum_{i=1}^{n}\|L(\mathbf{1}_{A_i}x_i)\|_U : 
Y =\sum_{i=1}^{n}\mathbf{1}_{A_i} x_i\in \mathcal{E}_H,A_i\text{ disjoint}, 
\mathbb{E}\left[\|Y\|_H^2\right]\le 1\}.
\end{align*} 
Here $\mathcal{E}_H$ denotes the set of simple functions in $H$. We say that a map $\hat{f}:L^2(\Omega, 
\mathcal{F}, \mathbb{P};H)\rightarrow U$ is $\Lambda$-continuously Fr\'echet differentiable if $\hat{f}$ 
is Fr\'echet differentiable with $D\hat{f}(X)\in \Lambda^{\mathbb{P}}_2(H,U)$ for all $X\in L^2(\Omega, 
\mathcal{F}, \mathbb{P};H)$, and $D\hat{f}:L^2(\Omega, \mathcal{F}, \mathbb{P};H)\rightarrow 
\Lambda_2^{\mathbb P} (H,U)$ is continuous.
\end{definition}

\noindent 
The proof that $(\Lambda_2^{\mathbb{P}}(H,U),\vertiii{\cdot}_{2,\mathbb{P}})$ is indeed a Banach space is 
postponed to Section \ref{LinOp}. Our main theorem is now formulated as follows:
	
\begin{theorem} 
\label{g}
Let $H$ be a separable real Hilbert space, $U$ be a Banach space and $f:\mathcal{P}_2(H)\rightarrow U$, 
such that the lift
\begin{align*}
\hat{f}:L^2(\Omega, \mathcal{F}, \mathbb{P};H)\rightarrow U, X\mapsto f(\mathcal{L}(X)),
\end{align*}
is $\Lambda$-continuously Fr\'echet differentiable. Then for any $\mu\in \mathcal{P}_2(H)$ there exists 
a measurable map $\frac{dm_{\mu}}{d\mu}:H\rightarrow L(H,U)$, such that $\frac{dm_{\mu}}{d\mu}(X) 
\in L^2(\Omega,\mathcal{F}, \mathbb{P}; L(H,U))$ for all $X\in L^2(\Omega, \mathcal{F}, \mathbb{P};H)$ with 
$\mathcal{L}(X)=\mu$, and for all $Y\in L^2(\Omega, \mathcal{F}, \mathbb{P};H)$ it holds 
$\frac{dm_{\mu}}{d\mu}(X)Y\in L^1(\Omega, \mathcal{F}, \mathbb{P};U)$ with
\begin{equation}
\label{RN}
\begin{aligned}
D\hat{f}(X)Y 
& = \int Ydm_{D\hat{f}(X)} \\
& = \int \mathbb{E}\left[Y|X\right]dm_{D\hat{f} (X)}\\
& = \mathbb{E}\left[\frac{dm_{\mu}}{d\mu}(X)Y\right].
\end{aligned}
\end{equation}
Furthermore we have for all $X\in L^2(\Omega, \mathcal{F}, \mathbb{P};H)$ with $\mathcal{L}(X)=\mu$
\begin{align*}
\mathbb{E}\left[\left\|\frac{dm_{\mu}}{d\mu}(X)\right\|_{op}^2\right]^{\frac{1}{2}} 
= \vertiii{D\hat{f}(X)}_{2,\mathbb{P}}.
\end{align*}
The map $\frac{dm_{\mu}}{d\mu}$ is a $\mu$-version of the Radon-Nikodym derivative of  
$m_{\mu}=m_{D\hat{f} (X)}\circ X^{-1}$ with respect to $\mu=\mathcal{L}(X)$, since by \eqref{RN} we have 
for all $A\in \mathcal{B}(H)$ and $u\in H$   
\begin{align*}
m_{\mu}(A)u = m_{D\hat{f} (X)}(X\in A)u =\int \mathbf{1}_A(x)\frac{dm_{\mu}}{d\mu}(x)ud\mu(x).
\end{align*}
We denote the equivalence class of $\frac{dm_{\mu}}{d\mu}$ in $L^2(H,\mathcal{B} (H), \mu;L(H,U))$, which 
is uniquely defined, by $\partial_{\mu}f(\mu)(\cdot)$.
\end{theorem}

\noindent
The assumption $D\hat{f}(X)\in \Lambda_2^{\mathbb{P}}(H,U)$ for all $X\in L^2 L^2(\Omega, \mathcal{F},
\mathbb{P};H)$ is indeed necessary which can be seen as follows. Assume that $\hat{f}:L^2(\Omega, 
\mathcal{F}, \mathbb{P};H)\rightarrow U$ is Fr\'echet differentiable and there exists some measurable map 
$\partial_{\mu}f(\mu_0)\in L^2(H,\mathcal{B}(H), \mu_0;L(H,U))$ with 
\begin{align*}
D\hat{f}(X_0)Y=\mathbb{E}\left[\partial_{\mu}f(\mu_0)(X_0)Y\right]
\end{align*}
for all $Y\in L^2(\Omega, \mathcal{F}, \mathbb{P};H)$. Then by Lemma \ref{variationg} the vector measure
\begin{align*}
m_{D\hat{f}(X_0)}(A):=D\hat{f}(X_0)(\mathbf{1}_A)
\end{align*} 
has finite variation given by
\begin{align*}
|m_{D\hat{f}(X_0)}| & :=\sup\{\sum_{i=1}^{n}\|m_{D\hat{f}(X_0)}(A_i)\|_{op}:n\in \mathbb{N},(A_i)_i  
\subset\mathcal{F} \text{ disjoint with }A_i\subset \Omega\} \\ 
& \, \, = \mathbb{E}\left[\|\partial_{\mu}f(\mu_0)(X_0)\|_{op}\right].
\end{align*}
Therefore by Lemma \ref{RN-iso} and Corollary \ref{cor1} we have
\begin{align*}
\vertiii{D\hat{f}(X_0)}_{2,\mathbb{P}}^2 
= \mathbb{E}\left[\|\partial_{\mu}f(\mu_0)(X_0)\|_{op}^2\right]<\infty.
\end{align*}
In Section \ref{Examples1} we will provide an example of a map $f:\mathcal{P}_2(H)\rightarrow U$ which is 
$L$-differentiable, i.e. $f$ has Fr\'echet differentiable lift, but does not satisfy $D\hat{f}(X_0) 
\in \Lambda_2^{\mathbb{P}}(H,U)$ for general $X_0\in L^2(\Omega,\mathcal{F}, \mathbb{P};H)$. This shows in 
particular that we have found a somehow new necessary condition which does not appear if $U$ is finite 
dimensional, but is non-trivial when $U$ is not finite dimensional. \\
	
\noindent 	
The rest of the paper will be structured as follows. In Section \ref{pre} we will recall some important 
concepts and results in the theory of vector measures. Section \ref{generalL} is concerned with the proof 
of our main theorem. As an application of our main result we will derive a mild Ito-formula for Mean-Field 
SPDEs in Section \ref{Ito} which provides the basis for our higher order Taylor expansion. Section 
\ref{outlook} provides some outlook for the possible applications of our theory in the context of mean-field 
SPDE control problems.

\section{Preliminaries} 
\label{pre}
	
Throughout the rest of the manuscript we fix a Polish space $\Omega$, denote with $\mathcal{F}$ its Borel 
$\sigma$-algebra, and let $\mathbb{P}$ be an atomless probability measure on $(\Omega , \mathcal{F})$.  
Moreover, we fix a separable real Hilbert space $H$ with ONB $(\varphi_i)_{i\in \mathbb{N}}$ and a Banach 
space $U$.
	
\subsection{Notation}
	
For $p\geq 1$ and a Banach space $(X,\|\cdot\|_X)$ we define the Bochner space $L^p (\Omega ,\mathcal{F}, 
\mathbb{P};X)$ as the quotient (with equality $\mathbb{P}$-almost everywhere) of the space 
\begin{align*}
\mathcal{L}^p(\Omega, \mathcal{F}, \mathbb{P}; X):=\{Y:\Omega\rightarrow X|Y\text{ is Bochner measurable, } 
\mathbb{E}\left[\|Y\|_X^p\right]<\infty\}.
\end{align*}
We also use the short hand notations $L^p (\Omega ; X, \mathbb{P})$ and $L^p (\Omega ; X)$, instead of 
$L^p (\Omega ,\mathcal{F}, \mathbb{P};X)$. For a random variable $Y:\Omega\rightarrow X$ we denote its 
distribution on $(X,\mathcal{B}(X))$ by 
\begin{align*}
\mathcal{L}(X) = \mathbb{P}\circ X^{-1},
\end{align*}
where $\mathcal{B}(X)$ denotes the Borel-sigma field on $X$. Furthermore by $\mathcal{P}_2(X)$ we denote 
the set of probability measures $\mu$ on $(X, \mathcal{B}(X))$ with 
\begin{align*}
\int \|x-x_0\|_X^2\mu(dx)<\infty,
\end{align*} 
for one, hence all, $x_0\in X$. On $\mathcal{P}_2(X)$ we consider the 2-Waserstein distance
\begin{align*}
W_2(\mu,\nu):=\inf_{\pi\in \Pi(\mu,\nu)}\left(\int_{X\times X}\|x-y\|_X^2\pi(d(x,y))\right)^{\frac{1}{2}},
\end{align*} 
where $\Pi(\mu,\nu)$ denotes the set of couplings w.r.t. $\mu,\nu\in \mathcal{P}_2(X)$. The set of 
continuous linear operators from $H$ to $U$ will be denoted by $L(H,U)$, which is equipped with the 
usual operator norm 
\begin{align*}
\| L\|_{op}:=\sup_{u\in H,\|u\|_H\le 1}\|Lu\|_U.
\end{align*}
Furthermore we denote the dual space of $U$ by $U'$ and by $\langle \cdot,\cdot\rangle$ the corresponding 
dual pairing. The space of Hilbert-Schmidt operators on $H$ is defined by
\begin{align*}
L_2(H,H):=\{L\in L(H,H):\|L\|^2_{L_2(H,H)}:=\sum_{i=1}^{\infty}\|L\phi_i\|_H^2<\infty\}.
\end{align*}
	
\subsection{Vector Measures}
	
In this subsection we will give a short introduction in the theory of vector measures and corresponding 
integration theory, based on \cite{Din66}. In order to improve readability, we will also give proofs of some 
results. The reader may skip this section if he is already familiar with this theory. In the following 
let $(X,\|\cdot\|_X)$ be a Banach space.

\begin{definition}\label{DefVecMeas}
A set function $m:\mathcal{F}\rightarrow X$ is called a measure, if 
\begin{align*}
m(\emptyset) 
& = 0,  \\
m(\bigcup_{i=1}^{\infty}A_i)&=\sum_{i=1}^{\infty}m(A_i),
\end{align*}
for any mutually disjoint sets $(A_i)_{i\in \mathbb{N}}\subseteq \mathcal{F}$.
\end{definition}
	
\begin{definition}
For a set function $m:\mathcal{F}\rightarrow X$ we define its variation $\bar{m}$ by
\begin{align*}
\overline{m}(A) := \sup\{\sum_{i=1}^{n}\|m(A_i)\|_X:n\in \mathbb{N},(A_i)_{i=1,...,n}\subset \mathcal{F} 
\text{ disjoint with }A_i\subset A\}, \, A\subset \Omega\, .
\end{align*}
The restriction of $\overline{m}$ to $\mathcal{F}$ will be denoted by $|m|$. 
We say that $m$ is of finite variation, if $|m|(A) < \infty$ for all $A\in\mathcal{F}$.
\end{definition}
	
\begin{remark}
If $m:\mathcal{F}\rightarrow X$  is a measure of finite variation, then $|m|$ defines a measure on 
$\mathcal{F}$, see e.g. \cite[§3.2~9]{Din66}.
\end{remark}
	
\begin{lemma} 
\label{variationg}
Let $E$ be a Banach space and $g:\Omega\rightarrow E$ be integrable with respect to $\mathbb{P}$. Then
\begin{align*}
n(A) := \mathbb{E}\left[g\mathbf{1}_A\right]
\end{align*}
defines a measure of finite variation and 
\begin{align*}
|n|(A)=\mathbb{E}\left[\|g\|_E\mathbf{1}_A\right].
\end{align*}
We will denote the measure $n$ by $g\mathbb{P}$.
\end{lemma}

\noindent 	
We will shortly recall the proof of the previous Lemma from  \cite[§10.5~Proposition 10]{Din66}.
	
\begin{proof}
It is easy to see that $n$ is indeed a measure. Now let 
\begin{align*}
f(\omega) := \begin{cases}
\frac{g(\omega)}{\|g(\omega)\|_E} & g(\omega)\not=0 \\
x_0\in E \text{ with }\|x_0\|_E=1 & g(\omega)=0.
\end{cases}
\end{align*}
Then $f(\omega)\|g(\omega)\|_E=g(\omega)$, $\|f(\omega)\|_E\equiv 1$ and therefore $f$ is integrable with 
respect to $\lambda = \|g\|_E\mathbb{P}$. Now let $A\in \mathcal{F}$ and $\epsilon>0$, then $f\mathbf{1}_A$ 
is $\lambda$ integrable and there exists a function $h=\sum_{i=1}^{n}x_i\mathbf{1}_{A_i}\in \mathcal{E}_E$ 
with $(A_i)_{i\in \mathbb{N}}\subset \mathcal{F}$ disjoint, $h\in \mathcal{L}^1(\Omega ; E,\lambda)$ and 
\begin{align*}
\int_A\|f-h\|_Ed\lambda<\frac{\epsilon}{2}.
\end{align*}
Therefore 
\begin{align*}
\lambda(A)=\int_A\|f\|_Ed\lambda \le \int_A\|h\|_Ed\lambda + \int_A\|f-h\|_E d\lambda  
<\int \|h\|_E d\lambda+\frac{\epsilon}{2}.
\end{align*}
We may assume that $A_i\subset A$, otherwise consider $A_i\cap A$. For every $i=1,...,n$ there exists 
some $x_i'\in E'$ with $\|x_i'\|_{E'}=1$ and $\langle x_i',x_i\rangle=\|x_i\|_{E}$. Now 
\begin{align*}
h':=\sum_{i=1}^{n}x_i'\mathbf{1}_{A_i}
\end{align*}	
satisfies 
\begin{align*}
\|h'(\omega)\|_{E'}=1,\langle h'(\omega),h(\omega)\rangle=\|h(\omega)\|_E,
\end{align*}
for all $\omega\in \Omega$. Furthermore we have 
\begin{align*}
\int \langle h', f\rangle d\lambda = \sum_{i=1}^n \langle x_i',f\lambda(A_i)\rangle 
=: \int h'd(f\lambda) .
\end{align*}
It follows
\begin{align*}
\lambda(A)-\frac{\epsilon}{2} 
& < \int_A\|h\|_E d\lambda = \int_A \langle h',h\rangle d\lambda \\ 
& \le \left|\int_A \langle h',f\rangle d\lambda\right| + \left|\int_A \langle h',f-h\rangle d\lambda\right| \\
& \le \left|\int_Ah'd(f\lambda)\right| + \int_A \|f-h\|_E\|h'\|_{E'}d\lambda \\
& \le \int_A \|h'\|_{E'}d|f\lambda|+\int_A\|f-h\|_Ed\lambda\\
& \le |f\lambda|(A)+\frac{\epsilon}{2}.
\end{align*}
Therefore $\lambda(A)\le |f\lambda|(A) + \epsilon$, hence $\lambda(A)\le |f\lambda|(A)$, since $\epsilon > 0$ 
was arbitrary, and thus 
\begin{align*}
\lambda\le |f\lambda|\le \|f\|_E\lambda = \lambda.
\end{align*}
In particular it holds $|f\lambda| = \lambda$. Now 
\begin{align*} 
|g\mathbb{P}| = |(f\|g\|_E)\mathbb{P}| = |f (\|g\|_E\mathbb{P})| = |f\lambda| = \lambda = \|g\|_E\mathbb{P}.
\end{align*}
\end{proof}

\subsubsection{Integration w.r.t.~Vector Measures of Finite Variation}
	
We will recall the definition of the integral with respect to vector measures of finite variation. Therefore 
we fix two Banach spaces $E,F$ and a bilinear map
\begin{align*}
\beta:X\times E\rightarrow F,(L,v)\mapsto \beta (L,v) =: Lv,
\end{align*} 
with 
\begin{align*}
\|\beta(L,v)\|_F=\|Lv\|_F\le \|L\|_X\|v\|_E.
\end{align*}
Furthermore let $m:\mathcal{F}\rightarrow X$ be a measure. For a simple function of the type
\begin{align*}
Y = \sum_{i=1}^{n}\mathbf{1}_{A_i}x_i,
\end{align*}
where $A_i\in \mathcal{F}$ and $x_i\in E$ we define 
\begin{align*}
\int_{}^{}Ydm:=\sum_{i=1}^{n}m(A_i)x_i\in F , 
\end{align*}
where the term $m(A_i)x_i$ is to be understood in the sense that $m(A_i)x_i=\beta(m(A_i),x_i)$. 
In the following we denote the set of simple functions with values in $E$ by $\mathcal{E}_E$ and simple 
functions vanishing outside a set $A\in \mathcal{F}$ by $\mathcal{E}_E(A)$. By \cite[§7.2~Proposition 2]{Din66} 
we have the following fundamental inequality.

\begin{lemma} 
\label{stepTriangle}
For any $Y\in\mathcal{E}_E$ it holds
\begin{align*}
\left\|\int Y dm\right\|_F\le \int \|Y\|_Ed|m|.
\end{align*}
\end{lemma}

\noindent 	
A sequence of simple functions $(Y_n)_n\subseteq \mathcal{E}_E$ is called a Cauchy sequence, if 
\begin{align*}
\lim\limits_{n,k\rightarrow \infty}\int_{}^{}\|Y_n-Y_k\|_Ed|m|=0.
\end{align*}

\noindent 
Before we can define the integral with respect to a vector measure of finite variation, we need the following 
uniqueness result from \cite[§7.3 Proposition~9]{Din66}.
	
\begin{lemma}
Let $(Y_n)_n,(Z_n)_n\subseteq \mathcal{E}_E$ be Cauchy sequences converging $|m|$-almost everywhere 
to the same limit $Y:\Omega\rightarrow E$, then 
\begin{align*}
\lim\limits_{n\rightarrow \infty}\int Y_ndm=\lim\limits_{n\rightarrow \infty}\int Z_n dm.
\end{align*}
\end{lemma}
 
\begin{definition}
Let $m:\mathcal{F}\rightarrow X$ be a measure of finite variation. A function $Y:\Omega\rightarrow E$ is said 
to be $m$-integrable, if there exists a Cauchy sequence $(Y_n)_n\subset \mathcal{E}_E$ converging to $Y$, 
$|m|$-almost everywhere. The integral of $Y$ with respect to $m$ is then defined by 
\begin{align*}
\int Ydm := \lim\limits_{n\rightarrow \infty}\int Y_n dm.
\end{align*}

\noindent 		
We will denote the set of $m$-integrable functions $Y:\Omega\rightarrow E$ by $\mathcal{L}^1(\Omega ; E,m)$. 
For a real valued positive measure $\mu:\mathcal{F}\rightarrow \mathbb{R}_+$ the set  
$\mathcal{L}^1 (\Omega ;E,\mu)$ is nothing else but the set of Bochner integrable functions.
\end{definition}
	
\subsubsection{Vector Measures of Finite (Semi) q-Variation}
	
In this subsection we will collect some results for measures of finite (semi) $q$-variation. In particular 
we extend the definition of the integral for measures of finite semi $q$-variation. To this end we consider 
two Banach spaces $E,F$ and a measure $m:\mathcal{F}\rightarrow X\subseteq L(E,F)$. (We can always find 
Banach spaces $E,F$, such that the Banach space $X$ can be isometrically embedded into $L(E,F)$, e.g. 
$E=X',F=\mathbb{R}$). The bilinear map $\beta:X\times E\rightarrow F$ is chosen naturally by $\beta(L,v)=L(v)$.
	
\begin{definition}
For $1\le q\le \infty$ and $\frac{1}{p}+\frac{1}{q}=1$ we define the $q$-variation of $m$ (with respect  
to $\mathbb{P}$) by
\small{
\begin{align*}
\overline{m}_q^{\mathbb{P}}(A) 
& :=\sup\{\sum_{i=1}^{n}\|m(A_i)\|_{op}|\alpha_i|:Y:=\sum_{i=1}^{n}\mathbf{1}_{A_i}\alpha_i 
\in \mathcal{E}_{\mathbb{R}}(A),(A_i)_{i=1,...,n}\subset \mathcal{F}\text{ disjoint},  
\mathbb{E}\left[|Y|^p\right]\le 1\}\\
& =\sup\{\sum_{i=1}^{n}\|m(A_i)x_i\|_F:Y:=\sum_{i=1}^{n}\mathbf{1}_{A_i} x_i 
\in \mathcal{E}_{E}(A),(A_i)_{i=1,...,n}\subset \mathcal{F}\text{ disjoint},\mathbb{E}\left[\|Y\|_E^p\right] 
\le 1\},
\end{align*}
}
for $A\subset \Omega$, and the semi $q$-variation by 
\begin{align*}
\tilde{m}_q^{\mathbb{P}}(A):=\sup\{\|\sum_{i=1}^{n}m(A_i)x_i\|_F: Y:= \sum_{i=1}^{n}\mathbf{1}_{A_i}x_i  
\in \mathcal{E}_E(A),(A_i)_{i=1,...,n}\subset \mathcal{F}\text{ disjoint},\mathbb{E}\left[\|Y\|^p_E\right] 
\le 1\},
\end{align*}
for $A\subset \Omega$.
\end{definition}
	
\begin{definition}\label{measqvar}
For $1\le q<\infty$ we define the space of measures with finite $q$-variation by
\begin{align*}
\mathcal{M}_q^{\mathbb{P}}(E,F) := \{m:\mathcal{F}\rightarrow L(E,F):m\text{ is a measure }, 
\overline{m}_q^{\mathbb{P}}(\Omega)<\infty\},
\end{align*}
which we will need later in the manuscript.
\end{definition}

\noindent 
Let us shortly recall the result from \cite[Corollary~3.6]{BCS15}.
\begin{lemma}
\label{2-var Banach}
The space $(\mathcal{M}_q^{\mathbb{P}}(E,F)$, w.r.t. the norm $\|m\|_{q,\mathbb{P}}  
:= \overline{m}_q^{\mathbb{P}}(\Omega))$ is a Banach space.
\end{lemma}

\noindent 	
We have the following interpolation result from \cite[§13.1 ~9`; §13.1~11; §13.2 11]{Din66}.
	
\begin{lemma}
\label{interpol}
Let $A\in \mathcal{F}$ and $1\le q\le \infty$, $\frac{1}{p}+\frac{1}{q}=1$, then it holds
\begin{align*}
|m|(A)&\le \overline{m}_1^{\mathbb{P}}(A)  
\le \mathbb{P}[A]^{\frac{1}{p}}\overline{m}^{\mathbb{P}}_q(A) 
\le \mathbb{P}[A]\overline{m}^{\mathbb{P}}_{\infty}(A) \\
\tilde{m}(A)&\le \tilde{m}^{\mathbb{P}}_1(A) \le \mathbb{P}[A]^{\frac{1}{p}}\tilde{m}^{\mathbb{P}}_q(A) 
\le \mathbb{P}[A]\tilde{m}^{\mathbb{P}}_{\infty}(A).
\end{align*}
\end{lemma}
	
\noindent 	
The next Lemma follows immediately from Lemma \ref{stepTriangle} and Lemma \ref{interpol}.
	
\begin{lemma} 
\label{YoungLoewe1}
Let $m:\mathcal{F}\rightarrow L(E,F)$ be a measure of finite semi variation and $1\le q\le \infty$, 
$\frac{1}{p}+\frac{1}{q}=1$. Then it holds for any $Y\in\mathcal{E}_E(A)$
\begin{align*}
\|\int_A Ydm\|_F\le \mathbb{E}\left[\|Y\|^p_E\right]^{\frac{1}{p}}\tilde{m}^{\mathbb{P}}_q(A)\le \infty.
\end{align*}
Similarly if $m$ is of finite variation, then it holds for any $Y\in \mathcal{E}_E(A)$
\begin{align*}
\|\int_A Ydm\|_F 
\le \int_A\|Y\|_Ed|m|\le\mathbb{E}\left[\|Y\|^p_E\right]^{\frac{1}{p}}\overline{m}^{\mathbb{P}}_q(A) 
\le \infty.
\end{align*}
\end{lemma}

\noindent 	
With the following elementary result we can define the integration of a random variable  
$Y\in \mathcal{L}^p(\Omega , \mathcal{F},\mathbb{P})$ with respect to a vector measure of finite semi 
$q$-variation.
	
\begin{lemma}
Let $E$ be a separable Banach space. For any $1\le p\le \infty$ the space $\mathcal{E}_E$ is dense in 
$\mathcal{L}^p(\Omega ; E,\mathbb{P})$ and if a sequence $(Y_n)_n$ converges towards $Y$ in  
$\mathcal{L}^p(\Omega ; E,\mathbb{P})$ with respect to the topology induced by the semi-norm  
$\mathcal{N}_p(Y):=\left(\int \|Y\|_E^p d\mathbb{P}\right)^{\frac{1}{p}}$, then there exists a  
subsequence $(Y_{n_k})_k$ converging to $Y$, $\mathbb{P}$-almost everywhere.
\end{lemma}
	
\begin{definition} 
\label{q-semiInt}
Let $m:\mathcal{F}\rightarrow L(F,F)$ be a measure with $\tilde{m}_q^{\mathbb{P}}(A)<\infty$, for all 
$A\in \mathcal{F}$, where $1\le q\le \infty$ and $\frac{1}{p}+\frac{1}{q}=1$. For every  
$Y\in \mathcal{L}^p(\Omega ; E,\mathbb{P})$ we define 
\begin{align*}
\int Ydm := \lim\limits_{n\rightarrow \infty}\int Y_ndm,
\end{align*}
where $(Y_n)_n\subset \mathcal{E}_E$ is an arbitrary Cauchy-sequence of simple functions in  
$\mathcal{L}^p(\Omega ; E,\mathbb{P})$ converging $\mathbb{P}$ almost everywhere to $Y$.
\end{definition}

\noindent 	
The next lemma, see \cite[§13.2 Proposition~4]{Din66}, shows that if $m\in \mathcal{M}_q^{\mathbb{P}}(E,F)$, 
then any function $Y\in \mathcal{L}^p(\Omega ; E, \mathbb{P})$ is also $m$-integrable and the integral  
coincides with the one defined in Definition \ref{q-semiInt}.
	
\begin{lemma} 
\label{RN-iso} 
Let $1\le p<\infty$ and $\frac{1}{p}+\frac{1}{q}=1$. Let $m:\mathcal{F}\rightarrow L(E,F)$ be of finite 
variation, with $|m|\ll \mathbb{P}$. If $A\in \mathcal{F}$ with $\overline{m}_q^{\mathbb{P}}(A)<\infty$, 
then $\mathcal{L}^p(A, E, \mathbb{P})\subset \mathcal{L}^1(A, E,m)$ and
\begin{align*}
\|\int Ydm\|_F\le \int \|Y\|_Ed|m|\le \left(\int \|Y\|^p_E d\mathbb{P}\right)^{\frac{1}{p}}  
\overline{m}_q^{\mathbb{P}}(A).
\end{align*}
Furthermore if $|m|=gd\mathbb{P}$, then $g\mathbf{1}_A\in \mathcal{L}^q(\Omega ; \mathbb{R},\mathbb{P})$  
and
\begin{align*}
\overline{m}_q^{\mathbb{P}}(A)=\left(\int |g|^q\mathbf{1}_Ad\mathbb{P}\right)^{\frac{1}{q}}.
\end{align*}
\end{lemma}
	
\begin{proof}
Let $Y\in \mathcal{L}^p(\Omega ;E,\mathbb{P})$ and $(Y_n)_{n\in \mathbb{N}}\subset \mathcal{E}_E(A)$ such  
that $(Y_n)_{n\in \mathbb{N}}$ converges almost surely and in $\mathcal{L}^p(\Omega ;E,\mathbb{P})$ to $Y$. 
Since $|m|\ll\mathbb{P}$, the sequence $(Y_n)_{n\in \mathbb{N}}$ converges $|m|$-almost everywhere. Now by 
Lemma \ref{YoungLoewe1}
\begin{align*}
\int \|Y_n-Y_m\|_Ed|m|\le \mathbb{E}\left[\|Y_n-Y_m\|_E^p\right]^{\frac{1}{p}}\overline{m}_q^{\mathbb P} (A),
\end{align*}
hence $(Y_n)_{n\in \mathbb{N}}$ is a Cauchy-sequence in $\mathcal{L}^1(\Omega ;E,|m|)$. Therefore  
$Y\in \mathcal{L}^1(\Omega ;E,m)$. Now for any $n\in \mathbb{N}$ it holds 
\begin{align*} 
\|\int Y_ndm\|_F\le \int \|Y_n\|_Ed|m|  
\le \mathbb{E}\left[\|Y_n\|_E^p\right]^{\frac{1}{p}}\overline{m}_q^{\mathbb P}(A).
\end{align*}
Passing to the limit yields
\begin{align*}
\|\int Ydm\|_F \le \int \|Y\|_E d|m| 
\le \mathbb{E}\left[\|Y\|_E^p\right]^{\frac{1}{p}}\overline{m}_q^{\mathbb P} (A).
\end{align*}
Now suppose $|m|=g\mathbb{P}$. Since $|g|\mathbf{1}_A\ge 0$ is measurable, there exists a positive, 
monotone increasing sequence $(g_n)_{n\in \mathbb{N}}\subset \mathcal{E}_{\mathbb{R}}$, such that 
$g_n\le |g|\mathbf{1}_A$ and 
\begin{align*} 
\lim\limits_{n\rightarrow \infty }g_n= |g|\mathbf{1}_A,
\end{align*}
$\mathbb{P}$-almost everywhere. Now since $\mathcal{L}^q(\Omega ,\mathcal{F},\mathbb{P}) 
= \mathcal{L}^p(\Omega ,\mathcal{F},\mathbb{P})'$ and $\mathcal{E}_{\mathbb{R}}$ is dense in 
$\mathcal{L}^p(\Omega ,\mathcal{F},\mathbb{P})$, we get for any $A\in \mathcal{F}$
\begin{align*}
\mathbb{E}\left[g_n^q\mathbf{1}_A\right]^{\frac{1}{q}}  
& = \sup\{\mathbb{E} \left[g_n |Y|\mathbf{1}_A\right]:Y=\sum_{i=1}^{n}\alpha_i\mathbf{1}_{A_i} 
   \in\mathcal{E}_{\mathbb{R}},\mathbb{E}\left[|Y|^p\right]\le 1\} \\
& \le \sup\{\mathbb{E}\left[|g| |Y|\mathbf{1}_A\right]:Y=\sum_{i=1}^{n}\alpha_i\mathbf{1}_{A_i} 
   \in\mathcal{E}_{\mathbb{R}},\mathbb{E}\left[|Y|^p\right]\le 1\} \\
& = \sup\{\int_A |Y|d|m|:Y=\sum_{i=1}^{n}\alpha_i\mathbf{1}_{A_i}  
   \in\mathcal{E}_{\mathbb{R}},\mathbb{E}\left[|Y|^p\right]\le 1\} \\
& \le \overline{m}_q^{\mathbb P} (A).
\end{align*} 
Therefore 
\begin{align*}
\mathbb{E}\left[|g|^q\mathbf{1}_A\right]^{\frac{1}{q}} 
= \lim_{n\to\infty} \mathbb{E}\left[g_n^q\mathbf{1}_A\right]^{\frac{1}{q}} 
\le \overline{m}_q^{\mathbb P} (A)  
\end{align*}
by monotone convergence, in particular $g\mathbf{1}_A \in \mathcal{L}^q(\Omega, \mathcal{F},\mathbb{P})$. 
The converse inequality, $\overline{m}_q^{\mathbb P} (A) \le \mathbb{E}\left[|g|^q\mathbf{1}_A 
\right]^{\frac{1}{q}} $, follows from the fact that for any $Y=\sum_{i=1}^{n}\alpha_i\mathbf{1}_{A_i}\in 
\mathcal{E}_{\mathbb{R}}(A)$ with $\mathbb{E}\left[|Y|^p\right]^{\frac{1}{p}}\le 1$ we have that 
\begin{align*}
\sum_{i=1}^{n}\|m(A_i)\||\alpha_i| 
& \le \sum_{i=1}^{n}|m|(A_i)|\alpha_i| =\int |Y|d|m| =\mathbb{E}\left[ |Y|g\mathbf{1}_A\right] \\
& \le \mathbb{E}\left[|Y|^p\right]^{\frac{1}{p}}\mathbb{E}\left[|g|^q\mathbf{1}_A\right]^{\frac{1}{q}} 
	\le \mathbb{E} \left[|g|^q\mathbf{1}_A\right]^{\frac{1}{q}},
\end{align*}
hence $\overline{m}_q^{\mathbb P} (A)\le \mathbb{E} \left[|g|^q\mathbf{1}_A\right]^{\frac{1}{q}}$.
\end{proof}

\subsection{Linear Operators on $\mathcal{L}^p$}
\label{LinOp}
	
In this section we will recall the fundamental relation between linear operators and vector measures  
stated in \cite{Din66}. The corresponding isomorphism between linear operators and vector measures, see 
Theorem \ref{meas1}, will be key for the proof of our main theorem. We fix again two Banach spaces $E,F$.
	
\begin{definition}
\label{varNorm}
Let $L:\mathcal{L}^p(\Omega ;E,\mathbb{P})\rightarrow F$ be a linear mapping. For $A\in \mathcal{F}$ we 
denote by $L_A$ the restriction of $L$ onto $\mathcal{L}^p ( A;E,\mathbb{P})$. Then we define
\small{
\begin{align*}
\|L_A\|_{p,\mathbb{P}}&:=\sup\{\|LY\|_F:Y\in \mathcal{L}^p(A;E,\mathbb{P})\text{ with }  
\mathbb{E}\left[\|Y\|_E^p\right]\le 1\} \\  
\vertiii{L_A}_{p,\mathbb{P}} & = \sup\{\sum_{i=1}^{n} \|L(\mathbf{1}_{A_i}x_i)\|_F:Y=\sum_{i=1}^{n}
\mathbf{1}_{A_i}x_i\in\mathcal{E}_E(A), (A_i)_{i=1,...,n}\subset \mathcal{F},   
A_i\text{ disjoint},\mathbb{E}\left[\|Y\|_E^p\right]\le 1\}.
\end{align*} 
}
\end{definition}
	
\begin{definition}
We write $\|L\|_{p,\mathbb{P}}=\|L_{\Omega}\|_{p,\mathbb{P}}$, $\vertiii{L}_{p,\mathbb{P}} 
= \| L_{\Omega}\|_{p,\mathbb{P}}$ and define 
\begin{align*}
\Lambda_p^{\mathbb{P}}(E,F) 
:= \{L\in L(\mathcal{L}^p(\Omega ; E,\mathbb{P}),F):\vertiii{L}_{p,\mathbb{P}} < \infty\}.
\end{align*}
\end{definition}

\noindent 
The following is an immediate consequence of the previous definition.

\begin{lemma}
Let $L\in L(\mathcal{L}^p(\Omega; E,\mathbb{P}),F)$. For any $A\in \mathcal{F}$ it holds 
\begin{align*}
\|L_A\|_{p,\mathbb{P}}\le \vertiii{L_A}_{p,\mathbb{P}}.
\end{align*}
\end{lemma}
	
\begin{lemma}
$\vertiii{\cdot}_{p,\mathbb{P}}$ defines a norm on $\Lambda_p^{\mathbb{P}}(E,F)$. Furthermore the space 
$(\Lambda_p^{\mathbb{P}}(E,F),\vertiii{\cdot}_{p,\mathbb{P}})$ is a Banach space.
\end{lemma}
	
\begin{proof}
It is easy to check that $\vertiii{\cdot}_{p,\mathbb{P}}$ defines a norm on $\Lambda_p^{\mathbb{P}}(E,F)$. 
Now let $(L_k)_{k\in \mathbb{N}}$ be a Cauchy-sequence in $\Lambda_p^{\mathbb{P}}(E,F)$. Since
\begin{align*}
\|L_m-L_k\|_{p,\mathbb{P}}\le \vertiii{L_m-L_n}_{p,\mathbb{P}},
\end{align*}
$(L_k)_{k\in \mathbb{N}}$ is a Cauchy-sequence in $L(\mathcal{L}^p(\Omega ; E,\mathbb{P}),F)$, hence 
there exists a $L\in  L(\mathcal{L}^p(\Omega; E,\mathbb{P}),F)$ with 
\begin{align*}
\|L-L_k\|_{p,\mathbb{P}}\rightarrow 0. 
\end{align*}
Now let $\epsilon>0$ and $N_{\epsilon}\in \mathbb{N}$, such that $\vertiii{L_m-L_k}_{p,\mathbb{P}} < \epsilon$
for all $m,k\ge N_{\epsilon}$. Then it holds for any $m\ge N_{\epsilon}$ and  
$Y=\sum_{i=1}^{n} x_i\mathbf{1}_{A_i}\in \mathcal{E}_E$ with $(A_i)_i$ disjoint
\begin{align*}
\sum_{i=1}^{n}\|L(\mathbf{1}_{A_i}x_i)-L_m(\mathbf{1}_{A_i}x_i)\|_F 
& =\lim\limits_{k\rightarrow \infty}\sum_{i=1}^{n}\|L_k(\mathbf{1}_{A_i}x_i)-L_m(\mathbf{1}_{A_i}x_i)\|_F \\ 
& \le \lim\limits_{k\rightarrow \infty}\vertiii{L_k-L_m}_{p,\mathbb{P}}  
   \mathbb{E}\left[\|Y\|^p\right]^{\frac{1}{p}} \\
& \le \epsilon \mathbb{E}\left[\|Y \|^p\right]^{\frac{1}{p}}.
\end{align*}
Therefore we obtain for all $Y=\sum_{i=1}^{n}x_i\mathbf{1}_{A_i}\in \mathcal{E}_E$ with $(A_i)_i$ disjoint
\begin{align*}
\sum_{i=1}^{n}\|L(\mathbf{1}_{A_i}x_i)\|_F 
& \le \epsilon \mathbb{E}\left[\|Y\|^p\right]^{\frac{1}{p}} + \vertiii{L_{N_{\epsilon}}}_{p,\mathbb{P}}
\end{align*}
and therefore $\vertiii{L}_{p,\mathbb{P}}<\infty$. Furthermore, $\vertiii{L-L_m}_{p,\mathbb{P}}\le 
\epsilon$ for all $m\ge N_{\epsilon}$, hence 
\begin{align*}
\lim_{m\to \infty} \vertiii{L-L_m}_{p,\mathbb{P}} = 0.
\end{align*}
\end{proof}

\noindent 
The following lemma, see \cite[§13.3 Proposition~5]{Din66}, shows that the assumptions of our main 
theorem are indeed not stronger as the ones in \cite{CGKPR22} if $U=\mathbb{R}$.
	
\begin{lemma} 
\label{variation1d}
Let $L:\mathcal{L}^p(\Omega ;E,\mathbb{P})\rightarrow F$ be a linear mapping. If $p=1$, or if $F=\mathbb{R}$,  
we have for any $A\in \mathcal{F}$
\begin{align*}
\|L_A\|_{p,\mathbb{P}} = \vertiii{L_A}_{p,\mathbb{P}}.
\end{align*}
\end{lemma}
	
\begin{proof}
We will shortly recall the proof from \cite{Din66} in the case $F=\mathbb{R}$. Since  
$\|L_A\|_{p,\mathbb{P}}\le\vertiii{L_A}_{p,\mathbb{P}}$ it remains to show  
$\vertiii{L_A}_{p,\mathbb{P}}\le \|L_A\|_{p,\mathbb{P}}$. To this end fix $Y\in \mathcal{E}_\mathbb{R}$,  
$Y = \sum_{i=1}^{n}\mathbf{1}_{A_i} x_i$. Let $\theta_i=1$ if $L(\mathbf{1}_{A_i}x_i)\ge 0$ and 
$\theta_i=-1$ if $L(\mathbf{1}_{A_i}x_i)<0$. Then 
\begin{align*}
\sum_{i=1}^{n}|L(\mathbf{1}_{A_i}x_i)|&=L(\sum_{i=1}^{n}\mathbf{1}_{A_i}\theta_ix_i) 
\le \|L_A\|_{p,\mathbb{P}} \mathbb{E}\left[\|Y\|_E^p\right]^{\frac{1}{p}} 
\end{align*}
which finishes the proof. 
\end{proof}

\noindent 
We now come to the important isomorphism between vector measures and linear operators, see  
\cite[§13.3 Theorem~1; §13.3 Corollary~1]{Din66}.
	
\begin{theorem} 
\label{meas1} 
Let $1\le p<\infty$ and $\frac{1}{p}+\frac{1}{q}=1$. Then there exists an isomorphism between the set 
of linear mappings $L:\mathcal{L}^p(\Omega ;E,\mathbb{P})\rightarrow F$ with $\|L_A\|_p<\infty$ for all 
$A\in \mathcal{F}$ and the set of measures $m:\mathcal{F}\rightarrow L(E,F)$ with 
$\tilde{m}_q^{\mathbb{P}}(A)<\infty$, for all $A\in \mathcal{F}$. The isomorphism is given by 
\begin{align*}
LY=\int Ydm,
\end{align*}
for $Y\in \mathcal{L}^p(\Omega ; E,\mathbb{P})$ and for all $A\in\mathcal{F}$ 
\begin{align*}
\|L_A\|_{p,\mathbb{P}}=\tilde{m}^{\mathbb{P}}_q(A),\quad \vertiii{L_A}_{p,\mathbb{P}} 
= \overline{m}^{\mathbb{P}}_q(A).
\end{align*}
\end{theorem}
	
\begin{corollary} 
\label{cor1}
Let $1\le p<\infty$ and $\frac{1}{p}+\frac{1}{q}=1$. Then there exists an isometric isomorphism between 
$(\Lambda_p^{\mathbb{P}}(E,F),\vertiii{\cdot}_{p,\mathbb{P}})$ and $(\mathcal{M}_q^{\mathbb{P}}(E,F), 
\|\cdot\|_{q,\mathbb{P}})$. The isomorphism is given by 
\begin{align*}
LY=\int Ydm,
\end{align*}
for $Y\in \mathcal{L}^p(\Omega ; E,\mathbb{P})$.
\end{corollary}
	
\begin{proof}
First we consider a measure $m$ with $\tilde{m}_q^{\mathbb{P}}(A)<\infty$ for all $A\in \mathcal{F}$. 
The map
\begin{align*}
L:\mathcal{L}^p(\Omega ; E,\mathbb{P})\rightarrow F,LY:=\int Ydm
\end{align*}
is clearly linear. For every $A\subset \Omega$ and $Y=\sum_{i=1}^{n}x_i\mathbf{1}_{A_i}\in\mathcal{E}_E(A)$ 
with $A_i$ disjoint, we have 
\begin{align*}
\|LY\|_F = \|\sum_{i=1}^{n}m(A_i)x_i\|_F
\end{align*}
and 
\begin{align*}
\sum_{i=1}^{n}\|L(\mathbf{1}_{A_i}x_i)\|_F=\sum_{i=1}^{n}\|m(A_i)x_i\|_F.
\end{align*}
Taking the supremum over all simple functions with $\mathbb{E}\left| \|Y\|_E^p\right]\le 1$ we deduce
\begin{align*}
\|L_A\|_p = \tilde{m}_q(A),\quad \vertiii{L_A}_p=\overline{m}_q(A).
\end{align*}
Conversely, if we consider a linear mapping $L:\mathcal{L}^p(\Omega ;E,\mathbb{P})\rightarrow F$ with 
$\|L\|_{p,\mathbb{P}}<\infty$, we define $m:\mathcal{F}\rightarrow L(E,F)$ by 
\begin{align*}
m(A)u:=L(\mathbf{1}_Au).
\end{align*}
For any $A\in \mathcal{F}$, $m(A):E\rightarrow F$ is obviously linear and it holds 
\begin{align*}
\|m(A)u\|_F&\le \|L_A\|_{p,\mathbb{P}}\|u\|_E\mathbb{E}\left( \mathbf{1}_A \right)^{\frac{1}{p}},
\end{align*}
hence $m(A)$ is continuous, i.e. $m:\mathcal{F}\rightarrow L(E,F)$. The mapping $m$ is clearly additive 
and for every simple function $Y=\sum_{i=1}^{n}u_i\mathbf{1}_{A_i}$ it holds 
\begin{equation} 
\label{linear}
LY = \sum_{i=1}^{n}L(\mathbf{1}_{A_i}u_i)=\sum_{i=1}^{n}m(A_i)u_i=\int Ydm.
\end{equation}
If we consider countably many mutually disjoint sets $(A_i)_{i\ge 1}\subset \mathcal{F}$, then 
\begin{align*}
Y_n^u := \sum_{i=1}^{n} u \mathbf{1}_{A_i} = u\mathbf{1}_{\bigcup_{i=1}^{n} A_i}\rightarrow Y^u 
:= \sum_{i=1}^{\infty} u \mathbf{1}_{A_i} = u\mathbf{1}_{\bigcup_{i=1}^{\infty}A_i}
\end{align*}
$\mathbb{P}$-a.s and in $\mathcal{L}^p(\Omega ;E,\mathbb{P})$, for any $u\in E$. In particular
\begin{align*}
\mathbb{E}\left[\|Y_n^u-Y^u\|^p_E\right] & =\|u\|^p_E\mathbb{P}\left[\bigcup_{i=n+1}^{\infty}A_i\right].
\end{align*} 
Therefore 
\begin{align*}
\sup_{u\in E,\|u\|_E\le 1}\| \sum_{i=1}^{n}m(A_i)u-m(\bigcup_{i=1}^{\infty}A_i)u\|_F 
& =\sup_{u\in E,\|u\|_E\le 1} \left\|\int (Y_n^u-Y^u)dm\right\|_F \\
& \le \mathbb{P}\left[ \bigcup_{i=n+1}^{\infty}A_i \right]^{\frac{1}{p}}\, \tilde{m}^{\mathbb{P}}_q 
(\Omega)\rightarrow 0.
\end{align*}
Thus $m$ is indeed countably additive and therefore a measure. From \eqref{linear} it is easy to deduce 
\begin{align*}
\vertiii{L_A}_{p,\mathbb{P}} = \bar{m}_q^{\mathbb{P}}(A),
\end{align*}
for every $A\in \mathcal{F}$, hence 
\begin{align*}
\bar{m}_q^{\mathbb{P}}(A) < \infty, 
\end{align*} 
for all $A\in \mathcal{F}$. Now for $Y\in\mathcal{L}^p(A ;E,\mathbb{P})$ consider a sequence  
$(Y_n)_{n\in \mathbb{N}}\subset \mathcal{L}^p(A;E,\mathbb{P})$ of simple functions converging to  
$Y$ in $\mathcal{L}^p(A;E,\mathbb{P})$ and almost surely. Then $LY_n\rightarrow LY$ and  
$\int Y_n\, dm\rightarrow \int Y\, dm$. Now $LY_n=\int Y_n\, dm$ for every $n\in \mathbb{N}$ implies that 
\begin{align*}
LY=\int Y\, dm.
\end{align*}
Finally the correspondence $m\leftrightarrow L$ is evidently linear and one-to-one, hence $m\mapsto L$ 
is a vector space isomorphism. 
\end{proof}
	

\section{Generalized L-Derivative, Proof of Theorem \ref{g}} 
\label{generalL}
	
We will now turn to the generalization of the L-derivative for $U$-valued functions  
$f:\mathcal{P}_2(H)\rightarrow U$. Similar to the finite dimensional case we want to prove that the 
notion of the $L$-derivative is intrinsic, and give some results concerning the structure of the 
$L$-derivative. Therefore we will derive a factorized gradient-type representation of $D\hat{f}(X_0)$ 
that only depends on $\mathcal{L}(X_0)$, but neither on the random variable $X_0$ itself nor on the 
underlying probability space $(\Omega, \mathcal{F},\mathbb{P})$. We will shortly explain the main 
intuition behind the proof of our main Theorem:
	
Let $f:\mathcal{P}_2(H)\rightarrow U$ be $L$-differentiable in $\mu_0\in \mathcal{P}_2(H)$ in the sense of 
Definition \ref{Ldiff1}. The directional derivative of the lift $\hat{f}$ at some $X_0\in L^2(\Omega; H)$ 
with $\mathcal{L}(X_0)=\mu_0$ in direction $u\mathbf{1}_A$,
\begin{align*}
D\hat{f}(X_0)(u\mathbf{1}_A)
& = \lim\limits_{\epsilon\downarrow 0}\frac{f(\mathcal{L}(X_0+\epsilon u\mathbf{1}_{A})) 
   -f(\mathcal{L}(X_0))}{\epsilon},
\end{align*} 
which describes the infinitesimal variation of the function $f$, if mass $A$ is shifted towards the point 
$u\in H$, is nothing but the vector measure $m_{D\hat{f}(X_0)}$ associated to $D\hat{f}(X_0)$ by Theorem 
\ref{meas1}, evaluated in $A$ and $u$, i.e.
\begin{align*}
m_{D\hat{f}(X_0)}(A)u=D\hat{f}(X_0)(u\mathbf{1}_A).
\end{align*}
In particular, for every $Y\in L^2(\Omega ;H)$ we have 
\begin{align*}
D\hat{f}(X_0)Y=\int Ydm_{D\hat{f}(X_0)}.
\end{align*}
As already mentioned in the introduction, our main strategy is to use the results on vector measures to 
derive our factorized gradient representation as a Radon-Nikodym derivative corresponding to $m_{D\hat{f}(X_0)}$ 
and $\mathbb{P}$. As a first naive approach we could consider the Radon-Nikodym derivative   
$\frac{dm_{D\hat{f}(X_0)}}{d\mathbb{P}}:\Omega\rightarrow L(H,U)$ of $m_{D\hat{f}(X_0)}$ with respect 
to $\mathbb{P}$, for the existence see e.g. \cite[§13,Theorem~8]{Din66}, which satisfies
\begin{align*}
D\hat{f}(X_0)Y=\int Ydm_{D\hat{f}(X_0)}&=\mathbb{E}\left[\frac{dm_{D\hat{f}(X_0)}}{d\mathbb{P}}Y\right],
\end{align*}
if $U$ is separable. The Radon-Nikodym derivative corresponds to the Riesz representation $D\hat{f}(X_0)\in 
L^2(\Omega ; H)$ in the case $U=\mathbb{R}$. Despite the required separability of $U$, the main problem is 
that it is not clear how $m_{D\hat{f}(X_0)}$, and in particular the Radon-Nikodym derivative, depends on the 
random variable $X_0$ and the underlying probability space. Therefore, similar to \cite{CD18} in the 
finite dimensional case, we would have to show that $\frac{dm_{D\hat{f}(X_0)}}{d\mathbb{P}}$ is 
$X_0$-measurable and that its 
factorization only depends on $\mathcal{L}(X_0)$. Instead we follow a different approach which is based on the 
disintegration of our vector measure. Therefore we take a step back and first observe that the image measure 
$m_{D\hat{f}(X_0)}\circ X_0^{-1}$ will only depend on $\mu_0=\mathcal{L}(X_0)$ and not on $X_0$ as a random 
variable itself, i.e. if $\mathcal{L}(X_0) = \mathcal{L}(\tilde{X}_0)$, then 
\begin{align*}
m_{D\hat{f}(X_0)}(X_0\in B)=m_{D\hat{f}(\tilde{X}_0)}(\tilde{X}_0\in B),
\end{align*} 
for all $B\in \mathcal{B}(H)$. Indeed, for $B\in \mathcal{B}(H)$, let $\phi (x) := x + u\mathbf{1}_B (x)$. Then 
$$ 
\mathcal{L}(X_0+u\mathbf{1}_{\{X_0\in B\}}) = \mathcal{L}(\phi (X_0)) 
= \mathcal{L}(\phi (\tilde{X}_0)) 
=  \mathcal{L}(\tilde{X}_0+u\mathbf{1}_{\{\tilde{X}_0\in B\}})
$$
Now given another lift of $f$
\begin{align*}
\tilde{f}:L^2(\tilde{\Omega};H)\rightarrow U,\tilde{X}\mapsto f(\mathcal{L}(\tilde{X})),
\end{align*}
based on the probability space $(\tilde{\Omega},\tilde{\mathcal{F}},\tilde{\mathbb{P}})$, which is 
Fr\'echet differentiable in $\tilde{X}_0\sim \mu_0$, we have
\begin{align*}
m_{D\hat{f}(X_0)}(X_0\in B)u&=D\hat{f}(X_0)(u\mathbf{1}_{\{X_0\in B\}}) \\
& = \lim\limits_{\epsilon\downarrow 0}\frac{f(\mathcal{L}(X_0+\epsilon u\mathbf{1}_{X_0\in B})) 
    -f(\mathcal{L}(X_0))}{\epsilon} \\ 
& = \lim\limits_{\epsilon\downarrow 0}\frac{f(\mathcal{L}(\tilde{X}_0 
	+\epsilon u\mathbf{1}_{\tilde{X}_0\in B}))-f(\mathcal{L}(\tilde{X}_0))}{\epsilon} \\
& = D\tilde{f}(\tilde{X}_0)(u\mathbf{1}_{\{\tilde{X}_0\in B\}}) \\
& = m_{D\hat{f}(\tilde{X}_0)}(\tilde{X}_0\in B)u.
\end{align*}
This means that $m_{D\hat{f}(X_0)}\circ X_0^{-1}$ does only depend on $\mathcal{L}(X_0)$, but neither on the 
random variable $X_0$ nor the lift considered for the L-derivative. We will write $m_{\mu_0}$ instead of  
$m_{D\hat{f}(X_0)}\circ X_0^{-1}$ in the following. The second observation will be that we can disintegrate 
$m_{D\hat{f}(X_0)}$ in the following way
\begin{equation} 
\label{condProbm}
\begin{aligned}
m_{D\hat{f}(X_0)}(A)u = \int \mathbb{P}[A|X_0]udm_{D\hat{f}(X_0)} 
& =\int \mathbb{P}[A|X_0=x]udm_{\mu_0}(x).
\end{aligned}
\end{equation}
This means in particular that $m_{D\hat{f}(X_0)}(A)$ only depends on the $\sigma(X_0)$-measurable 
$L^2$-projection of $A$. Once we have proven the disintegration, we can derive our representation as the 
Radon-Nikodym derivative $\frac{dm_{\mu_0}}{d\mu_0}$ of $m_{\mu_0}$ w.r.t. $\mu_0=\mathcal{L}(X_0)$, which 
will indeed	only depend on $\mu_0$. This will lead to 
\begin{align*}
m_{D\hat{f}(X_0)}(A)u=\int_H  \mathbb{P}[A|X_0=x] \frac{dm_{\mu_0}}{d\mu_0}(x)u d\mu_0(dx).
\end{align*}
It is then straightforward to derive the following representation 
\begin{align*}
D\hat{f}(X_0)Y
& = \int Ydm_{D\hat{f}(X_0)} \\
& = \int \frac{dm_{\mu_0}}{d\mu_0}(x) \mathbb{E}\left[Y|X_0=x\right] d\mu_0(x) \\
& = \mathbb{E}\left[\frac{dm_{\mu_0}}{d\mu_0}(X_0)Y\right].
\end{align*}

\noindent 	
The proof of our main theorem will be divided into two parts. In the first part we consider the simple case 
where $X_0$ only takes finitely many values. In the second part we will prove our main theorem for the general 
case. We will assume in the following that the assumptions of Theorem \ref{g} are satisfied.
	
\subsection{The Discrete Case} 
\label{discrete}
	
We will begin with the simple case, where $\mu_0=\mathcal{L}(X_0)\in \mathcal{P}_2(H)$ is of the type
\begin{align*}
\mathcal{L}(X_0)=\sum_{k=1}^{N}p_k\delta_{x_k},
\end{align*}
for some $n\in \mathbb{N}$, $(x_k)_{k=1,...,N}\subset H$ and $p_i>0$ with $\sum_{k=1}^{N}p_k=1$. Again we 
consider the corresponding vector measure $m_{D\hat{f}(X_0)}\in \mathcal{M}_2^{\mathbb{P}}(H,U)$ associated 
to $D\hat{f}(X_0)\in \Lambda_2^{\mathbb{P}}(H,U)$ by Corollary \ref{cor1}. Similar to \cite{CGKPR22} we can 
observe that for every $i=1,...,N$ and $A_1,A_2\subset \{X_0=x_i\}$ with $\mathbb{P}[A_1]=\mathbb{P}[A_2]$ 
it holds 
\begin{align*}
m_{D\hat{f}(X_0)}(A_1)u 
& = D\hat{f}(X_0)(u\mathbf{1}_{A_1}) \\
& = \lim\limits_{\epsilon\downarrow 0}\frac{f(\sum_{j\not=i}^{}p_j\delta_{x_j}  
	+ \mathbb{P}[A_1]\delta_{x_i+\epsilon u} + (p_i-\mathbb{P}[A_1])\delta_{x_i}) 
	-f(\sum_{i=1}^{N}p_i\delta_{x_i})}{\epsilon} \\
& = \lim\limits_{\epsilon\downarrow 0}\frac{f(\sum_{j\not=i}^{}p_j\delta_{x_j} + 
	\mathbb{P}[A_2]\delta_{x_i+\epsilon u}+(p_i-\mathbb{P}[A_2])\delta_{x_i})  
	-f(\sum_{i=1}^{N}p_i\delta_{x_i})}{\epsilon} \\
& = m_{D\hat{f}(X_0)}(A_2)u.
\end{align*}
In particular, if we shift the mass $A_1$ away from the point $x_i$ in direction of $u$, the directional 
derivative of the function $f$ in direction $u\mathbf{1}_{A_1}$ only depends on the probability 
$\mathbb{P}[A_1]$ of the set $A_1$, but not on the set $A_1$ itself as long as $A_1\subset \{X_0 = x_i\}$ 
is measurable. With the next Lemma we will see that this dependence is in general proportional to the  
conditional probability of $A_1$ given $X_0=x_i$.
	
\begin{lemma} 
\label{scale}
Let $m\in \mathcal{M}^{\mathbb{P}}_q(H,U)$ for some $1<q<\infty$. Furthermore let $B\in \mathcal{F}$ with 
$\mathbb{P}[B]>0$ and
\begin{align*}
m(A_1)u = m(A_2)u,\forall u\in H,
\end{align*}
for all $A_1,A_2\subset B$ with $\mathbb{P}[A_1]=\mathbb{P}[A_2]$. Then it holds
\begin{align*}
\mathbb{P}[A|B] m(B) = m(A\cap B),
\end{align*}
for all $A\in \mathcal{F}$.
\end{lemma}
	
\begin{proof}
For $u\in H$, $v\in U'$, we define the real-valued measure     	
\begin{align*}
m_{u,v}(A) := \langle m(A)u,v\rangle.
\end{align*}
Since 
\begin{align*}
|m_{u,v}|(A)&\le \|u\|_H\|v\|_{U'}|m|(A)
\end{align*}
and by Lemma \ref{interpol} $|m|  \ll\mathbb{P}$, we conclude that $|m_{u,v}| \ll\mathbb{P}$. Therefore the 
Lebesgue-Nikodym theorem implies the existence of an integrable function $g_{u,v}$ such that for every $A\in 
\mathcal{F}$
\begin{align*}
m_{u,v}(A)=\mathbb{E}\left[\mathbf{1}_A g_{u,v}\right].
\end{align*}
Since $m\in \mathcal{M}_q^{\mathbb{P}}(H,U)$ we have by Lemma \ref{RN-iso} 
\begin{align*}
g_{u,v}\in L^q(\Omega,\mathbb{P}).
\end{align*}
Now for any $A_1,A_2\subset B$ with $\mathbb{P}[A_1]=\mathbb{P}[A_2]$ it holds 
\begin{align*}
\mathbb{E}\left[\mathbf{1}_{A_1}g_{u,v}\right] 
 =\langle m(A_1)u,v\rangle 
 = \langle m(A_2)u,v\rangle 
 = \mathbb{E}\left[\mathbf{1}_{A_2}g_{u,v}\right].
\end{align*}
Therefore by \cite[Lemma~2]{WZ18}, $g_{u,v}$ is almost surely constant on $B$, hence for all $A\in \mathcal{F}$ 
\begin{align*}
m_{u,v}(B)\mathbb{P}[A|B] 
& = \mathbb{E}\left[\frac{m_{u,v}(B)}{\mathbb{P}[B]}\mathbf{1}_{A\cap B}\right] 
= \mathbb{E}\left[\frac{1}{\mathbb{P}[B]}\mathbb{E}\left[g_{u,v}\mathbf{1}_B\right] 
\mathbf{1}_{A\cap B}\right] \\
& = \mathbb{E}\left[g_{u,v}\mathbf{1}_{A\cap B}\right]
 = m_{u,v}(A\cap B)
 = \langle m(A\cap B)u,v\rangle.
\end{align*} 
Therefore 
\begin{align*}
\langle \mathbb{P}[A|B]m(B)u-m(A\cap B)u,v\rangle = 0,
\end{align*}
for every $u\in H,v\in U'$, hence $\mathbb{P}[A|B]m(B)=m(A\cap B)$.
\end{proof}

\noindent 
We can now prove our disintegration theorem in the discrete setting. Let us first state the following 
preparatory  
\begin{lemma}
For $A\in\mathcal F$, $u\in H$ it holds 
$$
D\hat{f} (X_0) (\mathbf{1}_{A\cap \{X_0 = x_k\}} u)  
= \mathbb{P} [A|X_0 = x_k] m_{D\hat{f} (X_0)}(X_0 = x_k)u. 
$$
\end{lemma}
	
\begin{proof} 
Based on the discussion at the beginning of this section we have for every $A_1$, $A_2\in\mathcal F$ with 
$A_1,A_2\subset\{X_0 = x_k\}$ and $\mathbb{P} [A_1] = \mathbb{P} [A_2]$ 
$$
m_{D\hat{f} (X_0)}(A_1)u = m_{D\hat{f} (X_0)}(A_2)u. 
$$
Therefore Lemma \ref{scale} implies the assertion. 
\end{proof} 

\begin{theorem}
\label{disintDiscrete}
Let $\mu_0=\mathcal{L}(X_0)\in \mathcal{P}_2(H)$ be of the type
\begin{align*}
\mathcal{L}(X_0) = \sum_{k=1}^{N}p_k\delta_{x_k},
\end{align*}
for some $n\in \mathbb{N}$, $(x_k)_{k=1,...,N}\subset H$ and $p_i>0$ with $\sum_{k=1}^{N}p_k=1$. If  
$m_{D\hat{f}(X_0)}$ is the vector measure associated to $D\hat{f}(X_0)$ by Corollary \ref{cor1}, then it 
holds for all $A\in \mathcal{F}$ and $u\in H$ 
\begin{align*}
m_{D\hat{f}(X_0)}(A)u=\int \mathbb{P}[A|X_0]udm_{D\hat{f}(X_0)}.
\end{align*}
\end{theorem}

\begin{proof}
The previous Lemma implies that 
\begin{align*}
D\hat{f} (X_0) (\mathbf{1}_A u)
& = \sum_{k=1}^N D\hat{f} (X_0) (\mathbf{1}_{A\cap\{X_0 = x_k\}} u) \\ 
& = \sum_{k=1}^N \mathbb{P} [A | X_0 = x_k] m_{D\hat{f} (X_0)} (X_0 = x_k)u \\ 
& = \int \sum_{k=1}^N \mathbb{P} [A | X_0 = x_k]u 
	\mathbf{1}_{\{X_0 = x_k\}} dm_{D\hat{f} (X_0)} \\ 
& = \int \mathbb{E} [\mathbf{1}_A | X_0 ]u dm_{D\hat{f} (X_0)}. 
\end{align*} 
\end{proof}

The Radon-Nikodym derivative of $m_{\mu_0}$ with respect to $\mu_0$ is now given as follows. 

\begin{lemma}
\label{RN-discrete}
Let $\mu_0=\mathcal{L}(X_0)\in \mathcal{P}_2(H)$ be of the type
\begin{align*} 
\mathcal{L}(X_0)=\sum_{k=1}^{N}p_k\delta_{x_k},
\end{align*}
for some $N\in \mathbb{N}$, $(x_k)_{k=1,...,N}\subset H$ and $p_i>0$ with $\sum_{k=1}^{N}p_k=1$. If  
$m_{D\hat{f}(X_0)}$ is the vector measure associated to $D\hat{f}(X_0)$ by Corollary \ref{cor1} and 
$m_{\mu_0}=m_{D\hat{f}(X_0)}\circ X_0^{-1}$, then 
\begin{align*}
\frac{dm_{\mu_0}}{d\mu_0}(x):=\sum_{k=1}^{N}\frac{m_{\mu_0}(\{x_k\})}{\mu_0(\{x_k\})}\mathbf{1}_{\{x_k\}}(x)
\end{align*}
satisfies 
\begin{align*}
\mathbb{E}\left[\mathbf{1}_A\frac{dm_{\mu_0}}{d\mu_0}(X_0)\right]=m_{D\hat{f}(X_0)}(A),
\end{align*}
for all $A\in \mathcal{F}$. 
\end{lemma}

\begin{proof}
First it is clear that $\frac{dm_{\mu_0}}{d\mu_0}(X_0)$ is Bochner integrable with respect to $\mathbb{P}$, 
since it is nothing but a simple function from $(\Omega, \mathcal{F},\mathbb{P})$ to $(L(H,U),\|\cdot\|_{op})$. 
Then we have by Theorem \ref{disintDiscrete}
\begin{align*}
\mathbb{E}\left[\mathbf{1}_A\frac{dm_{\mu_0}}{d\mu_0}(X_0)\right] 
& = \sum_{k=1}^{n} \frac{m_{D\hat{f}(X_0)}(X_0=x_k)}{\mathbb{P}[X_0=x_k]}\mathbb{P}(\{X_0=x_k\}\cap A) \\
& = \int \mathbb{P}\left[A|X_0\right]dm_{D\hat{f}(X_0)} \\
& = m_{D\hat{f}(X_0)}(A).
\end{align*}
\end{proof}
	
\subsection{The General Case}

We are now in the position to prove our main theorem, where we derive our factorized gradient representation 
of $D\hat{f}(X)$ when $\mu=\mathcal{L}(X)\in \mathcal{P}_2(H)$ is not necessarily supported on a finite set. 
	
\begin{proofThm}{\ref{g}}
We fix some $\mu\in \mathcal{P}_2(H)$ and some arbitrary $X\in L^2(\Omega ; H)$ with $\mathcal{L}(X)=\mu$. 
Let $(D_i^n)_{i\in \mathbb{N}}$ be a partition of Borel sets of $H$, such that $diam(D_i^n)<2^{-n}$, where 
$diam(A):=\sup\{\|x-y\|_H|x,y\in A\}$. Such a partition exists due to the separability of $H$. For $i=1,2,...$ 
we choose $d_i^n\in D_i^n$ and define
\begin{align*}
\xi_n:=\sum_{i=1}^{\infty}d_i^n\mathbf{1}_{D_i^n}(X).
\end{align*} 
It is not difficult to see that for every $\omega\in \Omega$ and $n\in \mathbb{N}$
\begin{align*}
\|\xi_n(\omega)-X(\omega)\|_H<2^{-n}.
\end{align*}
Now we can find an increasing sequence $(N_n)_{n\in \mathbb{N}}\subseteq \mathbb{N}$ with 
$N_n\rightarrow \infty$, such that the elementary approximation
\begin{align*}
\xi_n := \sum_{i=1}^{N_n}d_i^n\mathbf{1}_{D_i^n}(X)
\end{align*}
converges $\mathbb{P}$-a.s and in $L^2(\Omega ;H)$ to $X$. By Lemma \ref{RN-discrete} there are Bochner 
measurable functions $g_n:H\rightarrow L(H,U), g_n=\frac{dm_{\mu_n}}{d\mu_n}$, which only depend on
\begin{align*}
\mu_n = \mathcal{L}(\xi_n)&=\sum_{i=1}^{N_n}\mathcal{L}(X)(D_i^n)\delta_{d_i^n}
\end{align*}
and therefore only on $\mu=\mathcal{L}(X)$ and not on the random variable $X$ or the probability space 
$(\Omega,\mathcal{F},\mathbb{P})$, with 
\begin{align*}
\mathbb{E}\left[\mathbf{1}_Ag_n(\xi_n)u\right]&=m_{D\hat{f}(\xi_n)}(A)u 
= D\hat{f}(\xi_n)(u\mathbf{1}_A),
\end{align*}
for all $A\in \mathcal{F}$ and $g_n(\xi_n)\in L^2(\Omega ;L(H,U))$. Now we define 
\begin{align*}
\tilde{g}_n:=\sum_{i=1}^{N_n}g_n(d_i^n)\mathbf{1}_{D_i^n}.
\end{align*}
It then holds $\tilde{g}_n(X)=g_n(\xi_n)$, hence for any $u\in H$ and $A\in \mathcal{F}$
\begin{equation} 
\label{raddiff}
\begin{aligned}
\mathbb{E}\left[\mathbf{1}_A(\tilde{g}_n(X)-\tilde{g}_m(X))u\right]
& = \mathbb{E}\left[\mathbf{1}_A(g_n(\xi_n)-g_m(\xi_m))u\right] \\
& = D\hat{f}(\xi_n)(u\mathbf{1}_A)-D\hat{f}(\xi_m)(u\mathbf{1}_A)\\
& = (D\hat{f}(\xi_n)-D\hat{f}(\xi_m))(u\mathbf{1}_A).
\end{aligned}
\end{equation} 
Since $(D\hat{f}(\xi_n)-D\hat{f}(\xi_m))\in \Lambda_2^{\mathbb{P}}(H,U)$, using Corollary \ref{cor1}, we see 
that $m(A)u:= (D\hat{f}(\xi_n) -D\hat{f}(\xi_m))(u\mathbf{1}_A)$ defines a measure with finite 
$2$-variation 
\begin{align*}
\overline{m}_2(A) = \vertiii{(D\hat{f}(\xi_n)-D\hat{f}(\xi_m))_A}_{2,\mathbb{P}}.
\end{align*} 
Furthermore by \eqref{raddiff} and Lemma \ref{variationg} we have
\begin{align*}
|m|(A)= \mathbb{E}\left[ \mathbf{1}_A \|\tilde{g}_n(X)-\tilde{g}_m(X)\|_{op}\right].
\end{align*} 
Therefore by Lemma \ref{RN-iso} we have
\begin{align*}
\mathbb{E}\left[\|\tilde{g}_n(X)-\tilde{g}_m(X)\|_{op}^2\right]^{\frac 12} 
& = \overline{m}_2(\Omega) = \vertiii{D\hat{f}(\xi_n)-D\hat{f}(\xi_m)}_{2,\mathbb{P}}.
\end{align*}
Now by the continuity assumption on $D\hat{f}$ and since $(\Lambda_2^{\mathbb{P}}(H,U), 
\vertiii{\cdot }_2)$ is a Banach space, we can see that $(\tilde{g}_n(X))_{n\in \mathbb{N}}$ is a 
Cauchy-sequence in $L^2(\Omega ; L(H,U))$. Let 
\begin{align*}
\tilde{g}(X):=\lim\limits_{n\rightarrow \infty}\tilde{g}_n(X),
\end{align*}
in $L^2(\Omega ;L(H,U))$. Again due to the continuity assumption on $D\hat{f}$, it holds for any 
$A\in \mathcal{F}$ and $u\in H$
\begin{align*}
D\hat{f}(X)(u\mathbf{1}_A) 
& =\lim\limits_{n\rightarrow \infty}D\hat{f}(\xi_n)(u\mathbf{1}_A) 
= \lim\limits_{n\rightarrow \infty}\mathbb{E}\left[\mathbf{1}_A\tilde{g}_n(X)u\right] 
= \mathbb{E}\left[\mathbf{1}_A\tilde{g}(X)u\right].
\end{align*}
It is now easy to conclude 
\begin{align*}
D\hat{f}(X)Y=\mathbb{E}\left[\tilde{g}(X)Y\right],
\end{align*}
for any $Y\in L^2(\Omega ; H)$. Passing to a subsequence $\tilde{g}_{n_k}(X)$, converging almost surely, 
we define 
\begin{align*}
S := \{x\in H|\tilde{g}_{n_k}(x)\text{ converges in }L(H,U)\}.
\end{align*}
Now for $g:=\mathbf{1}_S\lim\limits_{k\rightarrow \infty}\tilde{g}_{n_k}$ we get 
\begin{align*}
\lim\limits_{k\rightarrow \infty}\tilde{g}_{n_k}(X)=g(X),\quad \mathbb{P}-a.s.,
\end{align*}
which yields $g(X)=\tilde{g}(X)$. 

Now indeed for all $Y\in L^2(\Omega ; H)$ we have 
\begin{align*}
D\hat{f}(X)Y 
& = \mathbb{E}\left[g(X)Y\right] = \mathbb{E}\left[g(X)\mathbb{E}\left[Y|X\right]\right] \\
& = D\hat{f}(X)\mathbb{E}\left[Y|X\right] = \int \mathbb{E}\left[Y|X\right]dm_{D\hat{f}(X)}.
\end{align*}
Furthermore since for all $B\in \mathcal{B} (H)$ and $u\in H$
\begin{align*}
m_{\mu}(B)u 
& = m_X(X\in B)u =\int u\mathbf{1}_{\{X\in B\}} dm_{D\hat{f}(X_0)} 
= \mathbb{E}\left[g(X)u\mathbf{1}_{B}(X)\right]\\
& = \int_{B} g(x)u\mu(dx),
\end{align*}
the function
\begin{align*}
\frac{dm_{\mu}}{d\mu}(x):=g(x),
\end{align*}
satisfies all properties stated in Theorem \ref{g} and is equal to the Radon-Nikodym derivative of $m_{\mu}$ 
with respect to $\mu$. By \cite[Theorem~3.1]{Ahm13} the Radon-Nikodym derivative is $\mu$-almost everywhere 
unique, hence $	\frac{dm_{\mu}}{d\mu}$ is $\mu$-almost everywhere uniquely defined. 
\end{proofThm}


\subsection{Structure of the L-Derivative}
	
Our main Theorem \ref{g} provides the factorized representation, necessary to give meaning to the 
$L$-derivative independent of the lifting. However, in addition to the usual $L$-differentiability 
from Definition \ref{Ldiff1} we need the continuity with respect to the $2$-variation norm. Therefore 
we will consider the following stronger definition of $L$-differentiability, based on our main result.

\begin{definition}
We say that a map $f:\mathcal{P}_2(H)\rightarrow U$ is absolutely continuous $L$-differentiable, if the lift 
\begin{align*}
\hat{f}:L^2(\Omega ;H)\rightarrow U,X\mapsto f(\mathcal{L}(X))
\end{align*}
is $\Lambda$-continuously Fr\'echet differentiable. We shall call the equivalence class $\partial_{\mu}f(\mu_0)$ 
given by Theorem \ref{g} the regular L-derivative of $f$ in $\mu_0$.
\end{definition}
	
\begin{remark}
If $U=\mathbb{R}$ every $L$-differentiable function is also absolutely continuous $L$-differentiable and the 
regular $L$-derivative coincides with the $L$-derivative in \cite{CGKPR22}. This follows immediately from Lemma 
\ref{variation1d}.
\end{remark}

\noindent 
Furthermore we will introduce the space of functions on the space of measures that depend on multiple variables 
and are two times differentiable, which we will need later on in Section \ref{Ito}.
	
\begin{definition}
We say that a function $f:[0,T]\times H\times \mathcal{P}_2(H)\rightarrow U$ is in $\mathcal{C}^{1,2,2}([0,T]
\times H\times\mathcal{P}_2(H),U)$, if the following holds:
\begin{enumerate}
\item For each fixed $x\in H,\mu\in \mathcal{P}_2(H)$, $f(\cdot,x,\mu)$ is continuously differentiable.
\item For each fixed $t\in [0,T],\mu\in \mathcal{P}_2(H)$, $f(t,\cdot,\mu)$ is two times continuously 
differentiable with jointly continuous derivatives 
\begin{align*}
(t,x,\mu)&\mapsto \partial_xf(t,x,\mu),\\
(t,x,\mu)&\mapsto \partial_{xx}f(t,x,\mu).
\end{align*}
\item For each fixed $t\in [0,T],x\in H$, $f(t,x,\cdot)$ is absolutely continuous $L$-differentiable and there 
exists a jointly continuous version of $(t,x,\mu,y)\mapsto \partial_{\mu}f(t,x,\mu)(y)$.
\item The jointly continuous version in 3. is differentiable in $y\in H$ in the classical sense, with jointly continuous derivative
\begin{align*}
(t,x,\mu,y)\in [0,T]\times \mathcal{P}_2(H)\times H&\mapsto \partial_y\partial_{\mu}f(t,x,\mu)(y)\in L(H,L(H,U)).\end{align*}
\end{enumerate}
If $f$ is not time dependent and satisfies 2. - 4. we say that $f$ is in $\mathcal{C}^{2,2} 
(H\times\mathcal{P}_2(H),U)$. Similarly, if $f$ does not depend on time and state and satisfies 3. - 4. we say that $f$ is in $\mathcal{C}^2(\mathcal{P}_2(H),U)$. If the corresponding derivatives are bounded we say that 
$f$ is in $\mathcal{C}_b^{1,2,2}([0,T]\times H\times\mathcal{P}_2(H),U)$.
\end{definition}

\noindent 	
In the rest of this section we will provide some important properties of the regular $L$-derivative. We start 
with the following remark on continuity/Lipschitz continuity of absolutely continuous $L$-differentiable maps.

\begin{remark}
Let $f:\mathcal{P}_2(H)\rightarrow U$ be absolutely continuous $L$-differentiable and $\mu,\mu_0\in 
\mathcal{P}_2(H)$. Then we have for any $X_0,X\in L^2(\Omega ;H)$ with $X_0\sim \mu_0,X\sim \mu$ 
\begin{align*}
f(\mu) & =f(\mu_0)+\mathbb{E}\left[\partial_{\mu}f(\mu_0)(X_0)(X-X_0)\right]+o(\|X-X_0\|_{L^2(\Omega ;H)}).
\end{align*}
Moreover we can find a $t\in [0,1]$, such that
\begin{align*}
\|f(\mu)-f(\mu_0)\|_U 
& = \|\mathbb{E}\left[\partial_{\mu}f(\mathcal{L}(tX+(1-t)X_0))(tX+(1-t)X_0) (X-X_0)\right]\|_U \\
& \le \mathbb{E}\left[\|\partial_{\mu}f(\mathcal{L}(tX+(1-t)X_0))(tX+(1-t)X_0)\|^2_{op} \right]^{\frac{1}{2}} 
  \mathbb{E}\left[\|X-X_0\|_H^2\right]^{\frac{1}{2}} \\
& = \vertiii{D\hat{f}(tX+(1-t)X_0)}_{2,\mathbb{P}} \|X-X_0\|_{L^2(\Omega ;H)}.
\end{align*}
If now $(\mu_n)_n\subset \mathcal{P}_2(H)$ is a sequence with 
\begin{align*}
\lim_{n\rightarrow \infty}W_2(\mu_n,\mu_0)=0,
\end{align*}
then by Skorokhod's theorem there exists a sequence of random variables $X_n\sim \mu_n,X\sim \mu$, 
such that $(X_n)_{n\in \mathbb{N}}$ converges almost surely towards $X$. Since $(\mu_n)_{n\in \mathbb{N}}$ 
is 2-uniformly integrable, see \cite[Theorem~5.5]{CD18}, i.e. for any $u_0\in H$ 
\begin{align*}
\lim_{r\rightarrow \infty}\sup_{n\geq 1}\int \|u_0-u\|_H^2\mathbf{1}_{\{\|u_0-u\|_H\geq r\}}d\mu_n(u)=0,
\end{align*}
we obtain
\begin{align*}
\lim_{n\rightarrow \infty}\|X_n-X\|_{L^2(\Omega ;H)}=0.
\end{align*}
Therefore, using the continuity of $D\hat{f}$, we obtain 
\begin{align*}
\lim_{n\rightarrow \infty}\|f(\mu_n)-f(\mu_0)\|_U=0,
\end{align*}
hence $f$ is continuous with respect to $W_2$. If $D\hat{f}$ is bounded in $\vertiii{\cdot}_{2,\mathbb{P}}$ by some 
constant $C>0$ we obtain
\begin{align*}
\|f(\mu)-f(\mu_0)\|_U&\le  C \|X-X_0\|_{L^2(\Omega, ;H)}.
\end{align*}
Thus taking the infimum over all $X,X_0$ we obtain 
\begin{align*}
\|f(\mu)-f(\mu_0)\|_U&\le C W_2(\mu,\mu_0).
\end{align*}
\end{remark}
	
\begin{lemma} 
\label{proj}
Let $U$ be a separable Hilbert space with ONB $(\tilde{\varphi}_i)_{i\in \mathbb{N}}$. If  
$f:\mathcal{P}_2(H)\rightarrow U$ is absolutely continuous $L$-differentiable, then for any $v\in U$
\begin{align*}
f_v:\mathcal{P}_2(H)\rightarrow \mathbb{R},\mu\mapsto \langle f(\mu),v\rangle_U
\end{align*}
is $L$-differentiable. Furthermore for every $\mu_0\in \mathcal{P}_2(H)$
\begin{align*}
\partial_{\mu}f_v(\mu_0)(y)=\partial_{\mu}f(\mu_0)(y)^{\ast}v
\end{align*} 
for all $v\in U$, $\mu_0$-almost everywhere, where $\partial_{\mu}f(\mu_0)(y)^{\ast}$ denotes the adjoint 
operator of $\partial_{\mu}f(\mu_0)(y)$. In particular, for $v=\tilde{\varphi}_i$ and $f_i (\mu )  
= \langle f(\mu),\tilde{\varphi}_i\rangle_U$ we have for $\mu_0$-almost every $y\in H$
\begin{align*}
\sum_{i=1}^{\infty}\langle \partial_{\mu}f_i(\mu_0)(y),w\rangle_H\tilde{\varphi}_i 
& = \sum_{i=1}^{\infty}\langle \tilde{\varphi}_i,\partial_{\mu}f(\mu_0)(y)w\rangle_U\tilde{\varphi}_i 
= \partial_{\mu}f(\mu_0)(y)w,
\end{align*}
i.e. the series converges for every $w\in H$, $\mu_0$-almost everywhere. Furthermore if $f\in  
\mathcal{C}^{2}(\mathcal{P}_2(H),U)$, then for any $v\in U$ $f_v\in \mathcal{C}^{2}(\mathcal{P}_2(H), 
\mathbb{R})$ and we have 
\begin{align*}
\partial_y\partial_{\mu}f_v(\mu)(y)z=(\partial_y\partial_{\mu}f(\mu_0)(y)z)^{\ast}v.
\end{align*}
\end{lemma}
	
\begin{proof}
Since the lift 
\begin{align*}	
\hat{f}:L^2(\Omega ;H)\rightarrow U, X\mapsto f(\mathcal{L}(X))
\end{align*}
is $\Lambda$-continuously Fr\'echet differentiable, it is immediate that
\begin{align*}
\hat{f}_v:L^2(\Omega ;H)\rightarrow \mathbb{R},X\mapsto \langle f(\mathcal{L}(X)),v\rangle_U
\end{align*}
is continuously Fr\'echet differentiable with
\begin{align*}
D\hat{f}_v(X)=D\hat{f}(X)^{\ast}v\in L^2(\Omega ;H),
\end{align*}
where we identified $D\hat{f}_v(X)\in L(L^2(\Omega ; H),\mathbb{R})$ with its dual element in $L^2(\Omega ;H)$. 
Now for all $\mu_0\in \mathcal{P}_2(H)$ and $X_0\in L^2(\Omega ;H)$ with $\mathcal{L}(X_0)=\mu_0$ we have 
\begin{align*}
\mathbb{E}\left[\langle D\hat{f}_v(X_0),Y\rangle_H\right] 
& = \mathbb{E}\left[\langle D\hat{f}(X_0)^{\ast}v,Y\rangle_H\right] \\
& = \langle D\hat{f}(X_0)Y,v\rangle_U \\
& = \langle \mathbb{E}\left[\partial_{\mu}f(\mu_0)(X_0)Y\right],v\rangle_U \\
& = \mathbb{E}\left[\langle \partial_{\mu}f(\mu_0)(X_0)Y,v\rangle_U\right] \\
& = \mathbb{E}\left[\langle Y,\partial_{\mu}f(\mu_0)(X_0)^{\ast}v\rangle_H\right],
\end{align*}
for all  $Y\in L^2(\Omega ;H)$. This yields 
\begin{align*}
D\hat{f}_v(X_0)=\partial_{\mu}f(\mu_0)(X_0)^{\ast}v
\end{align*}
$\mathbb{P}$-a.s., and since the element $\partial_{\mu}f_v(\mu_0)$ is uniquely defined in 
$L^2(H; H, \mu_0)$ we obtain 
\begin{align*}
\partial_{\mu}f_v(\mu_0)(x)=\partial_{\mu}f(\mu_0)(x)^{\ast}v,
\end{align*} 
for $\mu_0$-almost every $x$. Now let $x,y,z\in H$, then we have 
\begin{align*}
& \|\partial_{\mu}f_v(\mu_0)(y+z)-\partial_{\mu}f_v(\mu_0)(y)-(\partial_y\partial_{\mu} 
f(\mu_0)(y)z)^{\ast}v\|_H \\
& = \|(\partial_{\mu}f(\mu_0)(y+z)-\partial_{\mu}f(\mu_0)(y)-\partial_y\partial_{\mu} 
f(\mu_0)(y)z)^{\ast}v\|_H \\
& \le \|\partial_{\mu}f(\mu_0)(y+z)-\partial_{\mu}f(\mu_0)(y)-\partial_y\partial_{\mu}f(\mu_0)(y)z\|_{op} 
\|v\|_U \in o(\|y\|_H).
\end{align*}   
Therefore
\begin{align*}
\partial_y\partial_{\mu}f_v(\mu)(y)z=(\partial_y\partial_{\mu}f(\mu_0)(y)z)^{\ast}v.
\end{align*}
\end{proof}

\begin{remark}
As already mentioned in the introduction, if the lift of $f:\mathcal{P}_2(H)\rightarrow U$ is Fr\'echet 
differentiable, but not $\Lambda$-continuously Fr\'echet differentiable, then 
\begin{equation} 
\label{seriesdiv}
\sum_{i=1}^{\infty}\langle \partial_{\mu}f_i(\mu_0)(y),u\rangle_H\tilde{\varphi}_i
\end{equation}
will not necessary converge for all $u\in H$, $\mu_0$-almost everywhere. In particular, not every 
$L$-differentiable map is also absolutely continuous $L$-differentiable.
\end{remark}
	
\begin{lemma}
Let $U$ be a separable Hilbert space with ONB $(\tilde{\varphi}_i)_{i\in \mathbb{N}}$. If $f:\mathcal{P}_2(H)
\rightarrow U$ is absolutely continuous $L$-differentiable, then
\begin{align*}
\partial_{\mu}f(\mu_0)(y)^{\ast}v=\mathbb{E}\left[D\hat{f}(X_0)^{\ast}v|X_0=y\right].
\end{align*} 
$\mu_0$-almost everywhere, i.e. for all $u\in H$ and $v\in U$
\begin{align*}
\langle \partial_{\mu}f(\mu_0)(y)u,v\rangle_U=\langle u,\mathbb{E}\left[D\hat{f}(X_0)^{\ast}v|X_0=y\right]
\rangle_H,
\end{align*}
$\mu_0$-almost everywhere.
\end{lemma}
	
\begin{proof}
Let $B\in \mathcal{B}(H)$, then by the defining property of the regular $L$-derivative
\begin{align*}
\mathbb{E}\left[\langle \partial_{\mu}f(\mu_0)(X_0)u,v\rangle_U\mathbf{1}_B(X_0)\right] 
& = \mathbb{E} \left[\langle \partial_{\mu}f(\mu_0)(X_0)u\mathbf{1}_B(X_0),v\rangle_U\right] \\
& = \langle \mathbb{E}\left[\partial_{\mu}f(\mu_0)(X_0)u\mathbf{1}_B(X_0)\right],v\rangle_U \\
& = \langle D\hat{f}(X_0)(u\mathbf{1}_B(X_0)),v\rangle_U \\
& = \mathbb{E}\left[\langle u,D\hat{f}(X_0)^{\ast}v\rangle_H\mathbf{1}_B(X_0)\right].
\end{align*}
Thus 
\begin{align*}
\langle\partial_{\mu}f(\mu_0)(X_0)u,v\rangle_U 
= \mathbb{E}\left[\langle u,D\hat{f}(X_0)^{\ast}v\rangle_H|X_0\right]
\end{align*}
almost everywhere, which concludes the proof. 
\end{proof}

\begin{remark}
Again it is not clear if 
\begin{align*}
\sum_{i=1}^{\infty}\mathbb{E}\left[\langle u,D\hat{f}(X_0)^{\ast}\tilde{\varphi}_i\rangle_U|X_0=y\right] 
 \tilde{\varphi}_i
\end{align*}
will converge, if the lift of $f$ is not $\Lambda$-continuously Fr\'echet differentiable. In particular, it is 
not clear if $y\mapsto U(y)$, defined by 
\begin{align*}
\langle U(y)u,v\rangle_U :=\langle u,\mathbb{E}\left[D\hat{f}(X_0)^{\ast}v|X_0= y \right]\rangle_H , 
\end{align*}
is a map from $H$ to $L(H,U') \cong L(H,U)$.
\end{remark}
	
\subsection{Examples}
\label{Examples1}

As a prototypical example we illustrate the structure of the regular $L$-derivative in the linear case. 
This is the infinite dimensional analogue to \cite[Section 5.2.2 Example~1]{CD18}. We then present an example 
of a map which is $L$-differentiable but not absolutely continuous $L$-differentiable.
	
\begin{lemma}
Consider a Fr\'echet differentiable map $h:H\rightarrow H$, whose Fr\'echet derivative is a continuous map 
$Dh:H\rightarrow L(H)$ which satisfies the linear growth condition 
\begin{align*}
\|Dh(x)\|_{op} \le C(1+\|x\|_H),
\end{align*} 
for all $x\in H$. Let
\begin{align*}
f:\mathcal{P}_2(H)\rightarrow H, f(\mu):=\int_{H}^{}h(x)\mu(dx),
\end{align*}
where the above integral is the Bochner integral. Then for every $X_1,X_2\in L^2(\Omega ;H)$
\begin{align*}
& \vertiii{D\hat{f}(X_1)}_2^2\le \mathbb{E}\left[\|Dh(X_1)\|_{op}^2\right]<\infty, \\
& \vertiii{D\hat{f}(X_1)-D\hat{f}(X_2)}_2^2\le \mathbb{E}\left[\|Dh(X_1)-Dh(X_2)\|_{op}^2\right].
\end{align*}
In particular, $f$ is absolutely continuous L-differentiable with 
\begin{align*}
\partial_{\mu}f(\mu_0)(x)=Dh(x).
\end{align*}
\end{lemma}

\begin{proof} 
First we notice that since $D(h)$ is linearly bounded, it follows that $h$ is at most of quadratic growth, 
i.e. there exists a finite constant $C$ with $\|h(x)\|_H\le C(1 + \|x\|_H^2)$ for all $x\in H$. Therefore 
all appearing Bochner integrals are well defined. Indeed, since $H$ is separable, the 
map $h:H\rightarrow H$ is (strongly) measurable by Pettis' measurability theorem, and $\int\|x\|^2_H 
\mu (dx) < \infty$ for $\mu\in\mathcal{P}_2 (H)$. Similarly, $Dh(X_0)$ is Bochner integrable with respect to 
$(L(H),\|\cdot\|_{op})$. This follows directly from the continuity of $Dh$, hence the separability of 
$Dh(H)\subset L(H)$, and the linear growth assumption on $Dh$. Now thanks to the Fr\'echet differentiability 
of $h$ we have for all $X,Y\in L^2(\Omega ;H)$
\begin{align*}
\hat{f}(X+Y) 
& =\hat{f}(X)+\mathbb{E}\left[Dh(X)Y\right]+\mathbb{E}\left[\int_{0}^{1} 
(Dh(X+\lambda Y)-Dh(X))Yd\lambda\right].
\end{align*}
Furthermore
$$ 
\|\mathbb{E}\left[\int_{0}^{1}(Dh(X+\lambda Y)-Dh(X))Yd\lambda\right]\|_H 
\leq \mathbb{E}\left[\sup_{\lambda\in [0,1]}\|Dh(X+ \lambda Y)-Dh(X)\|_{op}	\|Y\|_H \right] .
$$
To obtain Fr\'echet differentiability of $\hat{f}$ at $X$ with Fr\'echet differential 
$D\hat{f} (X)Y = \mathbb{E}[Dh (X)Y]$ it suffices now to show that 
\begin{equation}
\label{Convergence} 
\mathbb{E}\left[\sup_{\lambda\in [0,1]} \|Dh(X+ \lambda Y)-Dh(X)\|_{op}	\|Y\|_H \right] 
= o( \|Y\|_{L^2(\Omega ;H)}). 
\end{equation}
To this end define the set $A := \{\|Y\|_H \ge \|Y\|^{1/2}_{L^2 (\Omega ; H)}\}$. Then 
\begin{equation} 
\label{Convergence1}
\begin{aligned}
\mathbb{E} & \big[1_A \sup_{\lambda\in [0,1]} \|Dh(X+ \lambda Y)-Dh(X)\|_{op}	\|Y\|_H \big] \\
& \le 2C \mathbb{E}\left[1_A (1+\|X\|_H + \|Y\|_H )\|Y\|_H \right] \\ 
& \le 2C \left( \mathbb{E}\left[ 1_A (1+\|X\|_H)^2\right]^{1/2} + \|Y\|_{L^2 (\Omega ; H)} \right) \|Y\|_{L^2 
(\Omega ; H)} = o (\|Y\|_{L^2 (\Omega ; H)})\, , 
\end{aligned} 
\end{equation} 
since $\mathbb{E}\left[ 1_A (1+\|X\|_H)^2\right] \to 0$ as $\|Y\|_{L^2 (\Omega ; H)} \to 0$ by Lebesgue's 
theorem, because 
$$ 
P[A] \le \mathbb{E}\left[\frac{\|Y\|_H^2}{\|Y\|_{L^2 (\Omega ; H)}}\right]  = \|Y\|_{L^2 (\Omega ; H)} \to 0, 
$$ 
thereby using Markov's inequality. On the complement $A^c$ we apply Cauchy-Schwartz inequality to obtain that  
\begin{equation} 
\label{Convergence2}
\begin{aligned}
\mathbb{E} & \Big[1_{A^c} \sup_{\lambda\in [0,1]} \|Dh(X + \lambda Y)-Dh(X)\|_{op}	\|Y\|_H \Big] \\
& \le \mathbb{E}\Big[1_{A^c} \sup_{y\in H:\|y\|_H\le \|Y\|^{1/2}_{L^2(\Omega ;H)}} 
 \|Dh(X+ y)-Dh(X)\|_{op}^2	\Big]^{1/2} \|Y\|_{L^2 (\Omega ; H)} = o (\|Y\|_{L^2 (\Omega ; H)})\, , 
\end{aligned} 
\end{equation} 
since 
$$
\mathbb{E}\Big[1_{A^c} \sup_{y\in H:\|y\|_H\le \|Y\|^{1/2}_{L^2(\Omega ;H)}} \|Dh(X+ y)-Dh(X)\|_{op}^2	
\Big] \to 0 , 
$$
as $\|Y\|_{L^2(\Omega ;H)}\to 0$, by continuity of $Dh$ and dominated convergence. 
\eqref{Convergence1} and \eqref{Convergence2} together imply \eqref{Convergence} and therefore 
$D\hat{f}(X)Y=\mathbb{E}\left[Dh(X)Y\right]$, hence $\vertiii{D\hat{f}(X)}_2 
= \mathbb{E}\left[\|Dh(X)\|_{op}^2\right]^{\frac{1}{2}}<\infty$, thanks to the linear growth assumption 
on $Dh$ and the square integrability of $X$.

The necessary condition $\vertiii{D\hat{f}(X)}_2 < \infty$ can also be derived directly, since for 
$Y = \sum_{i=1}^{n} \mathbf{1}_{A_i}x_i \in\mathcal{E}_H (\Omega )$, with $x_i\in H$ and $A_i$ disjoint,  
$$ 
\begin{aligned} 
\sum_{i=1}^{n}\|D\hat{f}(X)(\mathbf{1}_{A_i}x_i)\|_H
& = \sum_{i=1}^{n} \|\mathbb{E} \left[ \mathbf{1}_{A_i} Dh(X)x_i\right] \|_H  
\le  \mathbb{E}\left[\sum_{i=1}^{n}\mathbf{1}_{A_i}\|Dh(X)\|_{op}\|x_i\|_H\right] \\ 
& \le \mathbb{E}\left[\|Dh(X)\|_{op}^2\right]^{\frac{1}{2}}\mathbb{E}\left[\|Y\|_H^2 \right]^{\frac{1}{2}} 
\, , 
\end{aligned} 
$$ 
so that 
\begin{align*}
\vertiii{D\hat{f}(X)}_2 
& = \sup\{\sum_{i=1}^{n}\|D\hat{f}(X)(\mathbf{1}_{A_i}x_i)\|_H : Y = \sum_{i=1}^{n}\mathbf{1}_{A_i}x_i 
	\in \mathcal{E}_H(\Omega),A_i\text{ disjoint},\mathbb{E}\left[\|Y\|_H^2\right]\le 1\} \\
& \le \mathbb{E}\left[\|Dh(X)\|_{op}^2\right]^{\frac{1}{2}}.
\end{align*}
Similarly, the continuity of $D\hat{f}$ with respect to $\vertiii{\cdot}_2$ follows by taking 
$D\hat{f}(X_1)- D\hat{f}(X_2)$ instead of $ D\hat{f}(X)$ in the above considerations.
\end{proof}

\noindent 	
We will now construct an explicit example of a map $f:\mathcal{P}_2(H)\rightarrow U$ which is 
$L$-differentiable, i.e. its lift $\hat{f}$ is Fr\'echet differentiable, but not absolutely 
continuous $L$-differentiable. Our counterexample will include the following function as a building block: 
\begin{equation} 
\label{defH} 
H (\mu ) := \int_\mathbb{R} \frac{e^{-\frac{x^2}{2}}}{\sqrt{2\pi}}\left( \int_x^\infty 
\frac{e^{-\frac{y^2}{2}}}{\sqrt{2\pi}} (F_\mu (y) - \Phi (y))\, dy \right)^2 \, dx . 
\end{equation} 
Here $\Phi$ denotes the cumulative distribution function of the standard normal distribution 
$\mathcal{N} (0,1)$. The following Proposition collects some useful properties of $H$. 

\begin{proposition} 
\label{PropertiesH}
Let $H$ be as defined in \eqref{defH}. Then: 
\begin{itemize} 
\item[(i)] The lift $\hat{H} (X) = H(\mathcal{L} (X))$ can be represented as 
\begin{align*}
H(\mathcal{L} (X)) = \mathbb{E} \left[ H_0 (X,\tilde{X})\right] 
\end{align*} 
with 
\begin{equation} 
\label{KernelH}
\begin{aligned}
H_0 (X,\tilde{X}) 
& = - \frac{1}{12} \Phi^4 (X) - \frac{1}{2} \Phi^2 (X) - \frac{1}{3} \Phi^3 (X\vee \tilde{X}) + \Phi^2 (X) \Phi (X\vee\tilde{X}) 
 + \frac{1}{20} , 
\end{aligned}
\end{equation} 
where $\tilde{X}$ is an independent copy of $X$. 
\item[(ii)] $H$ is $L$-differentiable with
\begin{equation} 
\label{DerivativeH}
\begin{aligned}
\partial_{\mu} H (\mu)(x) 
& =  \Big( -\frac{1}{3} \Phi (x)^3 - \Phi (x) + \Phi^2 (x) F_\mu (x) \\  
& \qquad + 2 \Phi (x) \int_{x+}^\infty \Phi (y) \, \mu (dy) 
+ \int_{-\infty}^x \Phi^2 (y) \, \mu (dy) \Big) \frac{e^{-x^2/2}}{\sqrt{2\pi}} .
\end{aligned}
\end{equation} 
Here we set $\int_{x+}^\infty \Phi (x) \mu (dy) := \int 1_{(x, \infty)} (y) \Phi (y)\mu (dy)$. 

\item[(iii)] If $X\sim\mathcal{N}(0,1)$, or equivalently, $H(\mathcal{L}(X)) = 0$, 
then 
\begin{equation} 
\label{LipschitzAtNormalDistrH}
|H(\mathcal{L} (X+Y)) - H(\mathcal{L}(X))|\le L_H \|Y\|_{L^2 (\Omega ; \mathbb{R})}^2. 
\end{equation} 
\end{itemize} 
\end{proposition}
	
The proof of the Proposition is given in the Appendix \ref{AppendixA}. 

\begin{theorem} 
\label{mainEx}
Let 
\begin{align*}
f:\mathcal{P}_2(\mathbb{R})\rightarrow L^2(\mathbb{R}),f(\mu)(t) 
:= \int g_0 (t-\omega,H(\mu ))F_{\mu}(\omega)d\omega ,  
\end{align*}
where 
\begin{align*}
g_0 (x,y) := -\partial_x G_0 (x,y) = \frac{\text{ sign }(x)}{\sqrt{|x| + y}} \mathbf{1}_{\{|x|\le 1\}} ,
\end{align*}
with $G_0 : \mathbb{R} \times [0, \infty ) \rightarrow \mathbb{R}$,  
$G_0 (x,y) := 2 \left( \sqrt{1+y} - \sqrt{|x|+y} \right) 1_{\{|x|\le 1\}}$. 

\smallskip 
\noindent 
Then the lift $\hat{f} (\mathcal{L} (X))$ of $f$ is Fr\'echet differentiable with 
\begin{equation}
\label{RepresentationDifferential1} 
\begin{aligned} 
D\hat{f}(X)(Y) (t) 
& = -\mathbb{E} [g_0(t-X, H(\mathcal{L}(X))Y] \\ 
& \qquad 
+ \mathbb{E} [ \partial_y G_0 (t-X,H(\mathcal{L} (X)))]\, \mathbb{E} [\partial_\mu H(\mathcal{L} (X))(X)(Y)] .
\end{aligned} 
\end{equation}
In particular, $f$ is $L$-differentiable with 
\begin{equation}
\label{RepresentationLDerivative}  
\partial_\mu f(\mu )(x)(t) = - g_0 (t-x, H(\mu ))  
- \int \partial_y G_0 (t-\tilde{x}, H(\mu )) \mu (d\tilde{x}) \partial_\mu H(\mu )(x) . 
\end{equation} 
Moreover, in $\mu = \mathcal{N}(0,1)$ 
\begin{equation}
\label{RepresentationLDerivativeNormalDistribution}  
\partial_\mu f( \mathcal{N}(0,1))(x)(t) = - g_0 (t-x, 0) = -\frac{sign (t-x)}{\sqrt{|t-x|}} 
1_{\{|t-x|\le 1\}} \notin L^2 (\mathbb{R}) ,  
\end{equation} 
hence $f$ is not absolute continuous $L$-differentiable. 
\end{theorem}

The proof of the Theorem is given in the Appendix \ref{AppendixB}.

\section{Higher Order Expansion for Mean-Field SPDEs}
\label{Ito}

In this section we will develop a higher order Taylor expansion for Mean-Field SPDEs. Our expansion 
will be based on a mild Itô formula for flows of measures generated by Mean-Field SPDEs. A similar 
mild Itô formula for SPDEs can be found in \cite{DPJR19}.

\noindent 	
In the following we fix $T>0$ and consider a filtered probability space $(\Omega,\mathcal{F}, 
(\mathcal{F}_t)_{t\in [0,T]},\mathbb{P})$ and separable Hilbert spaces $H,U$ with given ONB  
$(\varphi_n)_{n\in \mathbb{N}}, (\tilde{\varphi}_n)_{n\in \mathbb{N}}$ respectively. For some 
$f\in\mathcal{C}^{1,2}([0,T]\times\mathcal{P}_2(H),U)$, the map $\partial_y\partial_{\mu} 
f(t,\mu)(y)\in L(H,L(H,U))\cong L(H\times H,U)$ is bilinear. Given a linear operator $B\in L(H,H)$, we can 
consider the composition $\partial_x\partial_{\mu}f(t,\mu)(y) (Bw_1)(Bw_2)$, $w_1, w_2\in L(H, L(H,U))$. 
If $B$ is Hilbert-Schmidt, this composition admits a well-defined trace, defined by    
$$ 
\begin{aligned} 
\textbf{tr}_H (\partial_x\partial_{\mu}f(t,\mu)(y)B^\ast B) 
& := \sum_{k=1}^\infty \partial_x\partial_{\mu}f(t,\mu)(y) (B\varphi_k) (B\varphi_k) ,    
\end{aligned} 
$$ 
and taking values in $U$, since 
$$
\begin{aligned} 
\left\|\sum_{k=1}^\infty \partial_x\partial_{\mu}f(t,\mu)(y) (B\varphi_k) (B\varphi_k) \right\|_U 
& \le \sum_{k=1}^\infty \| \partial_x\partial_{\mu}f(t,\mu)(y)(B\varphi_k) (B\varphi_k) \|_U \\
& \le \sum_{k=1}^\infty \| \partial_x\partial_{\mu}f(t,\mu)(y) (B\varphi_k) \|_{L(H,U)} \|B\varphi_k\|_H \\
& \le \sum_{k=1}^\infty \| \partial_x\partial_{\mu} f(t,\mu)(y)\|_{L(H, L(H,U))} \|B\varphi_k\|_H^2  \\
& = \| \partial_x\partial_{\mu} f(t,\mu)(y)\|_{L(H, L(H,U))} \|B\|_{L_2 (H,H)}^2 < \infty  . 
\end{aligned} 
$$

\subsection{Itô-Formula for Flow of Measures in Infinite Dimensions}

We start with a simple Ito formula for semi-martingale flows of measures in infinite dimensions, which is 
similar to \cite[Theorem~4.25]{CGKPR22}.
	
\begin{theorem}
\label{Ito1}
Consider the $H$-valued Itô process
\begin{align*}
du_t & = b_t dt + \sigma_t dW_t,\quad t\in [s,T]\\
u_s  & \in L^2(\Omega ;H),
\end{align*}
where $W=(W_t)$ is a cylindrical $Q$-Wiener process on a separable Hilbert space $\Xi$ and the covariance 
operator $Q:\Xi\rightarrow \Xi$ is assumed to be linear, bounded, positive-definite and self-adjoint. Let 
$b,\sigma$ be $(\mathcal{F}_t)$-progressively measurable with values in $H$ and $L_2(\Xi,H)$ respectively. 
Furthermore we assume that 
\begin{align*}
\int_0^T\mathbb{E}\left[\|b_t\|_H^2\right]dt < \infty,\quad \int_0^T 
\mathbb{E}\left[\|\sigma_t\|_{L_2(\Xi,H)}^2\right]dt < \infty.
\end{align*}
If $f:[0,T]\times\mathcal{P}_2(H)\rightarrow U$ is in $\mathcal{C}^{1,2}([0,T]\times \mathcal{P}_2(H),U)$, 
then we have with $\mu_t:=\mathcal{L}(u_t)$, $a_t:=\sigma_t\sigma_t^{\ast}$
\begin{align*}
f(t,\mu_t) 
& = f(s,\mu_s)+\int_{s}^t\partial_rf(r,\mu_r)dr+\int_s^t \mathbb{E}\left[\partial_{\mu} 
	f(r,\mu_r)(u_r)b_r\right] dr \\
&\quad + \frac{1}{2}\int_s^t\mathbb{E}\left[\textbf{tr}_H(\partial_x\partial_{\mu}f(r,\mu_r)(u_r) 
	\sigma_r Q\sigma_r^{\ast})\right]dr,
\end{align*}
for all $t\in [s,T]$.
\end{theorem}

\begin{proof}
Without loss of generality we assume that $Q=Id$. Let $f:[0,T]\times\mathcal{P}_2(H)\rightarrow U$ be in 
$\mathcal{C}^{1,2}([0,T]\times \mathcal{P}_2(H),U)$, then by Lemma \ref{proj}
\begin{align*}
f_v:[0,T]\times\mathcal{P}_2(H)\rightarrow \mathbb{R},f_v(t,\mu):=\langle f(t,\mu),v\rangle_U
\end{align*}
is in $\mathcal{C}^{1,2}([0,T]\times \mathcal{P}_2(H),\mathbb{R})$ for any $v\in U$ and 
\begin{equation}\label{der1}
\begin{aligned}
\partial_{\mu}f_v(t,\mu)(x) 
& =\partial_{\mu}f(t,\mu_0)(x)^{\ast}v, \\ 
\partial_x\partial_{\mu}f_v(t,\mu)(x)w 
& = (\partial_x\partial_{\mu}f(t,\mu)(x)w)^{\ast}v \quad\forall w\in H.
\end{aligned}
\end{equation}
Now the Ito formula in \cite[Theorem~4.25]{CGKPR22} yields 
\begin{align*}
f_v(t,\mu) 
& = f_v(s,\mu_s)+\int_s^t\partial_t f_v(r,\mu_r)dr+\int_s^t\mathbb{E}\left[\langle b_r,\partial_{\mu} 
	f_v(r,\mu_r)(u_r)\rangle_H\right]dr \\ 
& \quad + \frac{1}{2} \int_s^t \mathbb{E}\left[\mathbf{tr}_H (\sigma_r\sigma_r^{\ast}\partial_x\partial_{\mu} 
	f_v(r,\mu_r)(u_r))\right]dr.
\end{align*}
By \eqref{der1} we have 
\begin{align*}
\mathbb{E}\left[\langle b_r,\partial_{\mu}f_v(r,\mu_r)(u_r)\rangle_H\right] 
& = \mathbb{E}\left[\langle b_r, (\partial_{\mu}f(r,\mu_r)(u_r)^\ast v\rangle_H\right] 
= \mathbb{E}\left[\langle \partial_{\mu}f(r,\mu_r)(u_r)b_r,v\rangle_U\right]. 
\end{align*}
Now for any fixed $r\in [s,T]$ we have $a_r = \sigma_r\sigma_r^{\ast}$ is a positive semidefinite 
symmetric trace-class operator, hence we can find an orthonormal basis $(\varphi_k)_k$ of $H$ 
consisting of eigenvectors of $a_r$ with corresponding (nonnegative) eigenvalues $(\lambda_k)$. 
Using \eqref{der1} again, we have that 
$$
\begin{aligned} 
\langle\varphi_k, a_r\partial_x\partial_{\mu} f_v(r,\mu_r)(u_r)\varphi_k\rangle_H
& = \langle a_r\varphi_k, \partial_x\partial_{\mu} f_v(r,\mu_r)(u_r)\varphi_k\rangle_H \\
& = \langle\sqrt{\lambda_k} \varphi_k, \partial_x\partial_{\mu} f_v(r,\mu_r)(u_r) (\sqrt{\lambda_k} 
\varphi_k )\rangle_H \\ 
& = \langle \sqrt{\lambda_k} \varphi_k, \left( \partial_x\partial_{\mu} f (r,\mu_r)(u_r) (\sqrt{\lambda_k} 
\varphi_k )\right)^* v\rangle_H \\ 
& = \langle \partial_x\partial_{\mu} f (r,\mu_r)(u_r) (\sqrt{\lambda_k} \varphi_k ) 
(\sqrt{\lambda_k} \varphi_k), v\rangle_U \, . 
\end{aligned} 
$$
Summing up w.r.t. $k$ and taking expectation we obtain that 
\begin{align*}
\mathbb{E}\left[\mathbf{tr}_H (a_r\partial_x\partial_{\mu} f_v(r,\mu_r)(u_r))\right] 
& = \mathbb{E}\left[\langle \mathbf{tr}_H (\partial_x\partial_{\mu}f(r,\mu_r)(u_r) a_r ,v\rangle_U \right].
\end{align*} 
This concludes the proof.
\end{proof}

\subsection{A Mild Itô-Formula for Mean-Field SPDEs}
	
Given $\phi\in L^2(\Omega ;H)$ we consider the Mean-Field SPDE
\begin{equation} 
\label{MFSPDE}
\begin{aligned}
du_t & = \left[Au_t+b(u_t,\mathcal{L}(u_t))\right]dt+\sigma(u_t,\mathcal{L}(u_t))dW_t,\quad t\in [t_0,T] \\
u_{t_0} & = \phi\in L^2(\Omega ;H),
\end{aligned}
\end{equation}
where $W:[0,T]\times \Omega\rightarrow \Xi$ denotes a cylindrical $Q$-Wiener process on a separable Hilbert 
space $\Xi$. The covariance operator $Q:\Xi\rightarrow \Xi$ is assumed to be linear, bounded, positive-definite 
and self-adjoint, $(b,\sigma):H\times \mathcal{P}_2(H)\rightarrow  (H,L_2(\Xi ;H))$ are measurable and  
$A:D(A)\subset H\rightarrow H$ is a densely defined self-adjoint, negative definite linear operator with domain 
$D(A)$ and compact inverse (for example the Dirichlet-Laplace on a bounded domain). Thus $A$ is the 
generator of an analytic semigroup of contractions, in the following denoted by $(e^{tA})_{t\ge 0}$. 
Furthermore we assume the following.
	
\begin{enumerate}[start=0,label={(\bfseries A\arabic*):}]
\item The drift $b:H\times \mathcal{P}_2(H)\rightarrow H$ is globally Lipschitz, i.e. there exists some 
constant $C>0$ with
\begin{align*}
\|b(u_1,\mu_1)-b(u_2,\mu_2)\|_H&\le C (\|u_1-u_2\|_H+W_2(\mu_1,\mu_2)),
\end{align*} 
for all $u_1,u_2\in H$ and $\mu_1,\mu_2\in \mathcal{P}_2(H)$.
\item The noise coefficient $\sigma:H\times \mathcal{P}_2(H)\rightarrow L_2(\Xi,H)$ is globally Lipschitz, 
i.e. there exists some constant $C>0$ with
\begin{align*}
\|\sigma(u_1,\mu_1)-\sigma(u_2,\mu_2)\|_{L_2(\Xi,H)}&\le C (\|u_1-u_2\|_H+W_2(\mu_1,\mu_2)),
\end{align*}
for all $u_1,u_2\in H$ and $\mu_1,\mu_2\in \mathcal{P}_2(H)$.
\item There exists some constant $C>0$, such that for all $u\in H$ 
\begin{align*}
& \|b(u,\delta_0)\|_H\le C,\\
& \|\sigma(u,\delta_0)\|_{L_2(\Xi,H)}\le C,
\end{align*}
where $\delta_0$ denotes the dirac measure in $0\in H$.
\end{enumerate} 

\begin{definition}
A mild solution to \eqref{MFSPDE} is a continuous progressively measurable process $u=(u_t)_{t\in [t_0,T]}$ 
satisfying 
\begin{equation} 
\label{mild}
u_t = e^{(t-t_0)A}u_{t_0}+\int_{t_0}^te^{(t-s)A}b(u_s,\mathcal{L}(u_s))ds+\int_{t_0}^te^{(t-s)A}\sigma(u_s, 
\mathcal{L}(u_s))dW_s,\quad \forall t\in [t_0,T],\mathbb{P}-a.s.
\end{equation}
\end{definition}

\begin{theorem} 
\label{Ito2}
Under the assumptions (A0)-(A2) there exists a unique mild solution $u$ to \eqref{MFSPDE} with 
\begin{align*}
\mathbb{E}\left[\sup_{t\in [t_0,T]}\|u_t\|_H^2\right]<\infty.
\end{align*} 
For $f\in \mathcal{C}^{1,2,2}([0,T]\times H\times \mathcal{P}_2(H),U)$ we then have 
\begin{align*}
& f(t,u_t,\mathcal{L}(u_t))=f(t_0,e^{(t-t_0)A}u_{t_0},\mathcal{L}(e^{(t-t_0)A}u_{t_0})) 
	+ \int_{t_0}^t\partial_sf(s,e^{(t-s)A}u_s,\mathcal{L}(e^{(t-s)A}u_s))ds \\
& \quad + \int_{t_0}^t \tilde{\mathbb{E}}\left[\partial_{\mu}f(s,e^{(t-s)A}u_s,\mathcal{L}(e^{(t-s)A}u_s)) 
	(e^{(t-s)A}\tilde{u}_s)e^{(t-s)A}b(\tilde{u}_s,\mathcal{L}(u_s))\right]ds \\
& \quad + \frac{1}{2}\int_{t_0}^t \tilde{\mathbb{E}}\left[\textbf{tr}_H(\partial_y\partial_{\mu} f(s,u_s,
	\mathcal{L}(e^{(t-s)A}u_s))(e^{(t-s)A}\tilde{u}_s)e^{(t-s)A}\sigma(\tilde{u}_s, \mathcal{L}(u_s)) 
	Qe^{(t-s)A}\sigma(\tilde{u}_s, \mathcal{L}(u_s))^{\ast})\right]ds  \\
& \quad + \int_{t_0}^t\partial_x f(s,e^{(t-s)A}u_s,\mathcal{L}(e^{(t-s)A}u_s))e^{(t-s)A}b(u_s, 
	\mathcal{L}(u_s))ds \\
& \quad + \int_{t_0}^t\partial_x f(s,e^{(t-s)A}u_s,\mathcal{L}(e^{(t-s)A}u_s))e^{(t-s)A}\sigma(u_s, 
	\mathcal{L}(u_s))dW_s \\
& \quad + \frac{1}{2} \int_{t_0}^t \textbf{tr}_H(\partial_{xx}f(s,e^{(t-s)A}u_s, \mathcal{L}(e^{(t-s)A}
	u_s))e^{(t-s)A}\sigma(u_s,\mathcal{L}(u_s))Qe^{(t-s)A}\sigma(u_s,\mathcal{L}(u_s))^{\ast})ds.
\end{align*}
Here, $\tilde{u}$ denotes a copy of $u$ on some separate probability space $(\tilde{\omega}, \mathcal{F}, 
\tilde{\mathbb{P}})$ 
\end{theorem}

\begin{proof}
The existence and uniqueness of a mild solution is rather standard and follows from similar arguments as in 
the proof of \cite[Proposition~2.8]{CGKPR22}. For the mild Itô formula we define for fixed $t\in [t_0,T]$ 
and $t_0\le s\le t$
\begin{align*}
\overline{u}^{t}_s := e^{(t-t_0)A} u_{t_0} + \int_{t_0}^s e^{(t-r)A}b(u_r,\mathcal{L}(u_r))dr  
+ \int_{t_0}^se^{(t-r)A} \sigma(u_r, \mathcal{L}(u_r))dW_r.
\end{align*}
Now we first observe that
\begin{equation} 
\label{o1}
\overline{u}_t^t=u_t
\end{equation}
and
\begin{equation}
\label{o2}
\overline{u}^{t}_s=e^{(t-s)A}u_s
\end{equation}
almost surely. Let
\begin{align*}
g(s,x):=f(s,x,\mathcal{L}(\overline{u}_s^{t})),
\end{align*}
then by Theorem \ref{Ito1} we have for fixed $x\in H$
\begin{align*}
& g(s,x)=f(s,x,\mathcal{L}(\overline{u}_s^{t})) \\
& = f(t_0,x,\mathcal{L}(\overline{u}_{t_0}^{t}))+\int_{t_0}^s\partial_rf(r,x,\mathcal{L}(\overline{u}_r^t))dr 
	+ \int_{t_0}^s \mathbb{E}\left[\partial_{\mu}f(r,x,\mathcal{L}(\overline{u}_r^{t})) 
	(\overline{u}_r^{t})e^{A(t-r)}b(u_r,\mathcal{L}(u_r))\right] dr \\
& \quad +\frac{1}{2}\int_{t_0}^s\mathbb{E}\left[\textbf{tr}_H(\partial_y\partial_{\mu}f(r,x,\mathcal{L}
	(\overline{u}_r^{t}))(\overline{u}_r^{t})e^{A(t-r)}\sigma(u_r, \mathcal{L}(u_r)) Qe^{A(t-r)} 
	\sigma(u_r, \mathcal{L}(u_r))^{\ast})\right]dr.
\end{align*}
Thus due to the continuity assumptions on $b$ and $\sigma$, $g$ is differentiable in $s$ with 
\begin{equation} 
\label{der}
\begin{aligned}
\partial_sg(s,x) 
& = \partial_sf(s,x,\mathcal{L}(\overline{u}_r^t))+\mathbb{E}\left[\partial_{\mu}f(s,x, 
	\mathcal{L}(\overline{u}_s^{t}))(\overline{u}_s^{t})e^{(t-s)A}b(u_s,\mathcal{L}(u_s))\right] \\
& \quad + \frac{1}{2} \mathbb{E}\left[\textbf{tr}_H(\partial_y\partial_{\mu}f(s,x,\mathcal{L} 
	(\overline{u}_s^{t}))(\overline{u}_s^{t})e^{(t-s)A}\sigma(u_s, \mathcal{L}(u_s)) Qe^{(t-s)A} 
	\sigma(u_s, \mathcal{L}(u_s))^{\ast})\right].
\end{aligned}
\end{equation} 
By \eqref{o1} and the mild Itô formula \cite[Theorem~1]{DPJR19} we have 
\begin{align*}
& f(s,u_t,\mathcal{L}(u_t))=g(t,u_t) \\
& = g(t_0,e^{A(t-t_0)}u_{t_0})+ \int_{t_0}^t\partial_s g(s,e^{(t-s)A}u_s)ds \\
& \quad + \int_{t_0}^t\partial_x g(s,e^{(t-s)A}u_s)e^{(t-s)A}b(u_s,\mathcal{L}(u_s))ds  
	+ \int_{t_0}^t\partial_x g(s,e^{(t-s)A}u_s)e^{(t-s)A}\sigma(u_s,\mathcal{L}(u_s))dW_s \\
& \quad + \frac{1}{2}\int_{t_0}^t \textbf{tr}_H(\partial_{xx}g(s,e^{(t-s)A}u_s)e^{(t-s)A} 
	\sigma(u_s,\mathcal{L}(u_s))Qe^{(t-s)A}\sigma(u_s,\mathcal{L}(u_s))^{\ast})ds.
\end{align*}
With \eqref{der} and \eqref{o2} we arrive at 
\begin{align*}
& f(t,u_t,\mathcal{L}(u_t))=f(t_0,e^{A(t-t_0)}u_{t_0},\mathcal{L}(e^{A(t-t_0)}u_{t_0})) 
+ \int_{t_0}^t\partial_sf(s,e^{(t-s)A}u_s,\mathcal{L}(e^{(t-s)A}u_s))ds \\
& \quad + \int_{t_0}^t \tilde{\mathbb{E}}\left[\partial_{\mu}f(s,e^{(t-s)A}u_s,\mathcal{L}(e^{(t-s)A}u_s)) 
(e^{(t-s)A}\tilde{u}_s)e^{(t-s)A}b(\tilde{u}_s,\mathcal{L}(u_s))\right]ds \\
& \quad + \frac{1}{2}\int_{t_0}^t \tilde{\mathbb{E}}\left[\textbf{tr}_H(\partial_y\partial_{\mu}f(s,u_s,
	\mathcal{L}(e^{(t-s)A}u_s))(e^{(t-s)A}\tilde{u}_s)e^{(t-s)A}\sigma(\tilde{u}_s, \mathcal{L}(u_s))  
	Qe^{(t-s)A}\sigma(\tilde{u}_s, \mathcal{L}(u_s))^{\ast})\right]ds \\
& \quad + \int_{t_0}^t\partial_x f(s,e^{(t-s)A}u_s,\mathcal{L}(e^{(t-s)A}u_s))e^{(t-s)A}b(u_s, 
	\mathcal{L}(u_s))ds \\
& \quad + \int_{t_0}^t\partial_x f(s,e^{(t-s)A}u_s,\mathcal{L}(e^{(t-s)A}u_s))e^{(t-s)A}\sigma(u_s, 
	\mathcal{L}(u_s))dW_s \\
& \quad + \frac{1}{2} \int_{t_0}^t \textbf{tr}_H(\partial_{xx}f(s,e^{(t-s)A}u_s, \mathcal{L} 
	(e^{(t-s)A}u_s))e^{(t-s)A}\sigma(u_s,\mathcal{L}(u_s))Qe^{(t-s)A}\sigma(u_s,\mathcal{L}(u_s))^{\ast})ds.
\end{align*}
\end{proof}
	
\subsection{Stochastic Taylor Expansion for Mean-Field SPDEs}
	
Based on the Itô formula in Theorem \ref{Ito2} we can derive a (second order) stochastic Taylor expansion for 
Mean-Field SPDEs. For simplicity we consider the Mean-Field SPDE \eqref{MFSPDE} with $\sigma\equiv B$, for some 
$B\in L_2(H,H)$:
\begin{equation} 
\label{MFSPDE2}
\begin{aligned}
du_t & = \left[Au_t+b(\mathcal{L}(u_t))\right]dt + BdW_t,\quad t\in [t_0,T] \\
u_{t_0} & = \phi.
\end{aligned}
\end{equation}
Let $(e_i)_{i\in \mathbb{N}}$ be an ONB of $H$ such that $A$ is given by 
\begin{align*}
Au=\sum_{i=1}^{\infty}-\lambda_i\langle u,e_i\rangle_He_i,
\end{align*} 
with $\inf\{\lambda_i|i\in \mathbb{N}\}>-\infty$ and domain
\begin{align*}
D(A):=\{u\in H:\sum_{i=1}^{\infty}|\lambda_i|^2|\langle u,e_i\rangle_H|^2<\infty\}.
\end{align*} 
We fix $\kappa\ge 0$ with $\sup\{\kappa+\lambda_i|i\in \mathbb{N}\}>0$. Similar to \cite{JK10} we consider 
the following stronger setting.
	
\begin{enumerate}[start=0,label={(\bfseries R\arabic*):}]
\item There exists $\gamma\in ]0,1[$ and $\delta\in ]0,\frac{1}{2}]$ such that 
\begin{align*}
& \int_0^T\|(\kappa-A)^{\gamma}e^{As}B\|_{L_2(H,H)}^2<\infty,\\
& \int_0^t\|e^{As}B\|_{L_2(H,H)}^2ds\le C t^{2\delta}.
\end{align*}
\item The drift satisfies $b\in \mathcal{C}^2_b(\mathcal{P}_2(H),H)$ with
\begin{align*}
& \|\partial_{\mu}b(\mu_1)(X_1)-\partial_{\mu}b(\mu_2)(X_2)\|_{L^2(\Omega ; L(H,H))} 
	\le C\|X_1-X_2\|_{L^2(\Omega ;H)}, \\
& \|\partial_y\partial_{\mu}b(\mu_1)(X_1)-\partial_y\partial_{\mu}b(\mu_2)(X_2)\|_{L^2(\Omega ;L(L_2(H,H),H))} 
	\le C\|X_1-X_2\|_{L^2(\Omega ;H)},
\end{align*}
for all $\mu_1,\mu_2\in \mathcal{P}_2(H)$ and $X_1,X_2\in L^2(\Omega ;H)$ with $X_1\sim \mu_1,X_2\sim \mu_2$.
\item We have
\begin{align*}
\mathbb{E}\left[\|(\kappa-A)^{\gamma}\phi\|_H^p\right]<\infty,
\end{align*}
for all $p\ge 1$.
\end{enumerate}
	
\begin{theorem}
Let $u$ be the mild solution to \eqref{MFSPDE2}. Then it holds
\begin{align*}
u_t & = e^{A(t-t_0)}u_{t_0}+\int_{t_0}^te^{(t-s)A}b(\mathcal{L}(e^{A(s-t_0)}u_{t_0}))ds 
	+ \int_{t_0}^te^{(t-s)A}BdW_s \\
& \quad + \int_{t_0}^t\int_{t_0}^se^{(t-s)A}\tilde{\mathbb{E}}\left[\partial_{\mu}b(\mathcal{L} (e^{A(s-r)}u_r))  
	(e^{A(s-r)}\tilde{u}_r)e^{A(s-r)}b(\mathcal{L}(u_r))\right]drds \\
& \quad + \frac{1}{2}\int_{t_0}^t \int_{t_0}^se^{(t-s)A}\tilde{\mathbb{E}}\left[\textbf{tr}_H 
	(\partial_y\partial_{\mu}b(\mathcal{L}(e^{A(s-r)}u_r))(e^{A(s-r)}\tilde{u}_r)e^{A(s-r)}BQe^{A(s-r)}B^{\ast})
	\right]drds.
\end{align*}
Under the additional assumptions (R0)-(R2) we have 
\begin{equation} 
\label{Taylor1}
\begin{aligned}
u_t & = e^{A(t-t_0)}u_{t_0}+\int_{t_0}^te^{(t-s)A}b(\mathcal{L}(e^{A(s-t_0)}u_{t_0}))ds 
	+ \int_{t_0}^te^{(t-s)A}BdW_s \\
& \quad + \int_{t_0}^t\int_{t_0}^se^{(t-s)A}\tilde{\mathbb{E}}\left[\partial_{\mu}b(\mathcal{L}(e^{(s-r)A}
	u_{t_0}))  (e^{(s-r)A}\tilde{u}_{t_0})e^{(s-r)A}b(\mathcal{L}(u_{t_0}))\right]drds \\
& \quad + \frac{1}{2}\int_{t_0}^t \int_{t_0}^se^{(t-s)A}\tilde{\mathbb{E}} \left[\textbf{tr}_H 
	(\partial_y\partial_{\mu}b(\mathcal{L}(e^{(s-r)A}u_{t_0}))(e^{(s-r)A}\tilde{u}_{t_0}) 
	e^{(s-r)A}BQe^{(s-r)A}B^{\ast})\right]drds \\
& \quad + O((t-t_0)^{2+\min(\delta,\gamma)}).
\end{aligned} 
\end{equation}
\end{theorem}
	
\begin{remark}
It should be noted that in \eqref{Taylor1} we obtain an approximation of order $2+\min(\delta,\gamma)$ although 
there are no terms involving second order $L$-derivatives $\partial_{\mu}^2$. This is indeed not surprising since 
there are no such terms present in the Ito formula from Theorem \ref{Ito2}. This type of phenomena was already 
observed in the finite dimensional situation, in particular in the Pontryagin maximum principle for mean-field 
control problems in \cite{BDL11} and \cite{BLM16}.
\end{remark}

\begin{proof}
Starting with \eqref{mild} we have
\begin{align*}
u_t=e^{(t-t_0)A}u_{t_0}+\int_{t_0}^te^{(t-s)A}b(\mathcal{L}(u_s))ds+\int_{t_0}^te^{(t-s)A}BdW_s.
\end{align*}
Now by Theorem \ref{Ito2} we have 
\begin{align*}
u_t & = e^{(t-t_0)A}u_{t_0}+\int_{t_0}^te^{(t-s)A}b(\mathcal{L}(e^{(s-t_0)A}u_{t_0}))ds 
	+ \int_{t_0}^te^{(t-s)A}BdW_s \\
& \quad + \int_{t_0}^t\int_{t_0}^se^{(t-s)A}\tilde{\mathbb{E}}\left[\partial_{\mu}b(\mathcal{L}(e^{(s-r)A}u_r))
	(e^{(s-r)A}\tilde{u}_r)e^{(s-r)A}b(\mathcal{L}(u_r))\right]drds \\
& \quad + \frac{1}{2}\int_{t_0}^t \int_{t_0}^se^{(t-s)A}\tilde{\mathbb{E}} \left[\textbf{tr}_H 
	(\partial_y\partial_{\mu}b(\mathcal{L}(e^{(s-r)A}u_r))(e^{(s-r)A}\tilde{u}_r) 
	e^{(s-r)A}BQe^{(s-r)A}B^{\ast})\right]drds.
\end{align*}
Thanks to the smoothness assumptions the second part of the theorem follows immediately from
\begin{align*}
\|e^{(t-t_0)A}u_{t_0}\|_{L^2(\Omega ;H)} 
& \le C(t-t_0)^{\gamma}, \\
\|\int_{t_0}^te^{(t-s)A}\partial_{\mu}b(\mathcal{L}(u_s))(\tilde{u}_s)ds\|_{L^2(\Omega ;H)} 
& \le C(t-t_0),\\
\|\int_{t_0}^te^{(t-s)A}BdW_s\|_{L^2(\Omega ;H)} & \le C(t-t_0)^{\delta},
\end{align*}
which can be shown similar to the proof of \cite[Theorem~1]{JK10}.
\end{proof}
	
\section{Outlook on Control of Mean-Field SPDEs}
\label{outlook}
	
As mentioned in the introduction, the $L$-derivative is well suited to control infinitesimal perturbations by 
variations on the space of random variables, which naturally appear in the proof of the stochastic maximum 
principle. There is already some literature regarding the stochastic maximum principle for mean-field SPDEs, 
see e.g. \cite{OSD18}, \cite{Ahm14}, \cite{Ahm15}, \cite{Ahm16}, however the authors had to directly work with 
the lifted version of the controlled equation or had to formulate the maximum principle in the space of signed 
measures. This leads to several difficulties, in particular when it comes to the characterization of the 
adjoint state via the adjoint equation. The new extension of the $L$-derivative for the infinite dimensional 
setting now provides the right tool to directly derive the maximum principle for mean-field SPDEs in the spirit 
of \cite{CD18}. Furthermore there was some recent development regarding the Pontryagin maximum principle for 
SPDEs, see \cite{SW21}, with a new characterization of the second order adjoint state via a SPDE. The natural 
next step is to extend the Pontryagin maximum principle for mean-field SDEs in \cite{BDL11} and \cite{BLM16} to 
mean-field SPDEs, using the extension of the $L$-derivative and the novel characterization of the second order 
adjoint state from \cite{SW21}.
	
\appendix

\section{Proof of Proposition \ref{PropertiesH}} 
\label{AppendixA}

(i) The lift of $H$ is given by 
\begin{align*}
\hat{H} (X) 
& = \int_\mathbb{R}  \frac{e^{-x^2/2}}{\sqrt{2\pi}} \left( \int_x^\infty \frac{e^{-y^2/2}}{\sqrt{2\pi}} 
F_{\mathcal{L} (X)} (y)\, dy \right)^2 \, dx \\ 
& \qquad- 2\int_\mathbb{R}  \frac{e^{-x^2/2}}{\sqrt{2\pi}} \int_x^\infty \frac{e^{-y^2/2}}{\sqrt{2\pi}} 
F_{\mathcal{L} (X)} (y)\, dy \int_x^\infty  \frac{e^{-y^2/2}}{\sqrt{2\pi}} \Phi (y)\, dy\, dx \\   
& \qquad + \int_\mathbb{R} \frac{e^{-x^2/2}}{\sqrt{2\pi}} \left( \int_x^\infty \frac{e^{-y^2/2}} 
{\sqrt{2\pi}} \Phi (y) \, dy \right)^2\, dx \\ 
& = \int_\mathbb{R} \int_\mathbb{R}  \Phi (y_1\wedge y_2) \frac{e^{-(y_1^2 + y_2^2)/2}}{2\pi} 
F_{\mathcal{L} (X)} (y_1) F_{\mathcal{L} (X)} (y_2)\, dy_1 \, dy_2  \\
& \qquad - \int_\mathbb{R} \frac{e^{-y^2/2}}{\sqrt{2\pi}} (\Phi (y) - \frac 13 \Phi (y)^3 )
F_{\mathcal{L} (X)} (y)\, dy \\ 
& \qquad +  \frac 14 \int_\mathbb{R}  \frac{e^{-x^2/2}}{\sqrt{2\pi}} \left( 1- \Phi (x)^2\right)^2 \, dx \\ 
& =: I + II + III, \text{ say}. 
\end{align*}
For the first term we have 
\begin{align*} 
I & = 2\int_\mathbb{R} \int_\mathbb{R}  1_{\{y_1 \le y_2\}} \Phi (y_1) 
\frac{e^{-(y_1^2 + y_2^2)/2}}{2\pi} F_{\mathcal{L} (X)} (y_1) F_{\mathcal{L} (X)} (y_2)\, dy_1 \, dy_2 , 
\end{align*} 
and if $\tilde{X}$ is an independent copy of $X$ we can write  
$$ 
F_{\mathcal{L} (X)} (y_1) F_{\mathcal{L} (X)} (y_2) = \mathbb{E} \left[ 1_{\{X\le y_1\}}  
1_{\{\tilde{X} \le y_2\}}\right], 
$$
and Fubini's theorem then yields 
\begin{equation} 
\label{LH:eqn1} 
\begin{aligned} 
I & = 2 \mathbb{E}\left[ \int_X^\infty \int_{\tilde{X}}^\infty 1_{\{y_1 \le y_2\}} \Phi (y_1) 
\frac{e^{- (y_1^2 + y_2^2)/2}}{2\pi}\, dy_1\, dy_2 \right] \\
& = 2 \mathbb{E}\left[ \int_{\tilde{X}}^\infty \int_{X}^{y_2} \Phi (y_1)   
\frac{e^{-y_1^2/2}}{\sqrt{2\pi}} \, dy_1 \frac{e^{- y_2^2/2}}{\sqrt{2\pi}}\, dy_2 \right] \\
& = \mathbb{E}\left[ \int_{\tilde{X}}^\infty (\Phi^2 (y_2 ) - \Phi^2 (X))1_{\{y_2\ge X\}} 
\frac{e^{- y_2^2/2}}{\sqrt{2\pi}}\, dy_2 \right] 
 =  \mathbb{E}\left[ \int_{X\vee \tilde{X}}^\infty (\Phi^2 (y_2 ) - \Phi^2 (X) ) 
\frac{e^{- y_2^2/2}}{\sqrt{2\pi}}\, dy_2 \right] \\ 
& =  \mathbb{E}\left[ \frac{1}{3} (1-\Phi^3 (X\vee\tilde{X})  
- \Phi^2 (X) (1-\Phi (X\vee\tilde{X}))\right] . 
\end{aligned} 
\end{equation} 
For the second term we have 
\begin{equation} 
\label{LH:eqn2} 
\begin{aligned} 
II & = - \mathbb{E}\left[ \int_X^\infty \frac{e^{- y^2/2}}{\sqrt{2\pi}} 
(\Phi (y) - \frac 13 \Phi^3 (y)) \, dy \right] 
 = - \mathbb{E}\left[ \frac 12 (1 - \Phi^2 (X)) - \frac{1}{12} (1- \Phi^4 (X))\right],  
\end{aligned} 
\end{equation} 
and finally, using $\int_\mathbb{R} \frac{e^{-x^2/2}}{\sqrt{2\pi}} \Phi (x)^n\, dx = \frac 1{n+1}$, 
\begin{equation} 
\label{LH:eqn3} 
III = \frac 14 \int_\mathbb{R}  \frac{e^{-x^2/2}}{\sqrt{2\pi}} ( 1 -2 \Phi (x)^2 + \Phi (x)^4 )\, dx 
= \frac 14 (1 - \frac{2}{3} + \frac 15) = \frac{2}{15}.
\end{equation} 
Summing up the integrands in \eqref{LH:eqn1}, \eqref{LH:eqn2} and \eqref{LH:eqn3} yields the 
representation (i), since 
\begin{align*}
H_0 (X,\tilde{X}) 
& = \frac 13 (1-\Phi^3 (X\vee\tilde X )) - \Phi^2 (X) (1- \Phi (X\vee\tilde{X}))  \\ 
& \qquad - \frac 12 (1 - \Phi^2 (X)) + \frac{1}{12} (1- \Phi^4 (X)) + \frac{2}{15} \\ 
& =  - \frac{1}{12} \Phi^4 (X) - \frac{1}{2} \Phi^2 (X) - \frac{1}{3} \Phi^3 (X\vee \tilde{X}) 
+ \Phi^2 (X) \Phi (X\vee\tilde{X}) + \frac{1}{20}.
\end{align*} 

\medskip 
\noindent 
(ii): Let $X$, $Y\in L^2 (\Omega ; \mathbb{R})$ and $\tilde{X}$, $\tilde{Y}$ be independent copies. 
We will first derive the following decomposition 
\begin{equation}
\label{PropositionHDecomposition1}
H_0 (X+Y, \tilde{X} + \tilde{Y} ) - H_0 (X,\tilde{X}) = DH_0 (X,\tilde{X})(Y,\tilde{Y}) + R(X,\tilde{X}, Y, 
\tilde{Y}) 
\end{equation} 
with 
\begin{equation} 
\label{PropositionHDecomposition2}
\begin{aligned} 
DH_0 & (X,\tilde{X})(Y,\tilde{Y}) 
 = - (\frac{1}{3} \Phi (X)^3 + \Phi (X))\frac{e^{-X^2/2}}{\sqrt{2\pi}} Y   
+ 2\Phi (X) \Phi (X\vee\tilde{X}) \frac{e^{-X^2/2}}{\sqrt{2\pi}} Y \\
& \qquad + 1_{\{\tilde{X} > X\}} (\Phi (X)^2 - \Phi (\tilde{X})^2) \frac{e^{-\tilde{X}^2/2}}{\sqrt{2\pi}}
\tilde{Y} 
\end{aligned} 
\end{equation} 
and a remainder $R$ satisfying 
\begin{equation}
\label{PropositionHDecomposition3}
\mathbb{E} \left[ |R(X,\tilde{X}, Y, \tilde{Y}) | \right] 
\le C \mathbb{E} [Y^2]  = C\|Y\|^2_{L^2 (\Omega ;\mathbb{R})}
\end{equation} 
for some constant $C$ independent of $X$ and $Y$. 

To this end let us write $H_0 (X,\tilde{X}) = I (X) + II(X,\tilde{X}) + 1/20$, with 
$$
I(X)= - \frac{1}{12} \Phi (X)^4 - \frac{1}{2} \Phi (X)^2 
= -\int_{-\infty}^X ( \frac{1}{3} \Phi (y)^3 + \Phi (y))\frac{e^{-\frac{y^2}{2}}}{\sqrt{2\pi}} \, dy 
$$ 
which we can further decompose as 
\begin{equation} 
\label{PropertiesH:eqn1}
I (X+Y) - I (X) =  -\underbrace{( \frac{1}{3} \Phi (X)^3 + \Phi (X)) 
\frac{e^{-\frac{X^2}{2}}}{\sqrt{2\pi}}\, Y}_{ =: DH_{I} (X)(Y)} + R_{I} (X,Y) 
\end{equation}  
with remainder 
$$
\begin{aligned} 
R_{I} (X,Y) & = \int_X^{X+Y} ( \frac{1}{3} \Phi (y)^3 + \Phi (y))\frac{e^{-\frac{y^2}{2}}}{\sqrt{2\pi}}  
- ( \frac{1}{3} \Phi (X)^3 + \Phi (X))\frac{e^{-\frac{X^2}{2}}}{\sqrt{2\pi}} \, dy \\
& = \int_X^{X+Y} \int_X^y (\Phi (z)^2 + 1)\frac{e^{-z^2}}{2\pi}- (\frac{1}{3} \Phi (z)^3 + \Phi (z))z\frac{e^{-\frac{z^2}{2}}}{\sqrt{2\pi}} \, dz \, dy 
\end{aligned} 
$$ 
satisfying the estimate
\begin{equation} 
\label{PropertiesH:eqn2}
\begin{aligned}  
|R_{I} (X,Y)| & \le \sup_{z\in\mathbb{R}} |\Phi (z)^2 + 1)\frac{e^{-z^2}}{2\pi}- (\frac{1}{3} \Phi (z)^3 + \Phi (z))z\frac{e^{-\frac{z^2}{2}}}{\sqrt{2\pi}}|\,\cdot\, \frac{1}{2} Y^2 . 
\end{aligned} 
\end{equation} 
Concerning the second term  
$$
II(X, \tilde{X}) = - \frac{1}{3} \Phi (X\vee\tilde{X})^3 + \Phi (X)^2 \Phi (X\vee\tilde{X}) 
$$ 
we first decompose 
$$ 
\begin{aligned} 
II & (X+Y, \tilde{X} + \tilde{Y}) - II (X,\tilde{X}) 
 = - \frac{1}{3} \left( \Phi ((X+Y)\vee (\tilde{X} + \tilde{Y}))^3 - \Phi (X\vee \tilde{X})^3 \right) \\
& \qquad + \Phi (X)^2 \left( \Phi ((X+Y)\vee (\tilde{X} + \tilde{Y}))- \Phi (X\vee \tilde{X})\right) \\
& \qquad + (\Phi (X+Y)^2 - \Phi (X)^2)\Phi ((X+Y)\vee (\tilde{X} + \tilde{Y})) \\ 
& = \underbrace{\int_{X\vee\tilde{X}}^{(X+Y)\vee (\tilde{X} + \tilde{Y})} (\Phi^2 (X) - \Phi (y)^2 ) \frac{e^{-\frac{y^2}{2}}}{\sqrt{2\pi}}\, dy}_{=: II_a} + \underbrace{ (\Phi (X+Y)^2 - \Phi (X)^2)\Phi  
 ((X+Y)\vee (\tilde{X} + \tilde{Y}))}_{=: II_b}. 
\end{aligned} 
$$ 
Clearly, 
\begin{equation} 
\label{PropertiesH:eqn3}
II_b = \underbrace{2 \Phi (X)\Phi (X\vee\tilde{X}) \frac{e^{-\frac{X^2}{2}}}{\sqrt{2\pi}}Y}_{=: DH_{II_b} (X, \tilde{X})(Y, \tilde{Y})} + R_{II_b} 
\end{equation} 
with 
\begin{equation} 
\label{PropertiesH:eqn4}
\begin{aligned} 
|R_{II_b}| & = \big| (\Phi (X+Y)^2 - \Phi (X)^2) \Phi ((X+Y)\vee (\tilde{X} + \tilde{Y}))   
-  2 \Phi (X)\Phi (X\vee\tilde{X}) \frac{e^{-\frac{X^2}{2}}}{\sqrt{2\pi}} Y \big| \\ 
& \le |\Phi (X+Y)^2 - \Phi (X)^2||\Phi ((X+Y)\vee (\tilde{X} + \tilde{Y})) - \Phi (X\vee\tilde{X})|  \\
& \qquad + 2\left|\int_X^{X+Y} \Phi (y) \frac{e^{-\frac{y^2}{2}}}{\sqrt{2\pi}}  
- \Phi (X) \frac{e^{-\frac{X^2}{2}}}{\sqrt{2\pi}} \right| |\Phi (X\vee\tilde{X})|\\
& \le \frac{1}{\pi} |Y| |(X+Y)\vee (\tilde{X} + \tilde{Y})- X\vee \tilde{X}|
+  \sup_{z\in\mathbb{R}} |\frac{e^{-z^2}}{2\pi}- \Phi (z))z\frac{e^{-\frac{z^2}{2}}}{\sqrt{2\pi}}| \, Y^2 \\
& \le C (Y^2 + |Y||\tilde{Y}|), 
\end{aligned} 
\end{equation} 
thereby using $|(X+Y)\vee (\tilde{X} +\tilde{Y}) - X\vee\tilde{X} |\le |Y| + |\tilde{Y}|$. 
Concerning $II_a$ we can further decompose 
 \begin{equation} 
\label{PropertiesH:eqn5}
II_a = \underbrace{1_{\{\tilde{X} > X\}}(\Phi (X)^2 - \Phi (\tilde{X})^2) \frac{e^{-\frac{\tilde{X}^2}{2}}}
     {\sqrt{2\pi}} \tilde{Y}}_{=: DH_{II_a} (X, \tilde{X})(Y, \tilde{Y})} + R_{II_a} 
\end{equation} 
with 
\begin{equation} 
\label{PropertiesH:eqn6}
\begin{aligned} 
R_{II_a} & = 1_{\{X > \tilde{X}\}} \int_{X}^{(X+Y)\vee (\tilde{X} + \tilde{Y})}  
(\Phi (X)^2 - \Phi (y)^2) \frac{e^{-\frac{y^2}{2}}}{\sqrt{2\pi}} \, dy  \\ 
& \qquad + 1_{\{\tilde{X} > X\}} \int_{\tilde{X}}^{(X+Y)\vee (\tilde{X} + \tilde{Y})}  
\Phi (X)^2 (\frac{e^{-\frac{y^2}{2}}}{\sqrt{2\pi}} - \frac{e^{-\frac{\tilde{X}^2}{2}}}{\sqrt{2\pi}}) 
+ (\Phi (y)^2 \frac{e^{-\frac{y^2}{2}}}{\sqrt{2\pi}} - \Phi (\tilde{X})^2\frac{e^{-\frac{\tilde{X}^2}{2}}}
{\sqrt{2\pi}} )\, dy \\
& \qquad + 1_{\{\tilde{X} > X\}} ((X+Y)\vee (\tilde{X} + \tilde{Y})- \tilde{X} - \tilde{Y})
(\Phi (X)^2 - \Phi (\tilde{X})^2) \frac{e^{-\frac{\tilde{X}^2}{2}}}{\sqrt{2\pi}},  
\end{aligned} 
\end{equation} 
which can all be estimated in absolute value from above up to some generic constant by $Y^2 + |Y||\tilde{Y}| 
+ \tilde{Y}^2 \le 2Y^2 + 2 \tilde{Y}^2$. This is easy to see for the first two terms, for the third term it 
follows from 
$$
\begin{aligned} 
|1_{\{\tilde{X} > X\}} & ((X+Y)\vee (\tilde{X} + \tilde{Y})- \tilde{X} - \tilde{Y})
(\Phi (X)^2 - \Phi (\tilde{X})^2) \frac{e^{-\frac{\tilde{X}^2}{2}}}{\sqrt{2\pi}}|\\
& = 1_{\{\tilde{X} > X, X+Y > \tilde{X} + \tilde{Y}\}} (X-\tilde{X} + Y - \tilde{Y})
(\Phi (\tilde{X})^2 - \Phi (X)^2) \frac{e^{-\frac{\tilde{X}^2}{2}}}{\sqrt{2\pi}} \\
& \le 1_{\{\tilde{X} > X, X+Y > \tilde{X} + \tilde{Y}\}} \frac{1}{\pi}(Y - \tilde{Y})(\tilde{X}-X) \\
& \le 1_{\{\tilde{X} > X, X+Y > \tilde{X} + \tilde{Y}\}} \frac{1}{\pi}(Y - \tilde{Y})^2 \, .  
\end{aligned} 
$$
Summing up \eqref{PropertiesH:eqn1}, \eqref{PropertiesH:eqn3} and \eqref{PropertiesH:eqn5}, we arrive at the 
decomposition \eqref{PropositionHDecomposition1} with 
$$ 
\begin{aligned} 
DH_0 (X,\tilde{X})(Y,\tilde{Y}) 
& = - (\frac{1}{3} \Phi (X)^3 + \Phi (X))\frac{e^{-X^2/2}}{\sqrt{2\pi}} Y   
+ 2\Phi (X) \Phi (X\vee\tilde{X}) \frac{e^{-X^2/2}}{\sqrt{2\pi}} Y \\
& \qquad + 1_{\{\tilde{X} > X\}} (\Phi (X)^2 - \Phi (\tilde{X})^2) \frac{e^{-\tilde{X}^2/2}}{\sqrt{2\pi}}
\tilde{Y} 
\end{aligned} 
$$
and remainder   
$$ 
R(X,\tilde{X}, Y, \tilde{Y}) = R_I + R_{II_a} + R_{II_b}   
$$ 
satisfying the estimate \eqref{PropositionHDecomposition3}, using \eqref{PropertiesH:eqn2}, 
\eqref{PropertiesH:eqn4} and the discussion concerning $R_{II_a}$ given in \eqref{PropertiesH:eqn6}, 
as well as the fact that $Y$ and $\tilde{Y}$ have the same distribution.  

The decomposition \eqref{PropositionHDecomposition1} now implies that the lift $\hat{H} (X)$ is Fr\'echet 
differentiable on $L^2 (\Omega ;\mathbb{R})$ with derivative 
$$ 
\begin{aligned}
\mathbb{E} [DH_0 (X,\tilde{X}) (Y,\tilde{Y})] 
& = \mathbb{E} [  -(\frac{1}{3} \Phi (X)^3 + \Phi (X))\frac{e^{-X^2/2}}{\sqrt{2\pi}} Y   
+ 2\Phi (X) \Phi (X\vee\tilde{X}) \frac{e^{-X^2/2}}{\sqrt{2\pi}} Y] \\
& \qquad + \mathbb{E}[1_{\{\tilde{X} > X\}} \Phi (X)^2 - \Phi (\tilde{X})^2) \frac{e^{-\tilde{X}^2/2}}
{\sqrt{2\pi}}\tilde{Y} ] \\ 
& = \mathbb{E} [ (- \frac{1}{3} \Phi (X)^3 - \Phi (X)) + 2\Phi (X)^2 F_\mu (X) 
+ 2\Phi (X)\int_{X+}^\infty \Phi (y) \mu (dy) \\
& \qquad + \int_{-\infty}^X \Phi  (y)^2 \mu (dy) - \Phi^2 (X)F_\mu (X)) 
\frac{e^{-X^2/2}}{\sqrt{2\pi}} Y  ] 
\end{aligned}  
$$
thereby using 
$$
\begin{aligned} 
E[2 \Phi (X) & \Phi (X\vee\tilde{X}) \frac{e^{-X^2/2}}{\sqrt{2\pi}} Y] 
 = E[1_{\{\tilde{X}\le X\}} 2\Phi (X)^2 \frac{e^{-X^2/2}}{\sqrt{2\pi}} Y] 
+ E[1_{\{X < \tilde{X}\}} 2\Phi (X)\Phi (\tilde{X}) \frac{e^{-X^2/2}}{\sqrt{2\pi}} Y]  \\ 
& = E[2\Phi (X)^2 F_\mu (X) \frac{e^{-X^2/2}}{\sqrt{2\pi}} Y] 
+ E[2\int_{X+}^\infty \Phi (y)\, \mu (dy) \Phi (X) \frac{e^{-X^2/2}}{\sqrt{2\pi}} Y], 
\end{aligned} 
$$
because of independence of $\tilde{X}$ and $(X, Y)$, and 
$$
\begin{aligned} 
E[1_{\{\tilde{X} > X\}} (\Phi (X)^2 - \Phi (\tilde{X})^2) \frac{e^{-\tilde{X}^2/2}}{\sqrt{2\pi}}\tilde{Y}] 
& = E[1_{\{\tilde{X} \ge X\}} (\Phi (X)^2 - \Phi (\tilde{X})^2) \frac{e^{-\tilde{X}^2/2}}{\sqrt{2\pi}}
\tilde{Y}] \\
& = E[\int_{-\infty}^{\tilde{X}} \Phi (y)^2 \frac{e^{-\tilde{X}^2/2}}{\sqrt{2\pi}}\tilde{Y}] 
- E[ \Phi (\tilde{X})^2 F_\mu (\tilde{X}) \frac{e^{-\tilde{X}^2/2}}{\sqrt{2\pi}}\tilde{Y}] . 
\end{aligned} 
$$
In particular, $H$ is $L$-differentiable with
\begin{align*}
\partial_{\mu} H (\mu)(x) 
& =  \Big( -\frac{1}{3} \Phi (x)^3 - \Phi (x) + \Phi^2 (x) F_\mu (x) \\  
& \qquad + 2 \Phi (x) \int_{x+}^\infty \Phi (y) \, \mu (dy) 
+ \int_{-\infty}^x \Phi^2 (y) \, \mu (dy) \Big) \frac{e^{-x^2/2}}{\sqrt{2\pi}} .
\end{align*}

\medskip 
\noindent 
(iii) If $X\sim \mathcal{N}(0,1)$, hence $F_\mu = \Phi$, then \eqref{DerivativeH} reduces to 
$$ 
\begin{aligned} 
\partial_{\mu} H (\mu)(x) 
& =  \Big( -\frac{1}{3} \Phi (x)^3 - \Phi (x) + \Phi^3 (x)  \\  
& \qquad + 2 \Phi (x) \int_{x+}^\infty \Phi (y) \,  \Phi' (y) \, dy 
+ \int_{-\infty}^x \Phi^2 (y) \, \Phi' (y)\, dy \Big) \Phi' (x) \\ 
& = \Big( -\frac{1}{3} \Phi (x)^3 - \Phi (x) + \Phi^3 (x) + \Phi (x) (1-\Phi (x)^2) 
+ \frac{1}{3} \Phi^3 (x)\Big) \Phi' (x) = 0 
\end{aligned} 
$$
as to be expected, since $H(\mu )\ge H(\mathcal{N}(0,1))$ for all $\mu$. This implies 
$DH_0 (X,\tilde{X}) = 0$. It follows that 
$$
H(\mathcal{L} (X+Y)) - H(\mathcal{L} (X)) 
= \mathbb{E} [DH_0 (X,\tilde{X})(Y,\tilde{Y}) + R(X,\tilde{X}, Y, \tilde{Y})] 
\le L_H \|Y\|^2_{L^2 (\Omega ;\mathbb{R})} , 
$$
which proves \eqref{LipschitzAtNormalDistrH}.

\section{Proof of Theorem \ref{mainEx}} 
\label{AppendixB}

The proof of Theorem \ref{mainEx} requires a couple of further preparations:  

\begin{lemma} 	
\label{SomeIntegral} 
Let $x\neq\tilde{x}\in \mathbb{R}$. Then 
$$ 
\int_\mathbb{R} \frac {\text{ sign }(t-\tilde  {x})}{\sqrt{|t-\tilde{x}|}} 
\frac{\text{ sign }(t-x)}{\sqrt{|t-x|}} 1_{\{|t-\tilde{x}|\le 1, |t-x|\le 1\}}\, dt  
= \Psi_0 (x- \tilde{x})\, ,
$$
with
$$ 
\Psi_0 (s)  = 	
\begin{cases} 
 - \pi  +  4 \ln\left( \frac{1 + \sqrt{1-|s|}}{\sqrt{|s|}}\right)
& \text{ if } 0 < |s|\le 1  \\ 
2\arcsin\left( 1- \frac{2}{|s|}\right) & \text{ if } 1 < |s| \le 2 \\ 
0 & \text{ if } |s| > 2.  
\end{cases} 
$$
\end{lemma} 
	
\begin{proof}
If $|x- \tilde{x}| > 2$ it follows that $1_{\{|t-x|\le 1, |t-\tilde{x}|\le 1\}} \equiv 0$ and thus 
$$
\int_\mathbb{R} \frac {\text{ sign }(t-\tilde{x})}{\sqrt{|t-\tilde{x}|}} 
\frac{\text{ sign }(t-x)}{\sqrt{|t-x|}} 1_{\{|t-\tilde{x}|\le 1, |t-x|\le 1\}}\, dt = 0. 
$$
Next assume that $|x- \tilde{x}| \le 2$ and w.l.o.g. that $\tilde{x} < x$ and distinguish the two cases 
$1 < x-\tilde{x}\le 2$ and $x-\tilde{x}\le 1$.

\medskip 
\noindent 
\textbf{Case 1:} $1 < x-\tilde{x}$, hence $\tilde{x} < x-1 \le \tilde{x} + 1 < x$. Then 
$$
\begin{aligned} 
\int_\mathbb{R} & \frac {\text{ sign }(t-\tilde{x})}{\sqrt{|t-\tilde{x}|}} 
\frac{\text{ sign }(t-x)}{\sqrt{|t-x|}} 1_{\{|t-\tilde{x}|\le 1, |t-x|\le 1\}}\, dt
 = - \int_{x-1}^{\tilde{x}+1} \frac{1}{\sqrt{t-(\tilde{x})}\sqrt{x-t}}\, dt \\ 
& = \arcsin\left( \frac{-2t+ \tilde{x} + x}{x-\tilde{x}}\right)_{\mid_{t = x-1}^{t = \tilde{x} + 1}} 
 = - 2\arcsin\left( \frac{2 - (x- \tilde{x})}{x-\tilde{x}}\right) = \Psi_0 (|x-\tilde{x}|)\, ,   
\end{aligned} 
$$
using $\arcsin (-x) = -\arcsin (x)$. Also note that $\lim_{x-\tilde{x}\to 1} \Psi_0 (x-\tilde{x}) 
= -2 \arcsin (1) = - \pi$. 

\medskip 
\noindent 
\textbf{Case 2:} $0 < x-\tilde{x}\le 1$, hence $\tilde{x} - 1 <  x-1 \le \tilde{x} < x\le \tilde{x} + 1  
< x + 1$. Then 
$$
\begin{aligned} 
\int_\mathbb{R} & \frac{\text{ sign }(t-\tilde{x})}{\sqrt{|t-\tilde{x}|}} 
\frac{\text{ sign }(t-x)}{\sqrt{|t-x|}} 1_{\{|t-\tilde{x}|\le 1, |t-x|\le 1\}}\, dt  \\
& = \int_{x-1}^{\tilde{x}} + \int_{\tilde{x}}^x + \int_{x}^{\tilde{x}+1} \frac{\text{ sign }(t-\tilde{x})} 
{\sqrt{|t-\tilde{x}|}} \frac{\text{ sign }(t-x)}{\sqrt{|t-x|}}\, dt \\ 
& = \int_{x-1}^{\tilde{x}} \frac{1} {\sqrt{\tilde{x}-t}\sqrt{x-t}}\, dt  
- \int_{\tilde{x}}^x  \frac{1} {\sqrt{t-\tilde{x}}\sqrt{x-t}}\, dt 
+ \int_{x}^{\tilde{x}+1} \frac{1} {\sqrt{t-\tilde{x}}\sqrt{t-x}}\, dt \\
& = \ln \left( -2 \left( \sqrt{(\tilde{x}-t)(x-t)} + t\right) + \tilde{x} + x 
\right)_{\mid_{t = x-1}^{t = \tilde{x}}} 
- \arcsin\left(\frac{-2t + \tilde{x} + x}{x-\tilde{x}}\right)_{\mid_{t = \tilde{x}}^{t = x}}  \\ 
& \qquad + \ln \left( 2 \left( \sqrt{(t-\tilde{x})(t-x)} + t\right) - (\tilde{x} + x) 
\right)_{\mid_{t = x}^{t = \tilde{x}+1}} \\ 
& = \pi - \ln \left( -2 \sqrt{1- (x-\tilde{x})} + 2 - (x-\tilde{x})\right) 
 + \ln\left( 2 \sqrt{1- (x-\tilde{x})} + 2 - (x-\tilde{x})\right) \\
& = \pi + 2\ln \left( 1 + \sqrt{1 - (x-\tilde{x})}\right) 
- 2\ln \left( 1 - \sqrt{1 - (x-\tilde{x})}\right) \, ,   
\end{aligned} 
$$
using $\pm 2 \sqrt{1- (x-\tilde{x})} + 2 - (x-\tilde{x}) = (1 \pm \sqrt{1 - (x-\tilde{x})})^2$.
Finally, observe that 
$$ 
\begin{aligned} 
\ln \left( \frac{1 + \sqrt{1 - |s|}}{1 - \sqrt{1-|s|}} \right)  
& = \ln\left( \underbrace{\frac{\sqrt{|s|}}{1 - \sqrt{1-|s|}}}_{=: \kappa_1}\right)
+ \ln\left( \underbrace{\frac{1 + \sqrt{1-|s|}}{\sqrt{|s|}}}_{=: \kappa_2}\right) 
\end{aligned} 
$$ 
with 
$$ 
\frac{\kappa_1}{\kappa_2} = \frac{|s|}{1 - (1 - |s|)} \frac{|s|}{|s|} = 1, 
$$
so that $\kappa_1 = \kappa_2$, and therefore 
$$ 
\begin{aligned} 
2\ln \left( \frac{1+\sqrt{1- |s|}}{1 - \sqrt{1-|s|}} \right)  
& = 2\ln \kappa_1 + 2\ln \kappa_2 = 4 \ln \kappa_2 
= 4 \ln\left( \frac{1 + \sqrt{1-|s|}}{\sqrt{|s|}}\right) \, , 
\end{aligned} 
$$
which proves the assertion.  
\end{proof}

\begin{lemma}
\label{HoelderPsi} 
The function $\Psi_0$, as defined in Lemma \ref{SomeIntegral}, is locally Hölder-continuous on 
$\mathbb{R}\setminus\{0\}$. More specifically, for all $\alpha\in (0, 1)$ 
$$ 
\begin{aligned} 
|\Psi_0 (s ) - \Psi_0(\tilde{s})| 
& \le \sqrt{|s- \tilde{s}|} + \frac{2}{\alpha (s\wedge\tilde{s})^\alpha} |s- \tilde{s}|^\alpha  
\text{ for } s\cdot\tilde{s} > 0. 
\end{aligned} 
$$
\end{lemma}

\begin{proof}
W.l.o.g. we may assume that $0 < \tilde{s} \le s$. Estimating the derivative 
$$
\Psi_0' (s) = 
\begin{cases} 
-\frac{2}{\sqrt{1-s}(1+\sqrt{1-s})} - \frac{2}{s} & \text{ if } 0 < s\le 1 \\
\frac{1}{s\sqrt{1+s}}  & \text{ if } 1 < s \le 2 \\ 
0 & \text{ if } s > 2, 
\end{cases} 
$$ 
from above by 
$$
|\Psi_0' (s)| \le \frac{2}{\sqrt{1-s}}1_{\{0 < s \le 1\}} + \frac{2}{s} , 
$$ 
we obtain that 
$$
\begin{aligned} 
|\Psi_0 (s) - \Psi_0 (\tilde{s})| 
& \le \int_{\tilde{s}}^s \frac{2}{\sqrt{1-u}}1_{\{0 < u \le 1\}} + \frac{2}{u}\, du \\
& \le \sqrt{1-\tilde{s}\wedge 1} - \sqrt{1-s\wedge 1} + \frac{2}{\tilde{s}^\alpha}  
\int_{\tilde{s}}^s u^{\alpha - 1} \, du \\
& \le \sqrt{s\wedge 1 - \tilde{s}\wedge 1} + \frac{2}{\alpha\tilde{s}^\alpha} (s- \tilde{s})^\alpha \\ 
& \le \sqrt{s - \tilde{s}} + \frac{2}{\alpha\tilde{s}^\alpha} (s- \tilde{s})^\alpha. 
\end{aligned} 
$$ 
\end{proof}

\begin{lemma}
\label{Integrability} 
The following estimate holds:
$$
\frac{1}{\sqrt{2\pi}} \int_{\mathbb{R}} \frac{1}{|x-a|^\beta}e^{-\frac{x^2}{2}}\, dx 
\le 1 + \sqrt{\frac{2}{\pi}} \frac{1}{1-\beta} \text{ for all } a\in\mathbb{R}, \beta \in (0,1) .  
$$
\end{lemma}

\begin{proof}
Follows from 
$$
\begin{aligned} 
\frac{1}{\sqrt{2\pi}} \int_{\mathbb{R}} \frac{1}{|x-a|^\beta}e^{-\frac{x^2}{2}}\, dx 
& \le \frac{1}{\sqrt{2\pi}} \int_{|x-a|\ge 1} e^{-\frac{x^2}{2}}\, dx + 
\frac{1}{\sqrt{2\pi}} \int_{|x-a| < 1} \frac{1}{|x-a|^\beta}\, dx \\
& \le 1 + \sqrt{\frac{2}{\pi}} \frac{1}{1-\beta} \, .   
\end{aligned} 
$$
\end{proof}

\medskip 
\noindent 
We are now ready to give the proof of Theorem \ref{mainEx}.  

\begin{proof} 
\textbf{(of Theorem \ref{mainEx})}	
First note that we can represent the lift $\hat{f} (X) = f (\mathcal{L}(X))$ as 
\begin{equation}
\label{RepresentationLift}
\hat{f} (X) (t) = \mathbb{E} \left[ \int_\mathbb{R} g_0 (t-\omega , H(\mathcal{L}(X)))1_{\{X\le \omega \}} 
\, d\omega \right] = \mathbb{E} [G_0 (t-X, H(\mathcal{L} (X)))].  
\end{equation}
We then note that for all $X\in L^2 (\Omega ; \mathbb{R})$ and $\tilde{\mu}\in \mathcal{P}_2(\mathbb{R})$, 
the elementary estimate $|G_0(t -X ,H(\tilde{\mu}))|\le 2 \sqrt{1 + H(\tilde{\mu})} 1_{\{|t|\le 1\}}$ 
implies that 
$$
\begin{aligned} 
\int_\mathbb{R} \mathbb{E} [ |G_0 (t - X, H(\tilde{\mu}))| ]^2 \, dt 
& \le \int_\mathbb{R} 4 (1 + H(\tilde{\mu})) \mathbb{E} [ 1_{\{|t-X|\le 1\}}] \, dt 
= 8 (1+ H(\tilde{\mu})) < \infty  
\end{aligned} 
$$
and thus $\mathbb{E} [ |G_0 (\cdot - X, H(\tilde{\mu}))|] \in L^2 (\mathbb{R})$. 

\medskip 
\noindent 
Then note that $\hat{f} (X)$ can be written as the composition of the map 
$$
\begin{aligned} 
\tilde{G}_0 & : L^2 (\Omega ; \mathbb{R}) \times [0, \infty )\rightarrow L^2 (\mathbb{R}) \\
& \qquad (X,y) \mapsto \mathbb{E} \left[ G_0 (\cdot - X, y)\right] 
\end{aligned} 
$$
with the map 
$$
\begin{aligned} 
\tilde{H} & : L^2 (\Omega ; \mathbb{R}) \rightarrow  L^2 (\Omega ; \mathbb{R}) \times [0, \infty ) \\
& \qquad X \mapsto (X, H(\mathcal{L}(X))). 
\end{aligned} 
$$
It is easy to show that $\tilde{G}_0$ is Fr\'echet differentiable on $L^2(\Omega ; \mathbb{R})\times  
(0, \infty )$. Hence Proposition \ref{PropertiesH} and the chain rule imply Fr\'echet differentiability 
of $\hat{f} (X)$ in $X$ with $H( \mathcal{L} (X)) > 0$ with derivative \eqref{RepresentationDifferential1}. 

\medskip 
\noindent 
It thus remains to show that $\hat{f} (X)$ is Fr\'echet differentiable in $X$ with $H( \mathcal{L} (X)) = 0$, 
or equivalently with $X\sim \mathcal{N} (0,1)$. To this end we will show below for 
$Y\in L^2 (\Omega , \mathbb{R})$
\begin{equation}
\label{RepresentationDifferential3} 
\begin{aligned}  
G_0 (t-(X+Y), & H(\mathcal{L}(X+Y))) - G_0 (t-X,H(\mathcal{L}(X))) \\ 
& = \underbrace{G_0 (t-(X+Y),H(\mathcal{L}(X+Y)))  
    - G_0 (t-(X+Y),0)}_{ =: R_I (t)} \\ 
& \qquad + \underbrace{G_0 (t-(X+Y),0) -G_0 (t-X,0)}_{ =: II(t)}\, , 
\end{aligned} 
\end{equation}
the estimate 
\begin{equation}
\label{RemainderDifferentialI} 
\begin{aligned} 
\int_{\mathbb{R}} \mathbb{E} [ R_I (t)]^2  \, dt
\le C \|Y\|^3_{L^2 (\Omega ; \mathbb{R})} ,  
\end{aligned} 
\end{equation}
and the decomposition 
\begin{equation}
\label{RepresentationDifferentialII} 
\begin{aligned}  
II(t) = -g_0 (t-X, H(\mathcal{L} (X))Y + R_{II}(t) , 
\end{aligned} 
\end{equation}
together with the estimate 
\begin{equation}
\label{RemainderDifferentialII} 
\int_{\mathbb{R}} \mathbb{E} [ R_{II} (t)]^2  \, dt
\le C \left(\|Y\|_{L^2 (\Omega ; \mathbb{R})}^{5/2} +  \|Y\|_{L^2 (\Omega ; \mathbb{R})}^{2 + \alpha} \right)  
\end{equation} 
for $\alpha\in (0, 1/4)$. This implies the Fr\'echet differentiability at $X$ with the Fr\'echet differential 
\eqref{RepresentationDifferential1} in the case $H(\mathcal{L} (X)) = 0$, since 
$$
\begin{aligned}  
\mathbb{E} [G_0 (t-(X+Y), & H(\mathcal{L}(X+Y))) - G_0 (t-X,H(\mathcal{L}(X)))]
 -\mathbb{E} [g_0(t-X, H(\mathcal{L}(X))Y] \\ 
& =  R_I (t) + R_{II} (t), 
\end{aligned} 
$$
together with the estimates \eqref{RemainderDifferentialI} and \eqref{RemainderDifferentialII}, implies
\begin{equation}
\label{Differentiability1} 
\begin{aligned}  
\int_{\mathbb{R}} \mathbb{E} \big[ ( G_0 (t-(X+Y), H(\mathcal{L}(X+Y))) & - G_0 (t-X,H(\mathcal{L}(X)))  
- g_0(t-X, H(\mathcal{L}(X))Y \big]^2 \, dt  
\\ 
& = o (\|Y\|^2_{L^2 (\Omega ; \mathbb{R})}) . 
\end{aligned} 
\end{equation}
In particular, $f$ is $L$-differentiable with 
\begin{equation}
\partial_\mu f(\mu )(x) = - g_0 (t-x, H(\mu ))  
- \int \partial_y G_0 (t-\tilde{x}, H(\mu )) \mu (d\tilde{x}) \partial_\mu H(\mu )(x) . 
\end{equation} 
Since, however, for $\mu = \mathcal{N}(0,1)$ ,  
\begin{equation}
\int_{\mathbb{R}} \partial_\mu f(\mu )(x) = - g_0 (t-x, H(\mu )) , 
\end{equation} 
it is clear, that 
\begin{equation}
\partial_\mu f( \mathcal{N}(0,1))(x) = - g_0 (t-x, 0) = -\frac{sign (t-x)}{\sqrt{|t-x|}} 
1_{\{|t-x|\le 1\}} \notin L^2 (\mathbb{R}),  
\end{equation} 
hence $f$ is not absolute continuous $L$-differentiable in $\mathcal{N} (0,1 )$.

\bigskip 
\noindent 
It remains to prove estimates \eqref{RemainderDifferentialI} and \eqref{RemainderDifferentialII}. 

\smallskip 
\noindent
\textbf{Concerning $R_I$:} Since $H(\mathcal{L}(X+Y))\le C\|Y\|^2_{L^2 (\Omega ;\mathbb{R})}$, we can estimate  
\begin{align*}
|R_I(t)| & = 2 |\sqrt{1+H(\mathcal{L} (X+Y))} - \sqrt{|t-(X+Y)|+H(\mathcal{L} (X+Y))} \\
& \qquad\qquad -(1 - \sqrt{|t-(X+Y)|})| 1_{\{|t-(X+Y)|\le 1\}} \\ 
& \le 2 \left( \sqrt{|t-(X+Y)|+ H(\mathcal{L}(X+Y))}- \sqrt{|t-(X+Y)|} \right) \, 1_{\{|t-(X+Y)|\le 1\}}\\
& = \int_{0}^{H(\mathcal{L}(X+Y))} \frac{1}{\sqrt{|t-(X+Y)| + u}} \, du  \, 1_{\{|t-(X+Y)|\le 1\}} \\
& \le \left( \frac{H(\mathcal{L} (X+Y))}{\sqrt{|t-(X+Y)|}} \right)^\alpha 
\left( \int_{0}^{H(\mathcal{L}(X+Y))} \frac{1}{\sqrt{u}} \, du \right)^{1-\alpha}\, 1_{\{|t-(X+Y)|\le 1\}} \\
& =\left( \frac{H(\mathcal{L} (X+Y))}{\sqrt{|t-(X+Y)|}} \right)^\alpha 
(2H(\mathcal{L} (X+Y)))^{(1-\alpha )/2} \, 1_{\{|t-(X+Y)|\le 1\}} \\ 
& \le C \|Y\|^{1 + \alpha}_{L^2 (\Omega ;\mathbb{R})} \frac{1}{|t-(X+Y)|^{\alpha /2}} \, 1_{\{|t-(X+Y)|\le 1\}}
\end{align*}
for all $\alpha\in (0, 1)$. If we choose in particular $\alpha = 1/2$, then
\begin{equation} 
\label{EstimateI.5}
\begin{aligned}
|R_I (t)| \le C \|Y\|^{3/2}_{L^2 (\Omega ;\mathbb{R})} \frac{1}{|t-(X+Y)|^{1/4}} 1_{\{|t-(X+Y)|\le 1\}} \, , 
\end{aligned}
\end{equation} 
therefore 
$$
\begin{aligned} 
\int_{\mathbb{R}} \mathbb{E} [ R_{I} (t)]^2  \, dt 
\le C \mathbb{E} \left[ \int_{\mathbb{R}} \frac{1}{|t-(X+Y)|^{1/2}} 1_{\{|t-(X+Y)|\le 1\}}\, dt  \right] 
\|Y\|^3_{L^2 (\Omega ; \mathbb{R})} 
= 4C \|Y\|^3_{L^2 (\Omega ; \mathbb{R})} ,  
\end{aligned} 
$$ 
so that \eqref{RemainderDifferentialI} holds.  

\medskip 
\noindent
\textbf{Concerning $II$:} Using $H(\mathcal{L}(X)) = 0$, we have  
\begin{equation} 
\label{EstimateII.1}
\begin{aligned}
R_{II} (t) & = G_0 (t-(X+Y), 0) - G_0 (t-X,0) - \frac{\text{sign }(t-X)}{\sqrt{|t-X|}}  
1_{\{|t-X|\le 1\}}Y \\ 
& = \int_{X}^{X+Y} \frac{\text{sign }(t-s)}{\sqrt{|t-s|}}1_{\{|t-s|\le 1\}}  
- \frac{\text{sign }(t-X)}{\sqrt{|t-X|}}1_{\{|t-X|\le 1\}}\, ds .
\end{aligned} 
\end{equation}  
Let $(\tilde{X},\tilde{Y})$ be an independent copy of $(X,Y)$. Then 

\begin{equation} 
\label{EstimateII.2}
\begin{aligned}
\mathbb{E} [ R_{II} (t)]^2  
& = \mathbb{E} \Big[ \int_{X}^{X+Y} \int_{\tilde{X}}^{\tilde{X}+\tilde{Y}} 
\Big( \frac{\text{sign } (t-s)}{\sqrt{|t-s|}}1_{\{|t-s|\le 1\}} - \frac{\text{sign }(t-X)} 
{\sqrt{|t-X|}}1_{\{|t-X|\le 1\}}\Big) \\ 
& \qquad \qquad 
\Big( \frac{\text{sign }(t-\tilde{s})}{\sqrt{|t-\tilde{s}|}} 
1_{\{|t-\tilde{s}|\le 1\}} - \frac{\text{sign }(t-\tilde{X})}{\sqrt{|t-\tilde{X}|}} 
1_{\{|t-\tilde{X}|\le 1\}} \Big) \, ds \, d\tilde{s} \Big].
\end{aligned} 
\end{equation}  
Integration w.r.t. $t$ and using Lemma \ref{SomeIntegral} yields 
\begin{equation} 
\label{EstimateII.3}
\begin{aligned}
\int_{\mathbb{R}} \mathbb{E} [ R_{II} (t)]^2  \, dt 
& = \mathbb{E} \Big[ \int_{X}^{X+Y} \int_{\tilde{X}}^{\tilde{X}+\tilde{Y}} \Psi_0 ( s-\tilde{s}) - 
\Psi_0 (s-\tilde{X}) - \Psi_0 (X-\tilde{s}) \\ 
& \qquad\qquad 
+ \Psi_0 (X-\tilde{X})  \, ds \, d\tilde{s} \Big].
\end{aligned} 
\end{equation} 

\medskip 
\noindent 
We now rewrite \eqref{EstimateII.3} with the help 
of the new variables $u = s-X$ and $\tilde{u} = \tilde{s} - \tilde{X}$ as 
\begin{equation} 
\label{EstimateII.6}
\begin{aligned}
\int_{\mathbb{R}} \mathbb{E} [ R_{II} (t)]^2  \, dt 
& = \mathbb{E} \Big[ \int_{0}^{Y} \int_{0}^{\tilde{Y}} \Psi_0 ( u-\tilde{u} + X-\tilde{X})  
- \Psi_0 ( u + X-\tilde{X}) - \Psi_0 (-\tilde{u} + X-\tilde{X})\\ 
& \qquad\qquad + \Psi_0 (\tilde{X}-\tilde{X}) \, du \, d\tilde{u} \Big].  
\end{aligned} 
\end{equation}
Next define the set $A:= \{ 2(|Y|+ |\tilde{Y}|) \le |X-\tilde{X}|\}$. 
On this set $\min\{ |u-\tilde{u} + X-\tilde{X}|, |u + X-\tilde{X}|, |-\tilde{u} + X-\tilde{X}|, 
|X-\tilde{X}|\} \ge \frac{1}{2} |X-\tilde{X}|$ for all $|u|\le |Y|$ and all $|\tilde{u}|\le |\tilde{Y}|$. 
In particular, $u-\tilde{u} + X-\tilde{X}$, $-\tilde{u} + X-\tilde{X}$, $\tilde{u} + X-\tilde{X}$ and 
$X-\tilde{X}$ all have the same sign, and Lemma \ref{HoelderPsi} now implies for $\alpha\in (0, 1)$ that  
$$ 
|\Psi_0 (u-\tilde{u} + X-\tilde{X}) - \Psi_0 ( u + X-\tilde{X}) | 
\le  \sqrt{|\tilde{u}|} + \frac{2}{\alpha} \frac{2^\alpha}{|X-\tilde{X}|^\alpha} |\tilde{u}|^\alpha 
$$ 
and 
$$ 
|\Psi_0 (-\tilde{u} + X-\tilde{X}) - \Psi_0 (X-\tilde{X}) | 
\le \sqrt{|\tilde{u}|} + \frac{2}{\alpha} \frac{2^\alpha}{|X-\tilde{X}|^\alpha} |\tilde{u}|^\alpha   
$$ 
so that 
\begin{equation} 
\label{EstimateII.7}
\begin{aligned}
&\Big| \int_{0}^{Y} \int_{0}^{\tilde{Y}} \Psi_0 ( u-\tilde{u} + X-\tilde{X})  
- \Psi_0 ( u + X-\tilde{X}) - \Psi_0 (-\tilde{u} + X-\tilde{X}) \\ 
& \qquad\qquad + \Psi_0 (\tilde{X}-\tilde{X}) \, du \, d\tilde{u} \Big| \\ 
& \le  2 \Big| \int_{0}^{Y} \int_{0}^{\tilde{Y}} \sqrt{|\tilde{u}|} + \frac{2}{\alpha}  
\frac{2^\alpha}{|X-\tilde{X}|^\alpha} |\tilde{u}|^\alpha  \, du \, d\tilde{u} \Big| \\ 
& \le 2|Y| \left( \frac{2}{3} |\tilde{Y}|^{3/2} + \frac{2^{1+\alpha}}{\alpha (1+\alpha )} 
\frac{|\tilde{Y}|^{\alpha + 1}}{|X-\tilde{X}|^\alpha} \right) . 
\end{aligned} 
\end{equation} 
Taking expectation, and using independence of $Y$ and $\tilde{Y}$, we obtain that 
\begin{equation} 
\label{EstimateII.8}
\begin{aligned}
& \mathbb{E} \Big[ 1_A \int_{0}^{Y} \int_{0}^{\tilde{Y}} \Psi_0 ( u-\tilde{u} + X-\tilde{X})  
- \Psi_0 ( u + X-\tilde{X}) - \Psi_0 (-\tilde{u} + X-\tilde{X})\\ 
& \qquad\qquad + \Psi_0 (\tilde{X}-\tilde{X}) \, du \, d\tilde{u} \Big]\\ 
& \le \frac{4}{3} \mathbb{E} [ |Y|]\mathbb{E} [|\tilde{Y}|^{3/2}] 
+ \frac{2^{1+\alpha}}{\alpha (1+\alpha )} \mathbb{E} [ |Y|^2]^{1/2} \mathbb{E} [|\tilde{Y}|^2]^{1/2}  
\mathbb{E} [\frac{|\tilde{Y}|^{2\alpha}}{|X-\tilde{X}|^{2\alpha}}]^{1/2}  \\
& \le  C\left(\|Y\|_{L^2 (\Omega ; \mathbb{R})}^{5/2} +  \|Y\|_{L^2 (\Omega ; \mathbb{R})}^{2 + \alpha} 
\right) ,  
\end{aligned} 
\end{equation}
for $\alpha\in (0, 1/2)$, where the final step follows from integrating out the independent random 
variable $X$ 
\begin{equation} 
\label{EstimateII.9}
\mathbb{E} [\frac{|\tilde{Y}|^{2\alpha}}{|X-\tilde{X}|^{2\alpha}}] 
\le \left( 1 + \sqrt{\frac{2}{\pi}} \frac{1}{1-2\alpha} \right) \mathbb{E} [|\tilde{Y}|^{2\alpha}] 
= C \|Y\|_{L^2 (\Omega ; \mathbb{R})}^{\alpha} , 
\end{equation} 
thereby using Lemma \ref{Integrability}. 

\smallskip 
\noindent 
On the complement $A^c = \{ 2(|Y|+ |\tilde{Y}|) > |X-\tilde{X}|\}$, we use Hölder's inequality  
$$
\begin{aligned} 
\left| \int_0^y \Psi_0 (u+a)\, du \right| 
& \le \left( \int_0^y |\Psi_0 (u+a)|^{1/\alpha}\, du \right)^\alpha  |y|^{1- \alpha} 
\le C |y|^{1-\alpha} , 
\end{aligned} 
$$ 
for $\alpha \in (0,1)$ and some finite constant $C = C(\alpha )$. Then 
\begin{equation} 
\label{EstimateII.10}
\begin{aligned}
& \Big| \int_{0}^{Y} \int_{0}^{\tilde{Y}} \Psi_0 ( u-\tilde{u} + X-\tilde{X})  
- \Psi_0 ( u + X-\tilde{X}) - \Psi_0 (-\tilde{u} + X-\tilde{X})  + \Psi_0 (\tilde{X}-\tilde{X}) \, du \, d\tilde{u} \Big| \\ 
& \qquad 
\le C(|\tilde{Y}| |Y|^{1-\alpha} + |\tilde{Y}|^{1-\alpha} |Y| + |Y||\tilde{Y}| |\Psi_0 (X-\tilde{X})|) \\ 
& \qquad 
\le \frac{C}{|X-\tilde{X}|^{2\alpha}} ( |\tilde{Y}|^{1+\alpha} |Y| + |\tilde{Y}| |Y|^{1+\alpha}  
+ |\tilde{Y}|^{1-\alpha} |Y|^{1+2\alpha} + |\tilde{Y}|^{1+2\alpha} |Y|^{1-\alpha} ) , 
\end{aligned} 
\end{equation} 
thereby using $|\Psi_0 (X-\tilde{X})|\le C|X-\tilde{X}|^\alpha$ and $|X-\tilde{X}|\le 2(|Y| + |\tilde{Y}|)$ 
on the set $A^c$. 
Again taking expectation, using independence of $Y$ and $\tilde{Y}$ and integrating out independent variables 
$X$ (resp. $\tilde{X}$) similar to \eqref{EstimateII.9}, we obtain that 
\begin{equation} 
\label{EstimateII.11}
\begin{aligned}
& \mathbb{E} \Big[ 1_{A^c} \int_{0}^{Y} \int_{0}^{\tilde{Y}} \Psi_0 ( u-\tilde{u} + X-\tilde{X})  
- \Psi_0 ( u + X-\tilde{X}) - \Psi_0 (-\tilde{u} + X-\tilde{X})\\ 
& \qquad\qquad + \Psi_0 (\tilde{X}-\tilde{X}) \, du \, d\tilde{u} \Big]\\ 
& \le C\mathbb{E} [|\tilde{Y}|^2 |Y|^2]^{1/2}  
\left( \mathbb{E} \left[ \frac{|\tilde{Y}|^{2\alpha}}{|X-\tilde{X}|^{4\alpha}} \right]^{1/2} 
+ \mathbb{E} \left[ \frac{|Y|^{2\alpha}}{|X-\tilde{X}|^{4\alpha}} \right]^{1/2} \right) \\
& \qquad + C\left( \mathbb{E} [|\tilde{Y}|^2 |Y|^{2(1-\alpha)}]^{1/2}  + \mathbb{E} [|\tilde{Y}|^{2(1-\alpha )} 
|Y|^2]^{1/2} \right) \left( \mathbb{E} \left[ \frac{|\tilde{Y}|^{4\alpha}}{|X-\tilde{X}|^{4\alpha}} 
\right]^{1/2} 
+ \mathbb{E} \left[ \frac{|Y|^{4\alpha}}{|X-\tilde{X}|^{4\alpha}} \right]^{1/2} \right) \\
& \le  C\|Y\|_{L^2 (\Omega ; \mathbb{R})}^{2+\alpha )}, 
\end{aligned} 
\end{equation}
provided $\alpha\in (0, 1/4)$. 

\smallskip 
\noindent 
Adding up \eqref{EstimateII.8} and \eqref{EstimateII.11} leads to the final estimate 
\begin{equation} 
\label{EstimateII.8d}
\int_{\mathbb{R}} \mathbb{E} [ R_{II} (t)]^2  \, dt
\le C \left(\|Y\|_{L^2 (\Omega ; \mathbb{R})}^{5/2} +  \|Y\|_{L^2 (\Omega ; \mathbb{R})}^{2 + \alpha}\right) 
\end{equation} 
which then implies \eqref{RemainderDifferentialII}. 
\end{proof}


\printbibliography

@article{Ahm13,
	author = {Ahmed, N.U.},
	year = {2013},
	month = {04},
	pages = {},
	title = {A note on Radon-Nikodym theorem for operator valued measures and its applications},
	volume = {28},
	journal = {Commun. Korean Math. S.},
	doi = {10.4134/CKMS.2013.28.2.285}
}

@article{Ahm14,
	author = {Ahmed, N.U.},
	title = {Stochastic neutral evolution equations on Hilbert spaces with
	partially observed relaxed control and their necessary
	conditions of optimality},
	journal = {J. Nonlin. Anal.},
	year = {2014},
	volume = {101},
	pages = {66-79},
}

@article{Ahm15,
	author = {Ahmed, N.U.},
	title = {Optimal control of general {M}c{K}ean-{V}lasov stochastic evolution equations on {H}ilbert spaces and necessary conditions of optimality},
	journal = {Math. Differ. Incl. Control Optim. Discussiones Mathematicae.},
	year = {2015},
	volume = {35},
	number = {2},
	pages = {165-195},
}

@article{Ahm16,
	author = {Ahmed, N.U.},
	title = {A general class of {M}c{K}ean-{V}lasov stochastic evolution equations driven by {B}rownian motion and {L}\`evy process and controlled by {L}\`evy measure},
	journal = {Math. Differ. Incl. Control Optim. Discussiones Mathematicae.},
	year = {2016},
	volume = {36},
	number = {2},
	pages = {181-206},
}

@article{BCS15,
	title = {p-variations of vector measures with respect to vector measures and integral representation 
	         of operators},
	journal = {Banach J. Math. Anal.},
	volume = {9},
	number = {1},
	pages = {273-285},
	year = {2015},
	author = {Blasco, O. and Calabuig, J.M. and Sanchez-Perez, E.A.},
}

@article{BDL11,
	title = {A General Stochastic Maximum Principle for SDEs of Mean-field Type},
	journal = {Appl. Math. Optim.},
	volume = {64},
	number = {2},
	pages = {197-216},
	year = {2011},
	author = {Buckdahn, R. and Djehiche, B. and Li, J.},
}

@article{BLM16,
	title = {A Stochastic Maximum Principle for General Mean-Field Systems},
	journal = {Appl. Math. Optim.},
	volume = {74},
	number = {2},
	pages = {507–534},
	year = {2016},
	author = {Buckdahn, R. and Li, J. and Ma, J.},
}

@article{OSD18,
	author = {Oksendal, B. and Sulem, A. and Dumitrescu, R.},
	title = {Stochastic control of general mean-field {SPDE}s with jumps},
	journal = {J. Optim. Theory Appl.},
	year = {2018},
	volume = {176},
	number = {3},
	pages = {559-584},
}

@article{DPJR19,
	author = {Da Prato, G. and Jentzen,A.  and Röckner, M.},
	title = {A mild Ito formula for SPDEs},
	journal = {T. Am. Math. Soc.
},
	year = {2019},
	volume = {372},
	pages = {755–3807},
}

@article{CGKPR22,
	author = {Cosso, A. and Gozzi, F. and Kharroubi, I. and Pham, H. and Rosestolato, M.},
	title = {Optimal control of path-dependent McKean-Vlasov SDEs in infinite dimension},
	journal = {Ann. Appl. Probab.},
    volume = {33},
    number = {4},
	year = {2023},
    pages = {2863-2918},
}

@article{WZ18,
	author = {Wu, C. and Zhang, J.},
	title = {An elementary proof for the structure of Wasserstein derivatives},
	howpublished = {preprint},
	journal = {arXiv:1705.08046},
	year = {2018},
	url={https://arxiv.org/pdf/1705.08046.pdf},
}

@Book{Din66,
	author = {Dinculeanu, N.},
	title = {Vector Measures},
	publisher = {VEB Deutscher Verlag der Wissenschaften},
	year = {1966},
	place = {Leipzig},
}

@Book{CD18,
	author = {Carmona, R. and Delarue, F.},
	title = {Probabilistic Theory of Mean Field Games with Applications {I}},
	publisher = {Springer},
	year = {2018},
	doi = {https://doi.org/10.1007/978-3-319-58920-6},
	place = {Berlin},
}

@article{JK10,
	author = {Jentzen, A. and Kloeden, P.},
	title = {Taylor Expansions Of Solutions Of Stochastic Partial Differential Equations With Additive Noise},
	journal = {Ann. Probab.},
    volume ={38},
    number={2},
    sites={532-569},
	year = {2010},
}

@article{SW21,
	author = {Stannat, W. and Wessels, L.},
	title = {Peng's Maximum Principle for Stochastic Partial Differential Equations},
	journal = {SIAM, J. Control and Optimization},
    volume ={59},
    number={5},
    sites={3552-3573},
	year = {2021},
}
\end{document}